\documentclass[reqno,11pt]{amsart}

\usepackage[top=1.9cm,bottom=1.9cm,left=2.4cm,right=2.4cm]{geometry}

\usepackage{amsmath,amsthm,amssymb}
\usepackage{mathrsfs,dsfont,amsfonts}
\usepackage{bbm}
\usepackage{functan}
\usepackage{extarrows}
\usepackage{mathtools}
\usepackage{appendix}
\usepackage{array,tabularx,multirow}

\usepackage[colorlinks,
            linkcolor=purple,
            anchorcolor=purple,
            citecolor=purple]{hyperref}

\usepackage{xcolor}
\usepackage{latexsym}
\usepackage{enumerate}
\usepackage{float}
\usepackage{relsize}
\usepackage{marginnote}
\usepackage{stmaryrd}
\usepackage{esint}
\usepackage{tikz}
\usepackage{bm}
\usepackage{verbatim}
\usepackage{cite}
\InputIfFileExists{agents/tmp/tikz_externalization.tex}{}{}
\providecommand{\tikzsetnextfilename}[1]{}
\allowdisplaybreaks

\usepackage{aliascnt}

\theoremstyle{plain}
\newtheorem{theorem}{Theorem}[section]

\newaliascnt{corollary}{theorem}

\aliascntresetthe{corollary}

\newaliascnt{lemma}{theorem}
\newtheorem{lemma}[lemma]{Lemma}
\aliascntresetthe{lemma}

\newaliascnt{proposition}{theorem}
\newtheorem{proposition}[proposition]{Proposition}
\aliascntresetthe{proposition}

\newaliascnt{remark}{theorem}
\newtheorem{remark}[remark]{Remark}
\aliascntresetthe{remark}

\theoremstyle{definition}
\newtheorem{definition}{Definition}[section]

\numberwithin{equation}{section}

\renewcommand{\arraystretch}{1.1}
\def\R{\mathbb{R}}

\def\Z{\mathbb{Z}}

\def\T{\mathbb{T}}
\def\d{\mathrm{d}}

\def\e{\mathrm{e}}
\def\i{\mathrm{i}}

\def\eps{\varepsilon}

\def\p{\partial}
\def\loc{\operatorname{loc}}

\def\div{\operatorname*{div}}

\def\ls{\lesssim}

\newcommand{\ch}{\mathbbm{1}}
\newcommand{\Nop}{\mathcal{N}_{\beta}}

\DeclareSymbolFont{curvevectors}{OML}{cmm}{b}{it}
\DeclareMathSymbol{\curvetangent}{\mathord}{curvevectors}{"1C}
\DeclareMathSymbol{\curvenormal}{\mathord}{curvevectors}{`n}
\renewcommand{\t}{\curvetangent}
\newcommand{\n}{\curvenormal}

\newcommand{\x}{\mathtt{x}}

\title[$H^2$ SQG Patches]
{On SQG Patches with $H^2$ Regularity}

\author[S.-Q. Chen]{Shaoqing Chen}
\address{Department of Mathematics, Nanjing University, Nanjing 210093, P. R. China}
\email{sqchen@smail.nju.edu.cn}
\author[Y.-Z. Sun]{Yongzhong Sun}
\address{Department of Mathematics, Nanjing University, Nanjing 210093, P. R. China}
\email{sunyz@nju.edu.cn}

\begin{document}

\begin{abstract}
In this paper, we study patch solutions of the generalized surface quasi-geostrophic (gSQG) equations, which is an interpolation between the 2D Euler ($\beta=0$) and the surface quasi-geostrophic (SQG) equations ($\beta=1$).
We establish in the SQG case, the local well-posedness for patches with $H^2$ boundaries, which improves the previous result of $H^{2+s}$, $s>0$ by Gancedo--Nguyen--Patel \cite{Gancedo2022}. 
Also, in the gSQG case with $0<\beta<1$, the uniqueness of $H^2$ patch solutions is established.
Our proofs are based on the contour dynamics formulation introduced in \cite{Gancedo2008}.
\end{abstract}



\maketitle

\section{Introduction}
In this paper, we study the generalized surface quasi-geostrophic (gSQG) equation 
\begin{equation}\label{eq:gSQG}\tag{gSQG}
\left\{
\begin{aligned}
&\p_t\theta+u\cdot\nabla\theta=0 &&\text{in }(0,T)\times\R^2,\\
&u=\nabla^\perp (-\Delta)^{-1+\frac{\beta}{2} }\theta,
 &&\text{in }(0,T)\times\R^2,\\
&\theta\big|_{t=0}=\theta_0
\end{aligned}
\right.
\end{equation}
for $0\le \beta\le 1$. 
\footnote{Here we denote $v^\perp=(-v_2,v_1)$ for a vector $v=(v_1,v_2)$.}
The endpoints $\beta=0$ and $\beta=1$ formally give the two-dimensional Euler and surface quasi-geostrophic (SQG) equations.
The SQG equation arises in geophysical fluid dynamics as a model for the evolution of surface potential temperature in a rapidly rotating, stratified fluid with uniform interior potential vorticity; see, for example, \cite{Held1995}.
Constantin--Majda--Tabak \cite{Constantin1994a} emphasized its analogy with the three-dimensional incompressible Euler equations and studied the formation of strong fronts. 
The velocity $u$ can be expressed by the Biot--Savart law as\footnote{We omit the constant dependent on $\beta$ here and in the rest of this paper for simplicity. This can be recovered by rescaling.}
\begin{equation}\label{eq:Biot-Savart}
u(t,\x)=\int_{\R^2}\frac{(\x-y)^\perp}{|\x-y|^{2+\beta}}\theta(t,y)\,\d y.
\end{equation}
In the classical theory, the Sobolev regularity of $\theta$ gives control of $\nabla u$ in $L^\infty$. 
For $0\leq\beta\leq1$, the velocity operator has order $\beta-1$, and standard Sobolev energy estimates yield local well-posedness of the inviscid equation for sufficiently regular initial data (see \cite[Section~2.4]{Constantin1994a} for SQG). 
When $1<\beta<2$, the positive order of the velocity operator leads to a loss of derivatives in the direct energy estimate. 
Chae--Constantin--C\'ordoba--Gancedo--Wu \cite{Chae2012} used a commutator estimate to address this issue and obtain local well-posedness in $H^m$ with $m\geq4$.
For the dissipative SQG equation, Constantin--Wu \cite{Constantin1999} proved global well-posedness for smooth data in the subcritical regime.
At the critical dissipation exponent, Kiselev--Nazarov--Volberg \cite{Kiselev2006} obtained the same result for smooth periodic data.

Our main objects of study are patch solutions with initial data of the form 
$$
\theta_0(y)=\ch_{\Omega_0}(y)=\begin{cases}
1,&y\in\Omega_0,\\
0,&y\notin\Omega_0,
\end{cases}
$$
where $\Omega_0\subset\R^2$ is a bounded simply connected domain.
In the SQG setting, such patches model sharp temperature fronts, which attracted early attention through studies of geophysical frontogenesis \cite{Constantin1994a} and the dynamics of temperature discontinuities \cite{Held1995}.
Since in \eqref{eq:gSQG} the active scalar $\theta$ is transported by a divergence-free velocity field $u$, 
the solution is expected to remain an indicator function, that is, of the form $\theta(t,y)=\ch_{\Omega(t)}(y)$, and $|\Omega(t)|$ remains constant in time due to incompressibility. 
However, since $\theta \in L^1 \cap L^\infty$ is discontinuous, the regularity of $u$ now becomes a crucial issue.
For the Euler equations ($\beta=0$), the velocity is log-Lipschitz, in which case the flow map is unique and the patch structure is preserved; 
for $0<\beta<1$, H\"older continuity $u\in C^{1-\beta}$ can be expected, while for $\beta=1$, we only have $u\in\text{BMO}$. 
The classical Lagrangian theory then does not apply here directly for $\beta>0$.

We begin by specifying the class of patch solutions we consider in the present work, as a weak solution of \eqref{eq:gSQG}, which is standard in the literature; see for example \cite{Gancedo2008,Gancedo2021,Gancedo2022}.
\begin{definition}\label{def:eulerian-patch}
In this paper, a \emph{patch solution} of \eqref{eq:gSQG} on $[0,T]$ with initial domain $\Omega_0$ is a family $\{\Omega(t)\}_{t\in[0,T]}$ such that:
\begin{enumerate}[(i)]
\item for each $t\in[0,T]$, $\Omega(t)$ is a bounded simply connected domain with a simple closed boundary, 
and
$\theta(t,y)=\ch_{\Omega(t)}(y)$, with $\Omega(0)=\Omega_0$;
\item with $u(t)=\nabla^\perp(-\Delta)^{-1+\frac{\beta}{2}}\ch_{\Omega(t)}$, one has, for every $\varphi\in C_c^\infty([0,T)\times\R^2)$,
\begin{equation}\label{eq:weak}
\int_0^T\!\int_{\Omega(t)} \bigl(\p_t\varphi+u\cdot\nabla\varphi\bigr)(t,y)\,\d y\,\d t+\int_{\Omega_0}\varphi(0,y)\,\d y=0;
\end{equation}
\item $\theta\in C_w([0,T];L^2(\R^2))$.
\end{enumerate}
Let $\T=\R/(2\pi\Z)$, with $|h|$ denoting the absolute value of the shortest representative of $h\in\T$.
For a regular simple closed curve $\gamma:\T\to\R^2$, define its arc--chord constant by
$$
F(\gamma):=\sup_{x\in\T,\,h\in\T\setminus\{0\}}
\frac{|h|}{|\gamma(x)-\gamma(x-h)|}.
$$
We call a patch solution an \emph{$H^2$ patch solution} if, for each
$t\in[0,T]$, its boundary admits a regular
normalized arclength parametrization $\gamma(t):\T\to\R^2$ satisfying
\begin{equation}\label{eq:uniform-H2-boundary-control}
\sup_{0\leq t\leq T}\left\|\gamma(t)\right\|_{H^2(\T;\R^2)}<\infty,
\qquad
\sup_{0\leq t\leq T}F(\gamma(t))<\infty.
\end{equation}
\end{definition}

A natural question is whether the patch structure is preserved, and whether the regularity of the boundary presists, under a rough velocity field? 
The study of such questions originated from the vortex patch solutions of the 2D Euler equations.
For the 2D Euler equations on the whole plane, 
Chemin \cite{Chemin1993} showed the global persistence of $C^{k,\alpha}\ (k\ge 1,\alpha\in(0,1))$ regularity of the patch boundary by paradifferential techniques, 
while the same result was obtained by Bertozzi--Constantin \cite{Bertozzi1993} using geometric arguments. 
Though Euler flows preserve the $C^{1,\alpha}$ regularity of the patch boundary, Kiselev--Luo \cite{Kiselev2023} constructed a $C^2$ patch whose boundary loses $C^2$ regularity instantaneously. See also \cite{Depauw1999,Kiselev2016,Kiselev2019} for vortex patches problem in bounded domains and in the half-plane.
The patch problems also naturally arise in various other models with transport nature.
Examples include density patches for the inhomogeneous incompressible Navier--Stokes equations \cite{Danchin2018,Hao2024c},
interfaces between fluids of different densities in the incompressible porous media equation (the Muskat problem) \cite{Cordoba2007,Cordoba2011a},
and temperature patches for the viscous Boussinesq system without thermal diffusion \cite{Danchin2016,Chae2024,Lazar2026}.

As we discussed, the discontinuity of $\theta$ prevents the direct application of well-posedness theory for classical solutions. 
The patch problem of \eqref{eq:gSQG} hence becomes also more challenging than the vortex patch problem for 2D Euler. 
Gancedo, in \cite{Gancedo2008}, proposed a new contour dynamics formulation of the patch, which ensures that the parametrization $\gamma$ has constant speed $|\gamma'|=c(t)$ for all time. 
Within this framework, the local well-posedness of $H^k$ patches for $k\ge 3$ was established in \cite{Gancedo2008,Cordoba2018,Gancedo2021} for patches in the whole plane and $0<\beta\le 1$. 
Later, Gancedo--Nguyen--Patel \cite{Gancedo2022} lowered the regularity threshold to $H^{2+s}$ for any $s>0$ in the SQG case $\beta=1$, allowing unbounded curvature.

The global well-posedness or possibility of finite time singularity of gSQG patches is still wide open. 
Bifurcation \cite{Castro2016a,Cordoba2019,Hmidi2025} and KAM methods \cite{GomezSerrano2026a,Hassainia2025a} have been applied to construct periodic and quasi-periodic solutions.
In domains with boundaries such as the half-plane, the interaction between the patch and the boundary may develop singularities in finite time.
Kiselev--Ryzhik--Yao--Zlato\v{s} \cite{Kiselev2016,Kiselev2017} proved finite time singularity formation for patches in the half-plane for $0<\beta<1/12$, which was later extended to $0<\beta<1/3$ by Gancedo--Patel \cite{Gancedo2021}. See also \cite{Miao2024b} for analysis of this scenario.
Jeon--Zlato\v{s} \cite{Jeon2025} established local well-posedness of $H^2$ multi-patch solutions in the whole plane for $0<\beta\le 1/3$, allowing patches to touch at a point, and constructed solutions with smooth initial data that develop singularities in finite time.
Recently, G\'omez-Serrano--Pasquariello \cite{GomezSerrano2026} constructed sharp front solutions of gSQG with $7 / 4 <\beta<2$ that develop singularities in finite time.
See, for example, \cite{Rodrigo2005,Ai2026} for further work on the well-posedness of gSQG sharp fronts.

In the present paper, we focus on the local well-posedness of gSQG patches.
We summarize in \autoref{tab:local-patch-theory} the results regarding local existence, uniqueness, and ill-posedness relevant to this problem.
\begin{table}[tbp]
\centering
\caption{Local existence, uniqueness, and ill-posedness for gSQG patches.}
\label{tab:local-patch-theory}
\footnotesize
\setlength{\tabcolsep}{4pt}
\renewcommand{\arraystretch}{1.2}
\renewcommand{\tabularxcolumn}[1]{m{#1}}
\begin{tabularx}{\textwidth}{|>{\centering\arraybackslash}m{0.28\textwidth}
|>{\centering\arraybackslash}m{0.16\textwidth}
|>{\centering\arraybackslash}m{0.26\textwidth}
|>{\centering\arraybackslash}X|}
\hline
Reference & Range & Boundary regularity & Result\\
\hline
\multirow{2}{=}{\centering Gancedo \cite{Gancedo2008}}
& $0<\beta<1$
& $H^k$, $k\geq3$
& Local existence and uniqueness\\
\cline{2-4}
& $\beta=1$
& $H^k$, $k\geq3$
& Local existence\\
\hline
Chae--Constantin--C\'ordoba--Gancedo--Wu \cite{Chae2012}
& $1<\beta<2$
& $H^k$, $k\geq4$
& Local existence\\
\hline
Kiselev--Yao--Zlato\v{s} \cite{Kiselev2017}
& $0<\beta<1$
& $H^3$
& Local existence and uniqueness\\
\hline
C\'ordoba--C\'ordoba--Gancedo \cite{Cordoba2018}
& $\beta=1$
& $C^{2,\delta}$, $0<\delta<1$
& Uniqueness\\
\hline
\multirow{2}{=}{\centering Gancedo--Patel \cite{Gancedo2021}}
& $0<\beta<1$
& $H^2$
& Local existence\\
\cline{2-4}
& $1<\beta<2$
& $H^3$
& Local existence and uniqueness\\
\hline
Gancedo--Nguyen--Patel \cite{Gancedo2022}
& $\beta=1$
& $H^{2+s}$, $s>0$
& Local existence and uniqueness\\
\hline
\multirow{2}{=}{\centering Luo \cite{Luo2024c}}
& $0<\beta<1$
& $W^{2,p}$, $p>1/(1-\beta)$
& Uniqueness\\
\cline{2-4}
& $0<\beta<1/2$
& $C^{1,\delta}$, $2\beta<\delta<1$
& Uniqueness\\
\hline
\multirow{2}{=}{\centering Kiselev--Luo \cite{Kiselev2025}}
& \multirow{2}{=}{\centering $0<\beta<1$}
& $C^{2,\delta}$, $0\leq\delta<1$
& Strong ill-posedness\\
\cline{3-4}
&
& $W^{2,p}$, $p>1/(1-\beta)$, $p\ne2$
& Strong ill-posedness\\
\hline
Jeon--Zlato\v{s} \cite{Jeon2025}
& $0<\beta\leq1/3$
& $H^2$
& Local existence and uniqueness; touching multi-patches allowed\\
\hline
\multirow{2}{=}{\centering Maga\~na \cite{Magana2026}}
& \multirow{2}{=}{\centering $K(r)\propto r^{-\beta}$, $0<\beta<1$%
\footnotemark}
& $H^3$
& Local existence and uniqueness\\
\cline{3-4}
&
& $H^2$
& Local existence\\
\hline
\end{tabularx}
\end{table}
\footnotetext[\value{footnote}]{With several technical hypotheses, see\cite{Magana2026} for details.}
Our first main result establishes uniqueness of $H^2$ patch solutions for $0<\beta<1$. We will prove the following theorem in \autoref{sec:gsqg}.

\begin{theorem}[Uniqueness for $H^2$ gSQG patches]\label{thm:gsqg-uniqueness}
For $0<\beta<1$, suppose $\Omega_1(t)$ and $\Omega_2(t)$ are $H^2$
patch solutions of \eqref{eq:gSQG} on $[0,T]$, in the sense of
\autoref{def:eulerian-patch}. If
$\Omega_1(0)=\Omega_2(0)$, then
$
\Omega_1(t)=\Omega_2(t),\forall\, 0\leq t\leq T.
$
\end{theorem}

With estimates similar to the proof of \autoref{thm:gsqg-uniqueness}, we are able to establish the local well-posedness of SQG patches with $H^2$ boundaries. The proof will be presented in \autoref{sec:endpoint}.

\begin{theorem}[Local well-posedness for $H^2$ SQG patches]
\label{thm:SQG-wellposedness}
Let $\Omega_0$ be a bounded simply connected domain with $H^2$ boundary.
For initial datum $\theta_0=\ch_{\Omega_0}$, there is $T>0$ such that, \eqref{eq:gSQG} with $\beta=1$ admits a unique $H^2$ patch solution on $[0,T]$ in the sense of \autoref{def:eulerian-patch}, 
whose boundary admits a regular parametrization in $C([0,T];H^2)$.
\end{theorem}

The rest of the paper is organized as follows. 
In \autoref{sec:contour-dynamics}, we introduce the contour dynamics formulation and the constant-speed tangential gauge, which is the key tool for our analysis. 
In \autoref{sec:gsqg}, we prove the uniqueness \autoref{thm:gsqg-uniqueness} for $0<\beta<1$. 
In \autoref{sec:endpoint}, we prove \autoref{thm:SQG-wellposedness} for $\beta=1$. 
Finally, the rigorous derivation of the contour dynamics formulation and the constant-speed tangential gauge is given in \autoref{sec:patch-contour}.

\section{Contour dynamics and tangential gauges}\label{sec:contour-dynamics}

\subsection{The contour equation}

Before we proceed, we need to introduce some basic notations.
Denote the difference of a (vector-valued) function $f$ on $\T$ by
$$
\Delta_h f(x)=f(x)-f(x-h),\qquad h\in\T.
$$
In the following, we mainly deal with functions defined on $[0,T]\times\T$. To avoid too many subscripts, we write their derivatives by
$$
\gamma'(t,x)=\p_x\gamma(t,x),\qquad \gamma_t(t,x)=\p_t\gamma(t,x)
$$
and so on.
For a positively oriented regular curve $\gamma$, we denote its outward unit normal as $\n_\gamma=-(\gamma')^\perp/|\gamma'|$.

We first derive the contour equation formally. Let $\gamma(t,\cdot)$ be a positively oriented parametrization of $\p\Omega(t)$, and $n(t,y)$ be its outward unit normal.
The Biot--Savart law for gSQG \eqref{eq:Biot-Savart} and Green's formula give, for any $z\notin\p\Omega(t)$,
$$
\begin{aligned}
u(t,z)=\nabla_z^\perp\int_{\Omega(t)}\frac{1}{|z-y|^\beta}\,\d y
=-\int_{\p\Omega(t)}\frac{n(t,y)^\perp}{|z-y|^\beta}\,\d\sigma(y)
=-\int_\T\frac{\gamma'(t,x)}{|z-\gamma(t,x)|^\beta}\,\d x,
\end{aligned}
$$
where in the last equality we note that $n(t,\gamma(t,x))^\perp|\gamma'(t,x)|=\gamma'(t,x)$.
Writing $\n_\gamma(t,x)=n(t,\gamma(t,x))$, we therefore obtain by noting $\gamma'(t,x)\cdot\n_\gamma(t,x)=0$, that
$$
\begin{aligned}
\gamma_t(t,x)\cdot\n_\gamma(t,x)
&=u(t,\gamma(t,x))\cdot\n_\gamma(t,x)\\
&=-\int_\T \frac{\gamma'(t,x-h)\cdot\n_\gamma(t,x)}{|\Delta_h\gamma(t,x)|^\beta}\,\d h\\
&=\int_\T\frac{\Delta_h\gamma'(t,x)\cdot\n_\gamma(t,x)}{|\Delta_h\gamma(t,x)|^\beta}\,\d h=\Nop(\gamma)(t,x)\cdot\n_\gamma(t,x).
\end{aligned}
$$
Here the contour operator is defined as
$$
\Nop(\gamma)(x):=\int_{\T} \frac{\Delta_h\gamma'(x)}{|\Delta_h\gamma(x)|^\beta}\,\d h.
$$
Now $\gamma'(t,x)\cdot\n_\gamma(t,x)=0$ allows us to replace
$-\gamma'(t,x-h)$ by $\Delta_h\gamma'(t,x)$ in the normal component.
Thus, $\gamma_t-\Nop(\gamma)\parallel\gamma'$, and we can rewrite the contour equation as
$$
\gamma_t=\Nop(\gamma)+\lambda\gamma'.
$$
Here $\lambda$ is a scalar function, which can be chosen in a quite flexible way to specify the tangential motion of the parametrization.

\subsection{The constant-speed gauge}\label{sec:constant-speed-gauge}

We choose the tangential coefficient $\lambda_\gamma$ so that $|\gamma'(t,x)|$ is independent of $x$ and $\lambda_\gamma(t,-\pi)=0$. The following calculation is classical and standard since \cite[Sec.~5]{Gancedo2008}, see also \cite{Cordoba2018,Gancedo2022}.
Writing the contour equation as
\begin{equation}\label{eq:contour-equation}
\gamma_t=\Nop(\gamma)+\lambda_\gamma\gamma',
\end{equation}
we have
$$
\frac12\p_t|\gamma'|^2=\gamma'\cdot(\gamma_t)'
=\gamma'\cdot\bigl(\Nop(\gamma)\bigr)'+\lambda_\gamma'|\gamma'|^2+\lambda_\gamma\gamma'\cdot\gamma''.
$$
The last term vanishes because $\gamma'\cdot\gamma''=0$. This implies
$$
\frac{\p_t|\gamma'|^2}{2|\gamma'|^2}=\frac{\gamma'\cdot\bigl(\Nop(\gamma)\bigr)'}{|\gamma'|^2}+\lambda_\gamma'.
$$
Now we denote the contribution of the contour operator by
\begin{equation}\label{eq:g-definition}
g_\gamma:=\frac{\gamma'\cdot\bigl(\Nop(\gamma)\bigr)'}{|\gamma'|^2}.
\end{equation}
Since $\lambda_\gamma$ is periodic, averaging over $\T$ and then comparing with the previous identity, we obtain
$$
\frac{\p_t|\gamma'|^2}{2|\gamma'|^2}=\frac1{2\pi}\int_\T g_\gamma(y)\,\d y = \lambda_\gamma' + g_\gamma.
$$
We assume, for convenience, that $\lambda_\gamma(t,-\pi)=0$. Finally we can write the tangential coefficient as
\begin{equation}\label{eq:lambda-definition}
\lambda_\gamma(x)=\frac{x+\pi}{2\pi}\int_\T g_\gamma(y)\,\d y -\int_{-\pi}^x g_\gamma(y)\,\d y.
\end{equation}
The following proposition states the correspondence between patch solutions and contour solutions in this gauge. Its proof is given in \autoref{sec:patch-contour}.

\begin{proposition}\label{prop:patch-constant-speed-equivalence}
For $0<\beta<1$, the $H^2$ patch formulation in
\autoref{def:eulerian-patch} is equivalent to the constant-speed contour equation
\eqref{eq:contour-equation}-\eqref{eq:lambda-definition} in the following sense:
Every such patch solution $\{ \Omega(t) \}_{t\in[0,T]}$ admits an $H^2$ parametrization $\gamma(t):\T\to\R^2$ of $\p\Omega(t)$, which is continuous in time and satisfies \eqref{eq:contour-equation}-\eqref{eq:lambda-definition}.
Conversely, if a solution $\gamma(t)$ of \eqref{eq:contour-equation}-\eqref{eq:lambda-definition} has uniform-in-time $H^2$ and arc--chord bounds, then the enclosed domains $\Omega(t)$ form an $H^2$ patch solution of \eqref{eq:gSQG}.

For $\beta=1$, the correspondence holds for patches whose boundaries admit a regular parametrization in $C([0,T];H^2)$, and for contour solutions in $C([0,T];H^2)$.
\end{proposition}

\section{The gSQG Case $0<\beta<1$: Uniqueness}\label{sec:gsqg}
We recall the existence established by Gancedo and Patel \cite[Theorem~3]{Gancedo2021}: For $0<\beta<1$, constant-speed $H^2$ curves $\gamma,\zeta$ initiated from given $H^2$ curves $\gamma_0,\zeta_0$
exist on $[0,T]$, solving the contour equation \eqref{eq:contour-equation}-\eqref{eq:lambda-definition},
with uniform $H^2$ and arc--chord bounds
\begin{equation}\label{eq:gsqg_bound}
\sup_{t \in [0,T]}(\left\|\gamma(t) \right\|_{H^2} + \left\|\zeta(t) \right\|_{H^2}) < +\infty ,\quad F:=\sup_{t \in [0,T]}\max\{1, F(\gamma(t)), F(\zeta(t))\} < +\infty.
\end{equation}
This section is devoted comparing these solutions, and to prove \autoref{thm:gsqg-uniqueness}.
We denote $|X|\leq CY$ by $X\ls Y$ for a constant $C>0$ that may depend on $\beta$ in this section.

\subsection{Estimates for constant-speed curves}\label{sec:geometric-identities}
Before we proceed to the proof, we collect some basic estimates for constant-speed curves, 
which are well-known in the literature for example \cite{Gancedo2008}. 
First, constant speed gives
\begin{equation}\label{eq:speed-arc-chord}
|\gamma'|=(2\pi)^{-1/2}\left\|\gamma'\right\|_{L^2},\qquad
|\zeta'|=(2\pi)^{-1/2}\left\|\zeta'\right\|_{L^2},\qquad
|\gamma'|^{-1}\leq F(\gamma),\qquad |\zeta'|^{-1}\leq F(\zeta).
\end{equation}
Indeed,
\begin{equation}\label{eq:arc-chord-lower}
|\Delta_h\Gamma(x)|^{-q}\leq \frac{F(\Gamma)^q}{|h|^q} \qquad(q>0);
\end{equation}
For $f\in H^1(\T)$ and $0<h\leq\pi$, we have
\begin{align}
\left\|\Delta_h f\right\|_{L^2} &\leq\min\{2\left\|f\right\|_{L^2},\left\|f'\right\|_{L^2} h\}, \label{eq:basic-difference-L2}\\
|\Delta_h f(x)| &\leq h^{1/2}\left(\int_{x-h}^{x}|f'(r)|^2\,\d r\right)^{1/2}. \label{eq:basic-difference-local}
\end{align}
These inequalities follow simply from the fundamental theorem of calculus and Cauchy--Schwarz.
For $h<0$, the same bounds hold with $|h|$ and the interval between $x$ and $x-h$.

We also note that for constant-speed $\gamma,\zeta$, the reverse triangle inequality gives (note that the left-hand side is constant in $x$)
\begin{equation}\label{eq:speed-difference}
\bigl||\gamma'|-|\zeta'|\bigr| \leq (2\pi)^{-1/2}\left\|\gamma'-\zeta'\right\|_{L^2}.
\end{equation}

\begin{lemma}\label{lem:chord-identities}
For every constant-speed $H^2$ curve $\Gamma$,
\begin{align}
&\Delta_h\Gamma(x)\cdot\Delta_h\Gamma'(x)=\frac12\int_{x-h}^{x}\left(|\Gamma'(r)-\Gamma'(x-h)|^2 -|\Gamma'(r)-\Gamma'(x)|^2\right)\,\d r, 
\label{eq:chord-tangent-identity}\\
&|\Gamma'|^2h^2-|\Delta_h\Gamma(x)|^2 =\frac12\int_{x-h}^{x}\int_{x-h}^{x} |\Gamma'(r)-\Gamma'(s)|^2\,\d r\,\d s.
\label{eq:chord-defect-identity}
\end{align}
Consequently, 
\begin{align}
|\Delta_h\Gamma(x)\cdot\Delta_h\Gamma'(x)| &\ls |h|^2\left|\int_{x-h}^{x}|\Gamma''(r)|^2\,\d r\right|.
\label{eq:CV-bound}
\end{align}
\end{lemma}

\begin{proof}
It is enough to consider $h>0$. Since $\Delta_h\Gamma(x)=\int_{x-h}^{x}\Gamma'(r)\,\d r$,
$$
\begin{aligned}
2\Delta_h\Gamma(x)\cdot\Delta_h\Gamma'(x)
&=\int_{x-h}^{x} 2\Gamma'(r)\cdot\left(\Gamma'(x)-\Gamma'(x-h)\right)\,\d r\\
&=\int_{x-h}^{x} \left(|\Gamma'(r)-\Gamma'(x-h)|^2 -|\Gamma'(r)-\Gamma'(x)|^2\right)\,\d r.
\end{aligned}
$$
Here we used $|\Gamma'(x)|=|\Gamma'(x-h)|$. Similarly,
$$
\begin{aligned}
|\Delta_h\Gamma(x)|^2
&=\int_{x-h}^{x}\int_{x-h}^{x} \Gamma'(r)\cdot\Gamma'(s)\,\d r\,\d s\\
&=|\Gamma'|^2h^2-\frac12\int_{x-h}^{x}\int_{x-h}^{x} |\Gamma'(r)-\Gamma'(s)|^2\,\d r\,\d s.
\end{aligned}
$$
The inequality \eqref{eq:CV-bound} follows from the fundamental theorem of calculus and Cauchy--Schwarz inequalities.
\end{proof}

\begin{lemma}\label{lem:chord-square-difference}
Let $\gamma,\zeta\in H^2(\T;\R^2)$ have positive constant speeds.
For $0<|h|\leq\pi$,
\begin{equation}\label{eq:chord-defect-difference-bound}
\left| \bigl(|\gamma'|^2h^2-|\Delta_h\gamma|^2\bigr) -\bigl(|\zeta'|^2h^2-|\Delta_h\zeta|^2\bigr) \right|
\ls |h|^2 \left(\left\|\gamma''\right\|_{L^2}+\left\|\zeta''\right\|_{L^2}\right)
\left\|\gamma'-\zeta'\right\|_{L^2}.
\end{equation}
Consequently,
\begin{equation}\label{eq:normalized-chord-difference}
\left| |\Delta_h\gamma|^2-|\Delta_h\zeta|^2 \right|
\ls |h|^2\left(\left\|\gamma\right\|_{H^2}+\left\|\zeta\right\|_{H^2}\right) \left\|\gamma'-\zeta'\right\|_{L^2}.
\end{equation}
\end{lemma}

\begin{proof}
It is enough to consider $h>0$. By \eqref{eq:chord-defect-identity},
$$
\begin{aligned}
&\bigl(|\gamma'|^2h^2-|\Delta_h\gamma|^2\bigr) -\bigl(|\zeta'|^2h^2-|\Delta_h\zeta|^2\bigr)\\
&=\frac12\int_{x-h}^{x}\int_{x-h}^{x} \bigl[\gamma'(r)-\gamma'(s)+\zeta'(r)-\zeta'(s)\bigr]\cdot \bigl[(\gamma'-\zeta')(r)-(\gamma'-\zeta')(s)\bigr]\,\d r\,\d s.
\end{aligned}
$$
Since $|\gamma'(r)-\gamma'(s)|\leq h^{1/2}\left\|\gamma''\right\|_{L^2}$ for $r,s\in[x-h,x]$, and the same bound holds for $\zeta$, we obtain
$$
\begin{aligned}
&\left|
\bigl(|\gamma'|^2h^2-|\Delta_h\gamma|^2\bigr) -\bigl(|\zeta'|^2h^2-|\Delta_h\zeta|^2\bigr)\right|\\
\leq&\ \frac{h^{1/2}}2 \left(\left\|\gamma''\right\|_{L^2}+\left\|\zeta''\right\|_{L^2}\right)
\int_{x-h}^{x}\int_{x-h}^{x} \bigl(|\gamma'(r)-\zeta'(r)|+|\gamma'(s)-\zeta'(s)|\bigr)\,\d r\,\d s\\
=&\ h^{3/2} \left(\left\|\gamma''\right\|_{L^2}+\left\|\zeta''\right\|_{L^2}\right) \int_{x-h}^{x}|\gamma'(r)-\zeta'(r)|\,\d r\\
\leq&\ h^2 \left(\left\|\gamma''\right\|_{L^2}+\left\|\zeta''\right\|_{L^2}\right) \left\|\gamma'-\zeta'\right\|_{L^2},
\end{aligned}
$$
which proves \eqref{eq:chord-defect-difference-bound}. Here we used Cauchy--Schwarz in the last inequality.

For \eqref{eq:normalized-chord-difference}, note that by \eqref{eq:speed-difference},
$$
\left||\gamma'|^2 - |\zeta'|^2\right| \le (|\gamma'|+|\zeta'|)\bigl||\gamma'|-|\zeta'|\bigr| \ls (|\gamma'|+|\zeta'|) \left\|\gamma'-\zeta'\right\|_{L^2},
$$
we have
$$
\begin{aligned}
\left||\Delta_h\gamma|^2-|\Delta_h\zeta|^2\right|
&\leq h^2\left||\gamma'|^2 - |\zeta'|^2\right| + \left| \bigl(|\gamma'|^2h^2-|\Delta_h\gamma|^2\bigr) -\bigl(|\zeta'|^2h^2-|\Delta_h\zeta|^2\bigr)\right|\\
&\ls h^2\left(\left\|\gamma\right\|_{H^2}+\left\|\zeta\right\|_{H^2}\right) \left\|\gamma'-\zeta'\right\|_{L^2}.
\end{aligned}
$$
Here we used \eqref{eq:speed-arc-chord} and \eqref{eq:chord-defect-difference-bound} in the last inequality.
\end{proof}

We now estimate the difference of kernel values.
\begin{lemma}\label{lem:gsqg-flat-chord}
Let $\gamma,\zeta$ be constant-speed $H^2$ curves with positive speeds and arc--chord quantities bounded by $F$.  
For $\beta>0$ and $0<|h|\leq\pi$,
\begin{equation}\label{eq:normalized-kernel-difference}
\left||\Delta_h\gamma|^{-\beta}-|\Delta_h\zeta|^{-\beta}\right|
\ls F^{\beta+2}|h|^{-\beta} \left(\left\|\gamma\right\|_{H^2}+\left\|\zeta\right\|_{H^2}\right) \left\|\gamma'-\zeta'\right\|_{L^2}.
\end{equation}
\end{lemma}

\begin{proof}
The mean value theorem applied to $z^{-\beta/2}$ gives
$$
\begin{aligned}
&\left||\Delta_h\gamma|^{-\beta}-|\Delta_h\zeta|^{-\beta}\right|\\
\le\ &\frac\beta2 F^{\beta+2}|h|^{-\beta-2} \left||\Delta_h\gamma|^2-|\Delta_h\zeta|^2\right|\\
\ls\ & F^{\beta+2}|h|^{-\beta} \left(\left\|\gamma\right\|_{H^2}+\left\|\zeta\right\|_{H^2}\right) \left\|\gamma'-\zeta'\right\|_{L^2}.
\end{aligned}
$$
Here the derivative is bounded by $(\beta/2)F^{\beta+2}|h|^{-\beta-2}$ since by \eqref{eq:arc-chord-lower}, $|\Delta_h\gamma|^2\geq F^{-2}|h|^2$ and $|\Delta_h\zeta|^2\geq F^{-2}|h|^2$. 
The last inequality follows from \eqref{eq:normalized-chord-difference}.
\end{proof}

\begin{lemma}\label{lem:chord-tangent-difference}
Let $\gamma,\zeta$ be constant-speed $H^2$ curves with positive speeds. 
For $0<|h|\leq\pi$,
\begin{equation}\label{eq:H-relative}
\left\|\Delta_h\gamma\cdot\Delta_h\gamma' -\Delta_h\zeta\cdot\Delta_h\zeta'\right\|_{L^2}
\ls |h|^{3/2} \left(\left\|\gamma\right\|_{H^2}+\left\|\zeta\right\|_{H^2}\right) \left\|\gamma'-\zeta'\right\|_{L^2}.
\end{equation}
\end{lemma}

\begin{proof}
It is enough to consider $h>0$.  For $0<r\leq h$, H\"older's inequality gives
\begin{equation}\label{eq:tangent-square-difference}
\begin{aligned}
\left\||\Delta_r\gamma'|^2-|\Delta_r\zeta'|^2\right\|_{L^2}
&=\left\|(\Delta_r\gamma'+\Delta_r\zeta')\cdot \Delta_r(\gamma'-\zeta')\right\|_{L^2}\\
&\leq\left(\left\|\Delta_r\gamma'\right\|_{L^\infty} +\left\|\Delta_r\zeta'\right\|_{L^\infty}\right) \left\|\Delta_r(\gamma'-\zeta')\right\|_{L^2}\\
&\leq2r^{1/2} \left(\left\|\gamma''\right\|_{L^2}+\left\|\zeta''\right\|_{L^2}\right) \left\|\gamma'-\zeta'\right\|_{L^2}.
\end{aligned}
\end{equation}
Here we used 
$\left\|\Delta_r\gamma'\right\|_{L^\infty}\leq r^{1/2}\left\|\gamma''\right\|_{L^2}$
and $\left\|\Delta_r(\gamma'-\zeta')\right\|_{L^2} \leq2\left\|\gamma'-\zeta'\right\|_{L^2}$.
On the right hand side of \eqref{eq:chord-tangent-identity}, we put $s=x-h+r$ in the first term and $s=x-r$ in the second. This gives
$$
\begin{aligned}
\Delta_h\gamma(x)\cdot\Delta_h\gamma'(x)&=\frac12\int_{x-h}^x\left(|\gamma'(s)-\gamma'(x-h)|^2-|\gamma'(s)-\gamma'(x)|^2\right)\,\d s\\
&=\frac12\int_0^h\left(|\Delta_{-r}\gamma'(x-h)|^2-|\Delta_r\gamma'(x)|^2\right)\,\d r.
\end{aligned}
$$
Since $\Delta_{-r}f(x)=-\Delta_rf(x+r)$, we have by applying Minkowski's inequality that
$$
\begin{aligned}
&\left\|\Delta_h\gamma\cdot\Delta_h\gamma'-\Delta_h\zeta\cdot\Delta_h\zeta'\right\|_{L^2}\\
\leq&\ \frac12\int_0^h\left\||\Delta_{-r}\gamma'(\cdot-h)|^2-|\Delta_{-r}\zeta'(\cdot-h)|^2\right\|_{L^2}\,\d r+\frac12\int_0^h\left\||\Delta_r\gamma'|^2-|\Delta_r\zeta'|^2\right\|_{L^2}\,\d r\\
=&\ \int_0^h\left\||\Delta_r\gamma'|^2-|\Delta_r\zeta'|^2\right\|_{L^2}\,\d r\\
\ls&\ \left(\left\|\gamma''\right\|_{L^2}+\left\|\zeta''\right\|_{L^2}\right)\left\|\gamma'-\zeta'\right\|_{L^2}\int_0^h r^{1/2}\,\d r\\
\ls&\ h^{3/2}\left(\left\|\gamma\right\|_{H^2}+\left\|\zeta\right\|_{H^2}\right)\left\|\gamma'-\zeta'\right\|_{L^2}.
\end{aligned}
$$
This proves \eqref{eq:H-relative}. 
\end{proof}

\subsection{Estimates for gauge terms}
Now we establish the estimates for $g_\gamma$ and $\lambda_\gamma$.

\begin{lemma}\label{lem:single-gauge}
Let $\gamma\in H^2$ have positive constant speed and satisfy the arc--chord condition. Then
\begin{equation}\label{eq:gauge-exact}
g_\gamma(x)=-\frac{\beta}{2|\gamma'|^2}\int_\T \frac{\Delta_h\gamma(x)\cdot \gamma'(x)}{|\Delta_h\gamma(x)|^{\beta+2}} |\Delta_h\gamma'(x)|^2\,\d h.
\end{equation}
In particular,
\begin{equation}\label{eq:g-single-F}
\left\|g_\gamma\right\|_{L^\infty}+\left\|\lambda_\gamma'\right\|_{L^\infty}+\left\|\lambda_\gamma\right\|_{L^\infty}
\ls F(\gamma)^{\beta+2} \left\|\gamma''\right\|_{L^2}^2
\end{equation}
holds when $0<\beta<1$.
\end{lemma}

\begin{proof}
Write the kernel with independent parameters $x$ and $y$.
The differentiated numerator has zero scalar product with $\gamma'(x)$, since $\gamma'\cdot\gamma''=0$. Hence
$$
\begin{aligned}
\gamma'(x)\cdot (\mathcal{N}_\beta(\gamma))'(x)
&=\gamma'(x)\cdot \left( \int_\T \frac{\gamma'(x)-\gamma'(y)}{|\gamma(x)-\gamma(y)|^\beta}\,\d y \right)'\\
&=-\beta\int_\T \frac{(\gamma(x)-\gamma(y))\cdot\gamma'(x)} {|\gamma(x)-\gamma(y)|^{\beta+2}} (\gamma'(x)-\gamma'(y))\cdot\gamma'(x)\,\d y,
\end{aligned}
$$
hence
$$
\begin{aligned}
g_\gamma(x) &=(\mathcal{N}_\beta(\gamma))'(x) \cdot \frac{\gamma'(x)}{|\gamma'|^2} \\
&=-\frac{\beta}{|\gamma'|^2} \int_\T \frac{(\gamma(x) - \gamma(y)) \cdot \gamma'(x)}{|\gamma(x) - \gamma(y)|^{\beta + 2}} (\gamma'(x) - \gamma'(y))\cdot \gamma'(x) \, \d y \\
&=-\frac{\beta}{2|\gamma'|^2} \int_\T \frac{(\gamma(x) - \gamma(y)) \cdot \gamma'(x)}{|\gamma(x) - \gamma(y)|^{\beta + 2}} |\gamma'(x) - \gamma'(y)|^2 \, \d y,
\end{aligned}
$$
which proves \eqref{eq:gauge-exact}. Here we used the constant-speed condition in the last equality.

The $L^\infty$ bound for $g_\gamma$ follows from the following estimates: For $0<|h|\leq\pi$,
$$
|\Delta_h\gamma\cdot \gamma'(x)|\leq |\gamma'|^2|h|,
\qquad |\Delta_h\gamma'|^2 \leq |h|\int_{x-h}^{x}|\gamma''(r)|^2\,\d r,
$$
and \eqref{eq:arc-chord-lower}.  Thus the integrand in \eqref{eq:gauge-exact} satisfies
\begin{equation}\label{eq:gauge-integrand-local}
\begin{aligned}
\left|
-\frac{\beta}{2|\gamma'|^2}
\frac{\Delta_h\gamma(x)\cdot\gamma'(x)}
{|\Delta_h\gamma(x)|^{\beta+2}}
|\Delta_h\gamma'(x)|^2
\right|
&\ls F(\gamma)^{\beta+2}|h|^{-\beta}
\int_{x-h}^{x}|\gamma''(r)|^2\,\d r\\
&\leq F(\gamma)^{\beta+2}\left\|\gamma''\right\|_{L^2}^2|h|^{-\beta}.
\end{aligned}
\end{equation}
For $h<0$, these local integrals are taken over $[x,x-h]$.
Integrating over $h$ gives the desired bound. The bounds
for $\lambda_\gamma'$ and $\lambda_\gamma$ follow from \eqref{eq:lambda-definition}.
\end{proof}

\begin{proposition}\label{prop:gauge-stability}
Let $\gamma,\zeta$ be as in \autoref{lem:gsqg-flat-chord}.  Then
\begin{align}
\left\|g_\gamma-g_\zeta\right\|_{L^1}
&\ls F^{\beta+5}
\left(1+\left\|\gamma\right\|_{H^2}+\left\|\zeta\right\|_{H^2}\right)^4
\left\|\gamma'-\zeta'\right\|_{L^2},
\label{eq:gsqg-g-stability}\\
\left\|\lambda_\gamma-\lambda_\zeta\right\|_{L^\infty}
&\ls F^{\beta+5}
\left(1+\left\|\gamma\right\|_{H^2}+\left\|\zeta\right\|_{H^2}\right)^4
\left\|\gamma'-\zeta'\right\|_{L^2}.
\label{eq:gsqg-lambda-stability}
\end{align}
\end{proposition}

\begin{proof}
We assume $0<h\leq\pi$ without loss of generality. 
By \eqref{eq:basic-difference-L2} and
\eqref{eq:basic-difference-local}, we have
$$
\begin{aligned}
\left\|\Delta_h(\gamma-\zeta)\right\|_{L^\infty}
&\leq h^{1/2}\left\|\gamma'-\zeta'\right\|_{L^2},\\
\left\|\Delta_h\gamma'\right\|_{L^2}&\leq h\left\|\gamma''\right\|_{L^2}.
\end{aligned}
$$

By \eqref{eq:gauge-exact},
\begin{equation}\label{eq:gauge-difference}
\begin{aligned}
g_\gamma-g_\zeta
&=-\frac{\beta}{2|\gamma'|^2}\int_\T
\left[
\frac{\Delta_h\gamma\cdot\gamma'}{|\Delta_h\gamma|^{\beta+2}}
|\Delta_h\gamma'|^2
-\frac{\Delta_h\zeta\cdot\zeta'}{|\Delta_h\zeta|^{\beta+2}}
|\Delta_h\zeta'|^2
\right]\,\d h\\
&\quad-\frac{\beta}{2}
\left(|\gamma'|^{-2}-|\zeta'|^{-2}\right)
\int_\T
\frac{\Delta_h\zeta\cdot\zeta'}{|\Delta_h\zeta|^{\beta+2}}
|\Delta_h\zeta'|^2\,\d h.
\end{aligned}
\end{equation}

To deal with the first term, for a fixed $h$, 
we split the difference of the integrands in \eqref{eq:gauge-exact} as
\begin{equation}\label{eq:gauge-difference-1-split}
\begin{aligned}
&\frac{\Delta_h\gamma\cdot \gamma'}{|\Delta_h\gamma|^{\beta+2}}
|\Delta_h\gamma'|^2
-\frac{\Delta_h\zeta\cdot \zeta'}{|\Delta_h\zeta|^{\beta+2}}
|\Delta_h\zeta'|^2 = K_1 + K_2 + K_3 + K_4, \\
& K_1 :=\left(|\Delta_h\gamma|^{-\beta-2}-|\Delta_h\zeta|^{-\beta-2}\right)
(\Delta_h\gamma\cdot \gamma')|\Delta_h\gamma'|^2,\\
& K_2 :=|\Delta_h\zeta|^{-\beta-2}
\bigl(\Delta_h(\gamma-\zeta)\cdot\gamma'\bigr)
|\Delta_h\gamma'|^2,\\
& K_3 :=|\Delta_h\zeta|^{-\beta-2}
\bigl(\Delta_h\zeta\cdot(\gamma'-\zeta')\bigr)
|\Delta_h\gamma'|^2,\\
& K_4 :=|\Delta_h\zeta|^{-\beta-2}(\Delta_h\zeta\cdot \zeta')
\left(|\Delta_h\gamma'|^2-|\Delta_h\zeta'|^2\right).
\end{aligned}
\end{equation}
For $K_1$, we apply \autoref{lem:gsqg-flat-chord} with exponent $\beta+2$, noting that $|\Delta_h\gamma \cdot \gamma'| \le |\gamma'|^2 h$,
\begin{equation}\label{eq:gsqg-gauge-K1-bound}
\begin{aligned}
\left\|K_1\right\|_{L^1}
&\ls F^{\beta+4} h^{-2-\beta} (\|\gamma\|_{H^2} + \|\zeta\|_{H^2}) \|\gamma'-\zeta'\|_{L^2} \int_{\T} \left|\Delta_h\gamma\cdot\gamma'\right|\cdot |\Delta_h\gamma'|^2\,\d x\\
&\ls F^{\beta+4} h^{-1-\beta} (\|\gamma\|_{H^2} + \|\zeta\|_{H^2}) \|\gamma'-\zeta'\|_{L^2} \cdot |\gamma'|^2 \|\Delta_h\gamma'\|_{L^2}^2\\
&\ls F^{\beta+4} h^{1-\beta} (\|\gamma\|_{H^2} + \|\zeta\|_{H^2})^3 |\gamma'|^2 \|\gamma'-\zeta'\|_{L^2}.
\end{aligned}
\end{equation}
For $K_2$, the bounds are similar: Note that $|\Delta_h(\gamma-\zeta)\cdot\gamma'|\leq h^{1/2} |\gamma'| \|\gamma'-\zeta'\|_{L^2}$, we get
\begin{equation}\label{eq:gsqg-gauge-K2-bound}
\begin{aligned}
\left\|K_2\right\|_{L^1}
&\leq F^{\beta+2}h^{-3/2-\beta}
|\gamma'|\left\|\Delta_h\gamma'\right\|_{L^2}^2 \left\|\gamma'-\zeta'\right\|_{L^2}\\
&\leq F^{\beta+2}h^{1/2-\beta}|\gamma'|
\left\|\gamma''\right\|_{L^2}^2 \left\|\gamma'-\zeta'\right\|_{L^2}.
\end{aligned}
\end{equation}

For $K_3$, using $|\Delta_h\zeta|\leq|\zeta'|h$ and
H\"older's inequality, we obtain
\begin{equation}\label{eq:gsqg-gauge-K3-bound}
\begin{aligned}
\left\|K_3\right\|_{L^1}
&\leq F^{\beta+2}|\zeta'|h^{-1-\beta}
\int_\T|\gamma'-\zeta'|\,|\Delta_h\gamma'|^2\,\d x\\
&\leq F^{\beta+2}|\zeta'|h^{-1-\beta}
\left\|\gamma'-\zeta'\right\|_{L^2}
\left\|\Delta_h\gamma'\right\|_{L^2}\left\|\Delta_h\gamma'\right\|_{L^\infty}\\
&\leq F^{\beta+2}|\zeta'|h^{1/2-\beta}
\left\|\gamma''\right\|_{L^2}^2
\left\|\gamma'-\zeta'\right\|_{L^2}.
\end{aligned}
\end{equation}

The bound for $K_4$ needs finer treatment.
By an argument similar to \eqref{eq:tangent-square-difference}, we obtain
$$
\begin{aligned}
\left\||\Delta_h\gamma'|^2-|\Delta_h\zeta'|^2\right\|_{L^1}
&\leq\left\|\Delta_h\gamma'+\Delta_h\zeta'\right\|_{L^2}
\left\|\Delta_h(\gamma'-\zeta')\right\|_{L^2}\\
&\leq 2h
\left(\left\|\gamma''\right\|_{L^2}+\left\|\zeta''\right\|_{L^2}\right)
\left\|\gamma'-\zeta'\right\|_{L^2}.
\end{aligned}
$$
This, combined with $|\Delta_h\zeta \cdot \zeta'|\leq |\zeta'|^2 h$, leads us to
\begin{equation}\label{eq:gsqg-gauge-K4-bound}
\begin{aligned}
\left\|K_4\right\|_{L^1} & \le F^{\beta+2}h^{-1-\beta}|\zeta'|^2 \left\||\Delta_h\gamma'|^2-|\Delta_h\zeta'|^2\right\|_{L^1} \\ 
& \leq 2F^{\beta+2}h^{-\beta}|\zeta'|^2
\left(\left\|\gamma''\right\|_{L^2}+\left\|\zeta''\right\|_{L^2}\right)
\left\|\gamma'-\zeta'\right\|_{L^2}.
\end{aligned}
\end{equation}

We still need to estimate the latter term in \eqref{eq:gauge-difference}.
The mean value theorem and \eqref{eq:speed-difference} give
$$
\bigl||\gamma'|^{-2}-|\zeta'|^{-2}\bigr|
\leq2F^3\bigl||\gamma'|-|\zeta'|\bigr|
\ls F^3\left\|\gamma'-\zeta'\right\|_{L^2}.
$$
Hence, 
\begin{equation}\label{eq:gsqg-gauge-K5-bound}
\begin{aligned}
K_5(h):=\ &\left||\gamma'|^{-2}-|\zeta'|^{-2} \right|   
\int_\T\frac{|\Delta_h\zeta\cdot\zeta'|}{|\Delta_h\zeta|^{\beta+2}}
|\Delta_h\zeta'|^2\,\d x
\ls F^{3}\left\|\gamma'-\zeta'\right\|_{L^2} \int_\T\frac{|\Delta_h\zeta\cdot\zeta'|}{|\Delta_h\zeta|^{\beta+2}}
|\Delta_h\zeta'|^2\,\d x\\
\ls\ & F^{\beta+5}|\zeta'|^2\left\|\zeta''\right\|_{L^2}^2 h^{1-\beta}\left\|\gamma'-\zeta'\right\|_{L^2}.
\end{aligned}
\end{equation}

Conbining \eqref{eq:gsqg-gauge-K1-bound}-\eqref{eq:gsqg-gauge-K5-bound} together, we estimate the right hand side of \eqref{eq:gauge-difference} as
$$
\begin{aligned}
\left\|g_\gamma-g_\zeta\right\|_{L^1}
\ls\ & |\gamma'|^{-2} \int_{\T} \left\|K_1 \right\|_{L^1} \mathrm{d}h + F^2\int_{\T} \sum_{j=2}^4 \left\|K_j \right\|_{L^1} \mathrm{d}h + \int_{\T} K_5(h) \mathrm{d}h  \\ 
\ls\ & F^{\beta+5}
\left(1+\left\|\gamma\right\|_{H^2}+\left\|\zeta\right\|_{H^2}\right)^4
\left\|\gamma'-\zeta'\right\|_{L^2},
\end{aligned}
$$
which becomes
\eqref{eq:gsqg-g-stability}. \eqref{eq:gsqg-lambda-stability} follows directly from $\left\|\lambda_\gamma-\lambda_\zeta\right\|_{L^\infty}\le 2 \left\|g_\gamma-g_\zeta\right\|_{L^1}$.
\end{proof}

\subsection{The \texorpdfstring{$H^1$}{H1} contraction}
\label{sec:uniqueness}

Let $J$ be a common time interval on which $\gamma,\zeta$ are two solutions of
\eqref{eq:contour-equation}--\eqref{eq:lambda-definition} with the uniform
$H^2$ and arc--chord control supplied by
\autoref{prop:patch-constant-speed-equivalence}.

\begin{proposition}\label{prop:gsqg-H1-stability}
For almost every $t\in J$,
\begin{equation}\label{eq:gsqg-H1-stability}
\frac{\d}{\d t}\left\|\gamma-\zeta\right\|_{H^1}^2
\ls F^{\beta+5}
\left(1+\left\|\gamma\right\|_{H^2}+\left\|\zeta\right\|_{H^2}\right)^5
\left\|\gamma-\zeta\right\|_{H^1}^2.
\end{equation}
\end{proposition}

\begin{proof}
\textbf{(I) Zero-order estimate.}
We bound the $L^2$ norm at first. Pairing the difference of the two solutions with $\gamma-\zeta$ gives

\begin{equation}\label{eq:zero-order-split}
\begin{aligned}
\frac{1}{2}\frac{\d}{\d t}\left\|\gamma-\zeta\right\|_{L^2}^2
= \langle \mathcal{N}_\beta(\gamma) - \mathcal{N}_\beta(\zeta) , \gamma-\zeta \rangle + \langle \lambda_\gamma \gamma' - \lambda_\zeta \zeta', \gamma-\zeta \rangle \\
\le \left(\left\|\mathcal{N}_\beta(\gamma) - \mathcal{N}_\beta(\zeta)\right\|_{L^2} + \left\|\lambda_\gamma \gamma' - \lambda_\zeta \zeta'\right\|_{L^2}\right) \left\|\gamma-\zeta\right\|_{L^2}. 
\end{aligned}
\end{equation}

The normal difference here, denoted by $J_{01}$, is split into two parts as
\begin{equation}\label{eq:normal-difference}
\begin{aligned}
\Nop(\gamma)-\Nop(\zeta)
&=\int_\T\frac{\Delta_h(\gamma'-\zeta')}{|\Delta_h\gamma|^\beta}\,\d h+\int_\T\Delta_h\zeta'
\left(|\Delta_h\gamma|^{-\beta}-|\Delta_h\zeta|^{-\beta}\right)\,\d h \\
&=: J_{011}+J_{012}.
\end{aligned}
\end{equation}
The first term is bounded by Minkowski's inequality as
\begin{equation}\label{eq:J011-bound}
\begin{aligned}
\left\|J_{011}\right\|_{L^2}&\le \int_{\T} \left\|\frac{\Delta_h(\gamma'-\zeta')}{|\Delta_h\gamma|^\beta}\right\|_{L^2}\,\d h\\
&\ls F^\beta
\left\|\gamma'-\zeta'\right\|_{L^2}
\int_0^\pi h^{-\beta}\,\d h.
\end{aligned}
\end{equation}
For the second term, \autoref{lem:gsqg-flat-chord} and the bound
$$
\left\|\Delta_h\zeta'\right\|_{L^2}
\leq |h|\left\|\zeta''\right\|_{L^2}
$$
give, by Minkowski's and H\"older's inequalities,
\begin{equation}\label{eq:J012-bound}
\begin{aligned}
\left\|J_{012}\right\|_{L^2}
&\leq \int_{\T}\left\|\Delta_h\zeta'\right\|_{L^2}
\left\||\Delta_h\gamma|^{-\beta}-|\Delta_h\zeta|^{-\beta}\right\|_{L^\infty}\,\d h\\
&\ls F^{\beta+2}
\left(\left\|\gamma\right\|_{H^2}+\left\|\zeta\right\|_{H^2}\right)
\left\|\zeta''\right\|_{L^2}
\left\|\gamma'-\zeta'\right\|_{L^2}
\int_0^\pi h^{1-\beta}\,\d h.
\end{aligned}
\end{equation}

The gauge term satisfies, by \autoref{lem:single-gauge} and
\autoref{prop:gauge-stability},
\begin{equation}\label{eq:J02-bound}
\begin{aligned}
\left\|J_{02}\right\|_{L^2}:=\left\|\lambda_\gamma\gamma'-\lambda_\zeta\zeta'\right\|_{L^2}
&\leq \left\|\lambda_\gamma\right\|_{L^\infty}
\left\|\gamma'-\zeta'\right\|_{L^2}
+\left\|\lambda_\gamma-\lambda_\zeta\right\|_{L^\infty}
\left\|\zeta'\right\|_{L^2}\\
&\ls F^{\beta+5} \left(1+\left\|\gamma\right\|_{H^2}+\left\|\zeta\right\|_{H^2}\right)^5
\left\|\gamma'-\zeta'\right\|_{L^2}.
\end{aligned}
\end{equation}

Putting \eqref{eq:J011-bound}--\eqref{eq:J02-bound} into \eqref{eq:zero-order-split}
gives
\begin{equation}\label{eq:gsqg-L2-conclusion}
\frac{\d}{\d t}\left\|\gamma-\zeta\right\|_{L^2}^2
\ls F^{\beta+5}\left(1+\left\|\gamma\right\|_{H^2}+\left\|\zeta\right\|_{H^2}\right)^5
\left\|\gamma-\zeta\right\|_{H^1}^2.
\end{equation}

\textbf{(II) First-order estimate: Normal term.}
For the first-order part, write
\begin{equation}\label{eq:first-order-split}
\frac12\frac{\d}{\d t}\left\|\gamma'-\zeta'\right\|_{L^2}^2
=J_1+J_2,
\end{equation}
where
$$
\begin{aligned}
J_1&:=\int_\T(\gamma'-\zeta')\cdot
(\Nop(\gamma)-\Nop(\zeta))'\,\d x,\\
J_2&:=\int_\T(\gamma'-\zeta')\cdot
(\lambda_\gamma\gamma'-\lambda_\zeta\zeta')'\,\d x.
\end{aligned}
$$
We deal with $J_1$ first.
Differentiating the normal operator gives
$$
\bigl(\Nop(\gamma)\bigr)'
=\int_\T\frac{\Delta_h\gamma''}{|\Delta_h\gamma|^\beta}\,\d h
-\beta\int_\T
\frac{(\Delta_h\gamma\cdot\Delta_h\gamma')\Delta_h\gamma'}
{|\Delta_h\gamma|^{\beta+2}}\,\d h.
$$
Thus we can split $J_1 = J_{11} + J_{12}$ as
$$
\begin{aligned}
J_{11}&:=\iint_{\T^2}(\gamma'-\zeta')(x)\cdot
\left(\frac{\Delta_h\gamma''(x)}{|\Delta_h\gamma(x)|^\beta}
-\frac{\Delta_h\zeta''(x)}{|\Delta_h\zeta(x)|^\beta}\right)\,\d h\,\d x,\\
J_{12}&:=-\beta\iint_{\T^2}(\gamma'-\zeta')(x)\cdot
\left(\frac{(\Delta_h\gamma(x)\cdot\Delta_h\gamma'(x))\Delta_h\gamma'(x)}
{|\Delta_h\gamma(x)|^{\beta+2}}
-\frac{(\Delta_h\zeta(x)\cdot\Delta_h\zeta'(x))\Delta_h\zeta'(x)}
{|\Delta_h\zeta(x)|^{\beta+2}}\right)\,\d h\,\d x.
\end{aligned}
$$

\emph{Derivative falling on the numerator.}
The corresponding part of $J_1$ splits as
$$
\begin{aligned}
J_{11}&=\iint_{\T^2}(\gamma'-\zeta')(x)\cdot
\left(
\frac{\Delta_h\gamma''(x)}{|\Delta_h\gamma(x)|^\beta}
-\frac{\Delta_h\zeta''(x)}{|\Delta_h\zeta(x)|^\beta}
\right)\,\d h\,\d x\\
 &=\iint_{\T^2}(\gamma'-\zeta')(x)\cdot
\frac{\Delta_h(\gamma''-\zeta'')(x)}{|\Delta_h\gamma(x)|^\beta}
\,\d h\,\d x\\
&\quad+\iint_{\T^2}(\gamma'-\zeta')(x)\cdot\Delta_h\zeta''(x)
\left(|\Delta_h\gamma(x)|^{-\beta}-|\Delta_h\zeta(x)|^{-\beta}\right)
\,\d h\,\d x\\
&=: J_{111}+J_{112}.
\end{aligned}
$$

To overcome the derivative loss issue in $J_{111}$, we need to integrate by parts in $x$.
Write this term using the symmetry of the integrand,
\begin{equation}\label{eq:gsqg-J111-calculation}
\begin{aligned}
J_{111}
&=\iint_{\T^2}(\gamma'-\zeta')(x)\cdot
\frac{(\gamma''-\zeta'')(x) - (\gamma''-\zeta'')(y)}{|\gamma(x) - \gamma(y)|^\beta}
\,\d x\,\d y\\
&=\frac{1}{2}\iint_{\T^2}\left[(\gamma'-\zeta')(x) - (\gamma'-\zeta')(y)\right]\cdot
\frac{(\gamma''-\zeta'')(x) - (\gamma''-\zeta'')(y)}{|\gamma(x) - \gamma(y)|^\beta}
\,\d x\,\d y\\
&=\frac{1}{2}\iint_{\T^2}\Delta_h(\gamma'-\zeta')(x)\cdot
\frac{\Delta_h(\gamma''-\zeta'')(x)}{|\Delta_h\gamma(x)|^\beta}
\,\d h\,\d x\\
&=-\frac14\iint_{\T^2}|\Delta_h(\gamma'-\zeta')(x)|^2
\left(|\Delta_h\gamma(x)|^{-\beta}\right)'\,\d h\,\d x
\end{aligned}
\end{equation}
after integration by parts in $x$.
By \eqref{eq:CV-bound},
$$
\left|\left(|\Delta_h\gamma(x)|^{-\beta}\right)'\right|
=\beta \frac{\left|\Delta_h\gamma(x)\cdot \Delta_h\gamma'(x)\right|}{|\Delta_h\gamma(x)|^{\beta+2}}\ls F^{\beta+2}|h|^{-\beta}
\left|\int_{x-h}^{x}|\gamma''(r)|^2\,\d r\right|.
$$
Consequently, 
\begin{equation}\label{eq:gsqg-J111-conclusion}
|J_{111}|\ls F^{\beta+2} \left\|\gamma'' \right\|_{L^2}^2
\left\|\gamma'-\zeta'\right\|_{L^2}^2.
\end{equation}

By symmetry in $x,y$,
\begin{equation}\label{eq:gsqg-J112-calculation}
J_{112}=\frac12\iint_{\T^2}\Delta_h(\gamma'-\zeta')\cdot\Delta_h\zeta''\left(|\Delta_h\gamma|^{-\beta}-|\Delta_h\zeta|^{-\beta}\right)\,\d h\,\d x.
\end{equation}
For this term, \eqref{eq:normalized-kernel-difference} gives
$$
\begin{aligned}
|J_{112}|
&\leq\int_\T\left\|\gamma'-\zeta'\right\|_{L^2}\left\|\Delta_h\zeta''\right\|_{L^2}
\left\||\Delta_h\gamma|^{-\beta}-|\Delta_h\zeta|^{-\beta}\right\|_{L^\infty}\,\d h\\
&\ls F^{\beta+2}
\left(\left\|\gamma\right\|_{H^2}+\left\|\zeta\right\|_{H^2}\right)\left\|\zeta''\right\|_{L^2}
\left\|\gamma'-\zeta'\right\|_{L^2}^2\int_0^\pi h^{-\beta}\,\d h.
\end{aligned}
$$
\emph{Derivative falling on the denominator.}
The corresponding part of $J_1$ is
$$
\begin{aligned}
J_{12}=& -\beta\iint_{\T^2}(\gamma'-\zeta')(x)\cdot\left\{
\frac{(\Delta_h\gamma\cdot\Delta_h\gamma')\Delta_h\gamma'}
{|\Delta_h\gamma|^{\beta+2}}
-\frac{(\Delta_h\zeta\cdot\Delta_h\zeta')\Delta_h\zeta'}
{|\Delta_h\zeta|^{\beta+2}}
\right\}\,\d h\,\d x\\
=& -\beta\iint_{\T^2}(\gamma'-\zeta')(x)\cdot
\frac{\Delta_h\gamma\cdot\Delta_h\gamma'}
{|\Delta_h\gamma|^{\beta+2}}
\Delta_h(\gamma'-\zeta')\,\d h\,\d x\\
&{}-\beta\iint_{\T^2}(\gamma'-\zeta')(x)\cdot
\frac{\Delta_h\zeta'}{|\Delta_h\gamma|^{\beta+2}}
\left(\Delta_h\gamma\cdot\Delta_h\gamma'
-\Delta_h\zeta\cdot\Delta_h\zeta'\right)\,\d h\,\d x\\
&{}-\beta\iint_{\T^2}(\gamma'-\zeta')(x)\cdot
(\Delta_h\zeta\cdot\Delta_h\zeta')\Delta_h\zeta'
\left(|\Delta_h\gamma|^{-\beta-2}-|\Delta_h\zeta|^{-\beta-2}\right)
\,\d h\,\d x\\
=&: -\beta\left(J_{121}+J_{122}+J_{123}\right).
\end{aligned}
$$
For $J_{121}$, write the kernel in the variables $x,y$ and average under
$x\leftrightarrow y$, as in the calculation of $J_{111}$:
\begin{equation}\label{eq:gsqg-J121-calculation}
\begin{aligned}
\beta J_{121}
&=\beta\iint_{\T^2}(\gamma'-\zeta')(x)\cdot
\frac{(\gamma'-\zeta')(x)-(\gamma'-\zeta')(y)}
{|\gamma(x)-\gamma(y)|^{\beta+2}}\\
&\quad\cdot\left[(\gamma(x)-\gamma(y))\cdot(\gamma'(x)-\gamma'(y))\right]
\,\d x\,\d y\\
&=\frac{\beta}{2}\iint_{\T^2}
\frac{\left|(\gamma'-\zeta')(x)-(\gamma'-\zeta')(y)\right|^2}
{|\gamma(x)-\gamma(y)|^{\beta+2}}\\
&\quad\cdot\left[(\gamma(x)-\gamma(y))\cdot(\gamma'(x)-\gamma'(y))\right]
\,\d x\,\d y\\
&=\frac{\beta}{2}\iint_{\T^2}
|\Delta_h(\gamma'-\zeta')|^2
\frac{\Delta_h\gamma\cdot\Delta_h\gamma'}
{|\Delta_h\gamma|^{\beta+2}}\,\d h\,\d x\\
&=-\frac12\iint_{\T^2}
|\Delta_h(\gamma'-\zeta')|^2
\left(|\Delta_h\gamma|^{-\beta}\right)'\,\d h\,\d x
=2J_{111},
\end{aligned}
\end{equation}
so it obeys the same bound as \eqref{eq:gsqg-J111-conclusion}.
By symmetry in $x,y$,
\begin{equation}\label{eq:gsqg-J122-calculation}
J_{122}=\frac12\iint_{\T^2}\frac{\Delta_h(\gamma'-\zeta')\cdot\Delta_h\zeta'}{|\Delta_h\gamma|^{\beta+2}}\left(\Delta_h\gamma\cdot\Delta_h\gamma'-\Delta_h\zeta\cdot\Delta_h\zeta'\right)\,\d h\,\d x.
\end{equation}
For $J_{122}$, use \eqref{eq:H-relative} and
$\left\|\Delta_h\zeta'\right\|_{L^\infty}\leq|h|^{1/2}\left\|\zeta''\right\|_{L^2}$:
\begin{align}
|J_{122}|
&\leq F^{\beta+2}\int_\T |h|^{-\beta-2}
\left\|\gamma'-\zeta'\right\|_{L^2}\left\|\Delta_h\zeta'\right\|_{L^\infty} \cdot 
\left\|\Delta_h\gamma\cdot\Delta_h\gamma'
-\Delta_h\zeta\cdot\Delta_h\zeta'\right\|_{L^2}\,\d h\notag\\
&\ls F^{\beta+2}\left\|\zeta''\right\|_{L^2}
\left(\left\|\gamma\right\|_{H^2}+\left\|\zeta\right\|_{H^2}\right)
\left\|\gamma'-\zeta'\right\|_{L^2}^2\int_0^\pi h^{-\beta}\,\d h.
\label{eq:J122-bound}
\end{align}

For $J_{123}$, apply \autoref{lem:gsqg-flat-chord} with exponent $\beta+2$.
By \eqref{eq:CV-bound} and \eqref{eq:basic-difference-L2},
$$
\begin{aligned}
\left\|\Delta_h\zeta\cdot\Delta_h\zeta'\right\|_{L^\infty}
&\ls |h|^2\left\|\zeta''\right\|_{L^2}^2,\\
\left\|\Delta_h\zeta'\right\|_{L^2}
&\leq |h|\left\|\zeta''\right\|_{L^2}.
\end{aligned}
$$
H\"older's inequality gives
\begin{align}
|J_{123}|
&\leq\int_\T\left\|\gamma'-\zeta'\right\|_{L^2}
\left\|\Delta_h\zeta\cdot\Delta_h\zeta'\right\|_{L^\infty}
\left\|\Delta_h\zeta'\right\|_{L^2}\notag\\
&\qquad{}\cdot
\left\||\Delta_h\gamma|^{-\beta-2}-|\Delta_h\zeta|^{-\beta-2}\right\|_{L^\infty}
\,\d h\notag\\
&\ls F^{\beta+4}
\left(\left\|\gamma\right\|_{H^2}+\left\|\zeta\right\|_{H^2}\right)
\left\|\zeta''\right\|_{L^2}^3
\left\|\gamma'-\zeta'\right\|_{L^2}^2
\int_0^\pi h^{1-\beta}\,\d h.
\label{eq:J123-bound}
\end{align}

By \eqref{eq:gsqg-J121-calculation},
\begin{equation}\label{eq:normal-energy-decomposition}
J_1=-J_{111}+J_{112}-\beta J_{122}-\beta J_{123}.
\end{equation}
Combining \eqref{eq:gsqg-J111-conclusion}, the estimate for $J_{112}$ above, and
\eqref{eq:J122-bound}--\eqref{eq:J123-bound} gives
\begin{equation}\label{eq:J1-final}
|J_1|\ls F^{\beta+5}
\left(1+\left\|\gamma\right\|_{H^2}+\left\|\zeta\right\|_{H^2}\right)^4
\left\|\gamma'-\zeta'\right\|_{L^2}^2.
\end{equation}

\textbf{(III) First-order estimate: Tangential term.}
Split $J_2$ as
\begin{equation}\label{eq:tangential-decomposition}
\begin{aligned}
J_2
=\int_\T(\gamma'-\zeta')\cdot
\bigl(\lambda_\gamma(\gamma'-\zeta')\bigr)'\,\d x + \int_\T(\gamma'-\zeta')\cdot
\bigl((\lambda_\gamma-\lambda_\zeta)\zeta'\bigr)'\,\d x =: J_{21}+J_{22}.
\end{aligned}
\end{equation}
For the first integral, periodic integration by parts yields
$$
\begin{aligned}
J_{21}&=\int_\T(\gamma'-\zeta')\cdot
\bigl(\lambda_\gamma(\gamma'-\zeta')\bigr)'\,\d x\\
&=\int_\T(\lambda_\gamma')|\gamma'-\zeta'|^2\,\d x
+\frac12\int_\T\lambda_\gamma
\bigl(|\gamma'-\zeta'|^2\bigr)'\,\d x\\
&=\frac12\int_\T(\lambda_\gamma')
|\gamma'-\zeta'|^2\,\d x,
\end{aligned}
$$
hence the $L^\infty$ estimate in \eqref{eq:g-single-F} gives
\begin{align}
|J_{21}|
\leq\frac12\left\|\lambda_\gamma'\right\|_{L^\infty}
\left\|\gamma'-\zeta'\right\|_{L^2}^2 \ls F^{\beta+2}
\left\|\gamma''\right\|_{L^2}^2
\left\|\gamma'-\zeta'\right\|_{L^2}^2.
\label{eq:J21-bound}
\end{align}

For $J_{22}$, we integrate by parts before expanding the derivative, and use
$\gamma''\cdot\gamma'=\zeta''\cdot\zeta'=0$:
$$
\begin{aligned}
J_{22}&=\int_\T(\gamma'-\zeta')\cdot
\bigl((\lambda_\gamma-\lambda_\zeta)\zeta'\bigr)'\,\d x\\
&=-\int_\T(\gamma''-\zeta'') \cdot \zeta'(\lambda_\gamma-\lambda_\zeta)\,\d x\\
&=-\int_\T
\gamma''\cdot\zeta'(\lambda_\gamma-\lambda_\zeta)\,\d x = \int_\T
\gamma''\cdot(\gamma'-\zeta')(\lambda_\gamma-\lambda_\zeta)\,\d x,
\end{aligned}
$$
We obtain by Cauchy--Schwarz and \eqref{eq:gsqg-lambda-stability} that
\begin{align}
|J_{22}|
&\leq\left\|\lambda_\gamma-\lambda_\zeta\right\|_{L^\infty}
\left\|\gamma''\right\|_{L^2}
\left\|\gamma'-\zeta'\right\|_{L^2}\notag\\
&\ls F^{\beta+5}
\left(1+\left\|\gamma\right\|_{H^2}+\left\|\zeta\right\|_{H^2}\right)^4
\left\|\gamma''\right\|_{L^2}
\left\|\gamma'-\zeta'\right\|_{L^2}^2.
\label{eq:J22-bound}
\end{align}
Combining \eqref{eq:J21-bound} and \eqref{eq:J22-bound}, we obtain
\begin{equation}\label{eq:tangential-final}
|J_2|\ls F^{\beta+5}
\left(1+\left\|\gamma\right\|_{H^2}+\left\|\zeta\right\|_{H^2}\right)^5
\left\|\gamma'-\zeta'\right\|_{L^2}^2.
\end{equation}

Finally, combining \eqref{eq:gsqg-L2-conclusion}, \eqref{eq:first-order-split}, \eqref{eq:J1-final}, and \eqref{eq:tangential-final} proves \eqref{eq:gsqg-H1-stability}.
\end{proof}

\begin{proof}[Proof of \autoref{thm:gsqg-uniqueness}]
Assume two $H^2$ patch solutions $\Omega_1(t),\Omega_2(t)$ on $[0,T]$ share the same initial data $\Omega_0$.
By \autoref{prop:patch-constant-speed-equivalence}, we can choose their constant-speed parametrizations $\gamma,\zeta$ and $\gamma(0)=\zeta(0)$.
Then \eqref{eq:gsqg-H1-stability} and Gronwall's inequality give $\gamma=\zeta$ on $[0,T]$, hence $\Omega_1(t)=\Omega_2(t)$.
\end{proof}

\section{The SQG Case $\beta=1$: local well-posedness}
\label{sec:endpoint}

We recall that, as derived in \autoref{sec:contour-dynamics}, the contour dynamic equation for SQG in the constant-speed gauge is
\begin{equation}\label{eq:endpoint-contour-equation}
\gamma_t=\mathcal N(\gamma)+\lambda_\gamma\gamma',
\end{equation}
where
\begin{equation}\label{eq:endpoint-normal-operator}
\mathcal N(\gamma)(x) :=\int_\T \frac{\Delta_h\gamma'(x)}{|\Delta_h\gamma(x)|}\,\d h,
\end{equation}
and
\begin{equation}\label{eq:endpoint-lambda-definition}
\lambda_\gamma(x) =\frac{x+\pi}{2\pi}\int_\T g_\gamma(y)\,\d y -\int_{-\pi}^{x}g_\gamma(y)\,\d y,
\end{equation}
with
\begin{equation}\label{eq:endpoint-gauge-exact}
g_\gamma(x) =-\frac{1}{2|\gamma'|^2}\int_\T \frac{\Delta_h\gamma(x)\cdot\gamma'(x)} {|\Delta_h\gamma(x)|^3}|\Delta_h\gamma'(x)|^2\,\d h.
\end{equation}
This gauge imposes $\gamma'\cdot \gamma''=0$ and $\lambda_\gamma(t,-\pi)=0$.
For an $H^2$ curve with finite arc--chord quantity,
$$
\frac{|\Delta_h\gamma'(x)|}{|\Delta_h\gamma(x)|}
\leq F(\gamma)\left\|\gamma''\right\|_{L^2}|h|^{-1/2},
$$
so the integral defining $\mathcal N(\gamma)$ is absolutely convergent.

We claim the main result of this section, which is the local well-posedness of contour dynamics in $H^2$.

\begin{proposition}[Local well-posedness of the SQG contour system]
\label{prop:SQG-contour-wellposedness}
Let $\gamma_0\in H^2(\T;\R^2)$ be a positively oriented simple closed constant-speed curve. 
There is $T>0$, depending only on $\|\gamma_0''\|_{L^2}$ and $F(\gamma_0)$, such that the system \eqref{eq:endpoint-contour-equation}--\eqref{eq:endpoint-gauge-exact} with $\gamma(0)=\gamma_0$ has a unique simple closed constant-speed solution $\gamma$, satisfying
\begin{equation}\label{eq:SQG-contour-time-regularity}
\gamma\in C([0,T];H^2), \qquad \gamma_t\in C([0,T];C^0).
\end{equation}
The time $T$ can be chosen so that
\begin{equation}\label{eq:SQG-local-geometric-control}
 \sup_{0\leq t\leq T}\|\gamma''(t)\|_{L^2} \leq 2\|\gamma_0''\|_{L^2},
\qquad \sup_{0\leq t\leq T}F(\gamma(t))\leq 4F(\gamma_0).
\end{equation}
Moreover, if a sequence of such initial curves $\gamma_{0,n}$ satisfies $\sup_n(\|\gamma_{0,n}\|_{H^2}+F(\gamma_{0,n}))<\infty$ and $\gamma_{0,n}\to\gamma_0$ in $H^1$, 
then a common $T>0$ can be chosen and the corresponding solutions converge to $\gamma$ in $C([0,T];H^1)$.
\end{proposition}

This section is organized as follows.
\autoref{sec:SQG-preliminaries} gives the interpolation and compactness results used below.
\autoref{sec:SQG-a-priori} establishes the energy, gauge, and velocity estimates. 
In \autoref{sec:SQG-local-existence}, smooth initial approximations and logarithmic continuation give solutions with a uniform lifespan. 
Compactness and direct passage to each term then yield a contour solution with a classical time derivative.
\autoref{sec:SQG-uniqueness} proves uniqueness in a refined $H^1$ estimate.
Finally, in \autoref{sec:SQG-time-continuity} we use uniqueness to prove strong $H^2$ time continuity and complete the proof of \autoref{thm:SQG-wellposedness}.
We denote $|X|\leq CY$ by $X\ls Y$ for a constant $C>0$ that is absolute in this section.

\subsection{Preliminaries}
\label{sec:SQG-preliminaries}
Before we proceed to the a priori estimates, we record some standard analysis tools we will use then. 
Firstly, we include the following vector-valued Arzel\`a--Ascoli theorem  for completeness.
\begin{lemma}\label{lem:banach-valued-compactness}
Let $X$ be a Banach space and $T>0$. A family $\mathcal F\subset C([0,T];X)$ is relatively compact in the uniform topology, if and only if the following conditions hold:
\begin{enumerate}[(i)]
\item For every $t\in[0,T]$, the set $\{f(t):f\in\mathcal F\}$ is relatively compact in $X$.
\item For every $\varepsilon>0$, there is $\delta>0$ such that, $\|f(t)-f(s)\|_X<\varepsilon$ for all $f\in\mathcal F$ and all $s,t\in[0,T]$ with $|t-s|<\delta$.
\end{enumerate}
\end{lemma}

Now we introduce the following Gagliardo--Nirenberg interpolation inequality.

\begin{lemma}\label{lem:periodic-interpolation}
For every $f\in H^2(\T)$, it holds
$$
\|f'\|_{L^\infty} \ls \|f\|_{L^\infty}^{1/3}\|f''\|_{L^2}^{2/3} +\|f\|_{L^\infty}.
$$
\end{lemma}

\begin{proof}
Set $A=\|f\|_{L^\infty}$ and $B=\|f''\|_{L^2}$.
For $0<r\leq\pi$, the fundamental theorem of calculus and
Cauchy--Schwarz give
$$
f(x+r)-f(x) =rf'(x)+\int_0^r(r-s)f''(x+s)\,\d s,\quad \Rightarrow\quad
\|f'\|_{L^\infty} \leq\frac{2A}{r}+\frac{r^{1/2}}{\sqrt3}B.
$$
If $AB=0$, then $f'=0$. Otherwise, taking $r=\min\{\pi,(A/B)^{2/3}\}$ yields what we desired.
\end{proof}

The following lemma is key to our energy estimates. It is simply an equivalent formulation of embedding $\dot{H}^1\hookrightarrow \dot{B}_{\infty,2}^{1/2}$ in Besov space. A proof is included for completeness.

\begin{lemma}\label{lem:SQG-square-modulus}
For $f\in H^1(\T;\R^m)$ and $h\in [0,\pi]$, define
$\omega_f(h)=\sup_{|s|\leq h}\left\|\Delta_s f\right\|_{L^\infty}$.
Then
\begin{equation}\label{eq:SQG-square-modulus}
 \int_0^\pi\frac{\omega_f(h)^2}{h^2}\,\d h
 \ls\left\|f'\right\|_{L^2}^2 .
\end{equation}
\end{lemma}

\begin{proof}
The proof is based on a standard Littlewood--Paley decomposition.
Without loss of generality, assume $f$ has mean zero. 
Since $f\in H^1(\T)$, we can write standard Fourier coefficients and dyadic frequency decompositions by
$$
 \widehat f(n)=\frac{1}{2\pi}\int_\T f(x)\e^{-\i nx}\,\d x,
 \qquad
 f_j(x)=\sum_{2^j\leq|n|<2^{j+1}}\widehat f(n)\e^{\i nx},
 \quad j\geq0.
$$
Then $f=\sum_{j\geq0}f_j$ in $H^1(\T)$ and hence uniformly. The Parseval's theorem gives
$$
\sum_{j\geq0}2^{2j}\left\|f_j\right\|_{L^2(\T)}^2 \leq2\pi\sum_{n\ne0}|n|^2|\widehat f(n)|^2=\left\|f'\right\|_{L^2(\T)}^2.
$$
Now we make estimates for each block $f_j$. With trivial estimate $|1-\e^{-\i ns}|\leq\min\{2,|n||s|\}$ and Cauchy--Schwarz, we have
$$
\begin{aligned}
\left\|\Delta_s f_j\right\|_{L^\infty(\T)}
&\leq\sum_{2^j\leq|n|<2^{j+1}} |1-\e^{-\i ns}|\,|\widehat f(n)|\\
&\ls\min\{1,2^j|s|\}\,2^{j/2} \left( \sum_{2^j\leq|n|<2^{j+1}}|\widehat f(n)|^2 \right)^{1/2}\\
&\ls\min\{1,2^j|s|\}\,2^{j/2} \left\|f_j\right\|_{L^2(\T)}.
\end{aligned}
$$
Set $a_j=2^j\left\|f_j\right\|_{L^2(\T)}$. Fix $h\in(0,\pi]$ and let $k\ge 0$ be such that $\pi2^{-k-1}<h\le \pi2^{-k}$.
Summing over $j\le k$ and $j>k$ gives
$$
2^{\frac{k}{2}}\omega_f(h)\ls\sum_{j=0}^k2^{-\frac{k-j}{2}}a_j+\sum_{j>k}2^{-\frac{j-k}{2}}a_j=\sum_{j\ge 0}2^{-\frac{|j-k|}{2}}a_j.
$$
Since $\sum_{\ell\in\Z}2^{-|\ell|/2}<\infty$, Cauchy--Schwarz now implies
$$
\left(\sum_{j\ge 0}2^{-\frac{|j-k|}{2}}a_j\right)^2 \ls \sum_{j\ge 0}2^{-\frac{|j-k|}{2}}a_j^2,
$$
and hence we reach at
$$
\frac{\omega_f(h)^2}{h}\ls \sum_{j\geq0}2^{-\frac{|j-k|}{2}}a_j^2.
$$
This doesn't directly applies to the integral in \eqref{eq:SQG-square-modulus} due to a logarithmic issue.
Instead, integrating over the dyadic $h$ intervals and summing, we obtain
$$
\begin{aligned}
\int_0^\pi\frac{\omega_f(h)^2}{h^2}\,\d h
&\ls \sum_{k\geq0} \int_{\pi2^{-k-1}}^{\pi2^{-k}}\frac{\omega_f(h)^2}{h^2}\,\d h \ls \sum_{k\geq0} \sum_{j\geq0}2^{-\frac{|j-k|}{2}}a_j^2\int_{\pi2^{-k-1}}^{\pi2^{-k}}\frac{\d h}{h} \\
&\ls \sum_{k\geq0} \sum_{j\geq0} 2^{-\frac{|j-k|}{2}}a_j^2 \ls \sum_{j\geq0} a_j^2 \ls \left\|f'\right\|_{L^2}^2.
\end{aligned}
$$
This proves \eqref{eq:SQG-square-modulus}.
\end{proof}

We will also use the following product estimate. The proof is standard.

\begin{lemma}\label{lem:SQG-moser}
Let $f,g\in H^2(\T)$, then
$$
\|(fg)''\|_{L^2}\leq C\left( \|f''\|_{L^2}\|g\|_{L^\infty} +\|f\|_{L^\infty}\|g''\|_{L^2}\right).
$$
\end{lemma}

\subsection{A priori estimates}
\label{sec:SQG-a-priori}

Now we are able to close the energy estimates for the SQG contour dynamics.
Similar to the gSQG case with $\beta<1$, we first establish the estimates for the tangential terms.

\begin{lemma}\label{lem:sqg-gauge-estimate}
Let $\gamma$ be a regular constant-speed $H^2$ curve, and let $F$ bound its arc--chord quantity.
Then the functions $g_\gamma$ and $\lambda_\gamma$, defined in \eqref{eq:endpoint-gauge-exact} and \eqref{eq:endpoint-lambda-definition}, satisfy
\begin{equation}\label{eq:sqg-gauge-derivative}
\begin{aligned}
&\left\|g_\gamma\right\|_{L^\infty}+\left\|\lambda_\gamma'\right\|_{L^\infty}+\left\|\lambda_\gamma\right\|_{L^\infty}
\ls F^3\left\|\gamma''\right\|_{L^2}^2.
\end{aligned}
\end{equation}
\end{lemma}

\begin{proof}
By $|\Delta_h\gamma|\le |\gamma'||h|$,
$$
\begin{aligned}
|g_\gamma(x)|&\ls \frac{1}{|\gamma'|^2}\int_\T\frac{|\Delta_h\gamma(x)||\gamma'(x)|}{|\Delta_h\gamma(x)|^{3}}|\Delta_h\gamma'(x)|^2\,\d h \\
&\ls F^3\int_\T\frac{|\Delta_h\gamma'(x)|^2}{h^2}\,\d h \ls F^3\left\|\gamma''\right\|_{L^2}^2.
\end{aligned}
$$
In the last step we used \autoref{lem:SQG-square-modulus} with $f=\gamma'$. The bound for $\lambda_\gamma$ follows from $\left\|\lambda_\gamma'\right\|_{L^\infty}\leq2\left\|g_\gamma\right\|_{L^\infty}$.
\end{proof}

Now we establish the energy estimates.
\begin{proposition}
\label{prop:SQG-second-order}
Let $\gamma\in C_tH^5$ be a regular constant-speed solution of \eqref{eq:endpoint-contour-equation}\footnote{The extra requirement of regularity ensures the existence of derivatives and may be unnecessary.}, and let $F$ bound its arc--chord quantity. Then it holds
\begin{equation}\label{eq:SQG-second-order-energy}
\left|\frac{\d}{\d t}\left\|\gamma''\right\|_{L^2}^2\right|
\ls F^3\left\|\gamma''\right\|_{L^2}^4,
\end{equation}
and
\begin{equation}\label{eq:SQG-velocity-bound}
\left\|\gamma_t\right\|_{L^\infty}
\ls F\left\|\gamma''\right\|_{L^2}+F^2\left\|\gamma''\right\|_{L^2}^2.
\end{equation}
\end{proposition}
\begin{proof}
\textbf{(I) Preliminaries.}
We first prove \eqref{eq:SQG-second-order-energy}. The time derivative of energy is
\begin{equation}\label{eq:SQG-second-order-decomposition}
\begin{aligned}
\frac12\frac{\d}{\d t}\left\|\gamma''\right\|_{L^2}^2
&=\int_\T\gamma''\cdot(\mathcal N(\gamma))''\,\d x
+\int_\T\gamma''\cdot(\lambda_\gamma\gamma')''\,\d x\\
&=N_1+N_2+N_3+N_4,
\end{aligned}
\end{equation}
where
$$
\begin{aligned}
N_1&:=\iint_{\T^2}\frac{\gamma''\cdot\Delta_h\gamma'''}{|\Delta_h\gamma|}\,\d h\,\d x,\\
N_2&:=2\iint_{\T^2}\gamma''\cdot\Delta_h\gamma''\left(|\Delta_h\gamma|^{-1}\right)'\,\d h\,\d x,\\
N_3&:=\iint_{\T^2}\gamma''\cdot\Delta_h\gamma'\left(|\Delta_h\gamma|^{-1}\right)''\,\d h\,\d x,\\
N_4&:=\int_\T\gamma''\cdot(\lambda_\gamma\gamma')''\,\d x.
\end{aligned}
$$
We write the derivatives of the kernel:
$$
\begin{aligned}
\left(|\Delta_h\gamma|^{-1}\right)'&=-\frac{\Delta_h\gamma\cdot\Delta_h\gamma'}{|\Delta_h\gamma|^3},\\
\left(|\Delta_h\gamma|^{-1}\right)''&=3\frac{(\Delta_h\gamma\cdot\Delta_h\gamma')^2}{|\Delta_h\gamma|^5}
-\frac{|\Delta_h\gamma'|^2}{|\Delta_h\gamma|^3}
-\frac{\Delta_h\gamma\cdot\Delta_h\gamma''}{|\Delta_h\gamma|^3}.
\end{aligned}
$$
Thus $N_3=N_{31}+N_{32}+N_{33}$, where
$$
\begin{aligned}
N_{31}&:=3\iint_{\T^2}\frac{(\gamma''\cdot\Delta_h\gamma')(\Delta_h\gamma\cdot\Delta_h\gamma')^2}{|\Delta_h\gamma|^5}\,\d h\,\d x,\\
N_{32}&:=-\iint_{\T^2}\frac{(\gamma''\cdot\Delta_h\gamma')|\Delta_h\gamma'|^2}{|\Delta_h\gamma|^3}\,\d h\,\d x,\\
N_{33}&:=-\iint_{\T^2}\frac{(\gamma''\cdot\Delta_h\gamma')(\Delta_h\gamma\cdot\Delta_h\gamma'')}{|\Delta_h\gamma|^3}\,\d h\,\d x.
\end{aligned}
$$
We use the modulus $\omega_{\gamma'}$ from \autoref{lem:SQG-square-modulus}. By \eqref{eq:chord-tangent-identity}, $|\Delta_h\gamma\cdot\Delta_h\gamma'| \leq |h|\omega_{\gamma'}(|h|)^2$. Together with \eqref{eq:arc-chord-lower}, this gives
\begin{equation}\label{eq:SQG-existence-kernel-x}
\left\|\left(|\Delta_h\gamma|^{-1}\right)'\right\|_{L^\infty} \leq F^3\frac{\omega_{\gamma'}(|h|)^2}{h^2}.
\end{equation}

\medskip

\textbf{(II) Estimates of $N_1$ and $N_2$.}
The symmetrization process used in \eqref{eq:gsqg-J111-calculation} give
\begin{equation}\label{eq:sqg-N1-calculation}
\begin{aligned}
N_1
&=\iint_{\T^2}\gamma''(x)\cdot\frac{\gamma'''(x) - \gamma'''(y)}{|\gamma(x) - \gamma(y)|}\,\d x\,\d y\\
&=\frac{1}{2}\iint_{\T^2}\left[\gamma''(x) - \gamma''(y)\right]\cdot\frac{\gamma'''(x) - \gamma'''(y)}{|\gamma(x) - \gamma(y)|}\,\d x\,\d y\\
&=\frac{1}{2}\iint_{\T^2}\Delta_h\gamma''(x)\cdot\frac{\Delta_h\gamma'''(x)}{|\Delta_h\gamma(x)|}\,\d h\,\d x
=\frac{1}{4}\iint_{\T^2}\frac{\left(|\Delta_h\gamma''(x)|^2\right)'}{|\Delta_h\gamma(x)|}\,\d h\,\d x \\ 
&=-\frac14\iint_{\T^2}|\Delta_h\gamma''(x)|^2 \left(|\Delta_h\gamma(x)|^{-1}\right)'\,\d h\,\d x
\end{aligned}
\end{equation}
Thus, \eqref{eq:SQG-existence-kernel-x} and \autoref{lem:SQG-square-modulus} give
\begin{equation}
\label{eq:SQG-N1-bound}
\begin{aligned}
|N_1| &\ls F^3\int_\T\left\|\Delta_h\gamma''\right\|_{L^2}^2 \frac{\omega_{\gamma'}(|h|)^2}{h^2}\,\d h\\
&\ls F^3\left\|\gamma''\right\|_{L^2}^2
\int_\T\frac{\omega_{\gamma'}(|h|)^2}{h^2}\,\d h \ls F^3\left\|\gamma''\right\|_{L^2}^4.
\end{aligned}
\end{equation}
The estimate for $N_2$ is completely similar. Indeed applying Cauchy--Schwarz gives the control in the second line in \eqref{eq:SQG-N1-bound}. 
We obtain
$$
|N_1|+|N_2| \ls F^3\left\|\gamma''\right\|_{L^2}^4.
$$

\medskip

\textbf{(III) Estimate of $N_3$.}
The first two terms can be treated in similar manner: After putting everything into absolute value, it holds
$$
\begin{aligned}
|N_{31}| + |N_{32}| &\ls \iint_{\T^2} \frac{|\gamma''||\Delta_h\gamma'|^3}{|\Delta_h\gamma|^3}\,\d h\,\d x\\
&\ls F^3 \iint_{\T^2}\frac{|\gamma''||\Delta_h\gamma'|^3}{|h|^3}\,\d h\,\d x \ls F^3 \iint_{\T^2}\frac{|\gamma''|^2|\Delta_h\gamma'|^2}{|h|^2}\,\d h\,\d x \\
&\ls F^3\left\|\gamma''\right\|_{L^2}^2 \int_\T\frac{\omega_{\gamma'}(|h|)^2}{h^2}\,\d h \ls F^3\left\|\gamma''\right\|_{L^2}^4.
\end{aligned}
$$

For $N_{33}$, the constant-speed gives $\gamma'\cdot\gamma''=0$ and hence
$$
\Delta_h\gamma(x)\cdot\Delta_h\gamma''(x) =(\Delta_h\gamma(x)-\gamma'(x)h)\cdot\Delta_h\gamma''(x)-h\Delta_h\gamma'(x)\cdot\gamma''(x-h).
$$
The fundamental theorem of calculus also gives $\left\|\Delta_h\gamma-h\gamma'\right\|_{L^\infty} \leq|h|\omega_{\gamma'}(|h|)$, thus
\begin{equation}\label{eq:SQG-Dh-Dh2-L2}
\|\Delta_h\gamma\cdot\Delta_h\gamma''\|_{L^2} \ls |h|\omega_{\gamma'}(|h|)\|\gamma''\|_{L^2}.
\end{equation}
By applying Cauchy--Schwarz, $|\Delta_h\gamma'|\le \omega_{\gamma'}(|h|)$, and \autoref{lem:SQG-square-modulus}, we obtain
$$
\begin{aligned}
|N_{33}| &\ls F^3 \iint_{\T^2} \frac{|\gamma''||\Delta_h\gamma'||\Delta_h\gamma\cdot\Delta_h\gamma''|}{|h|^3}\,\d h\,\d x\\
&\ls F^3\left\|\gamma''\right\|_{L^2}^2 \int_\T\frac{\omega_{\gamma'}(|h|)^2}{h^2}\,\d h
\ls F^3\left\|\gamma''\right\|_{L^2}^4.
\end{aligned}
$$

\medskip

\textbf{(IV) Estimate of $N_4$ and velocity.}
It remains to estimate the tangential term $N_4$. Using $\gamma''\cdot\gamma'=0$ and integration by parts, as in the calculation of $J_{21}$, we obtain
$$
\begin{aligned}
N_4
&=2\int_\T\lambda_\gamma'|\gamma''|^2\,\d x+\frac12\int_\T\lambda_\gamma\bigl(|\gamma''|^2\bigr)'\,\d x\\
&=\frac32\int_\T\lambda_\gamma'|\gamma''|^2\,\d x.
\end{aligned}
$$
By \autoref{lem:sqg-gauge-estimate},
$$
|N_4|\leq\frac32\left\|\lambda_\gamma'\right\|_{L^\infty} \left\|\gamma''\right\|_{L^2}^2 \ls F^3\left\|\gamma''\right\|_{L^2}^4.
$$
Combining the preceding bounds with \eqref{eq:SQG-second-order-decomposition} proves \eqref{eq:SQG-second-order-energy}.

Finally we bound the velocity. \eqref{eq:endpoint-normal-operator}, \eqref{eq:arc-chord-lower}, and \eqref{eq:basic-difference-local} give
$$
\begin{aligned}
\left\|\mathcal N(\gamma)\right\|_{L^\infty}
&\leq F\int_\T\frac{\left\|\Delta_h\gamma'\right\|_{L^\infty}}{|h|}\,\d h\\
&\ls F\left\|\gamma''\right\|_{L^2},
\end{aligned}
$$
since $|h|^{-1/2}$ is integrable on $\T$. For the tangential velocity, \eqref{eq:endpoint-gauge-exact} and $|\Delta_h\gamma\cdot\gamma'|\leq|\Delta_h\gamma||\gamma'|$ yield
$$
\begin{aligned}
|\gamma'|\left\|g_\gamma\right\|_{L^\infty}
&\leq\frac{F^2}{2}\int_\T\frac{\omega_{\gamma'}(|h|)^2}{h^2}\,\d h\\
&\ls F^2\left\|\gamma''\right\|_{L^2}^2,
\end{aligned}
$$
where we used \eqref{eq:arc-chord-lower} and \autoref{lem:SQG-square-modulus}. Since $|\gamma'|$ is independent of $x$, \eqref{eq:endpoint-lambda-definition}
gives
$$
\left\|\lambda_\gamma\gamma'\right\|_{L^\infty}
\ls |\gamma'|\left\|g_\gamma\right\|_{L^\infty}
\ls F^2\left\|\gamma''\right\|_{L^2}^2.
$$
These two bounds and \eqref{eq:endpoint-contour-equation}
prove \eqref{eq:SQG-velocity-bound}.
\end{proof}

\subsection{Existence of contour solutions}
\label{sec:SQG-local-existence}

Now we construct solutions of the contour system \eqref{eq:endpoint-contour-equation}--\eqref{eq:endpoint-gauge-exact}.
With the a priori estimates established in \autoref{sec:SQG-a-priori}, we can try to approximate the solution we desire by smooth solutions, which has been extensively studied in \cite{Gancedo2008,Gancedo2022}.

\medskip

\textbf{(I) Approximation of the initial curve.} 
First we construct the smooth constant-speed curves to approximate the initial $H^2$ curve $\gamma_0$.
The following lemma is from \cite[Lemma~8.2]{Peletier2009}. For completeness, we give a proof in \autoref{app:constant-speed-approximation}.

\begin{lemma}[Smooth approximation of constant-speed curves]
\label{lem:constant-speed-approximation}
Let $\gamma_0\in H^2(\T;\R^2)$ be a simple closed constant-speed curve. 
There exists a sequence of smooth simple closed curves $\gamma_0^n$ with constant speed, such that
\begin{equation}\label{eq:SQG-initial-approximation}
 \gamma_0^n\longrightarrow\gamma_0\quad\hbox{in }H^2 \quad(n\to\infty),
 \qquad F(\gamma_0^n)\leq2F(\gamma_0)\quad\text{for all }n\geq1.
\end{equation}
\end{lemma}
From now on we choose a sequence of smooth initial curves $\gamma_0^n$ as in \autoref{lem:constant-speed-approximation}, with $\gamma_0^n\to\gamma_0$ in $H^2$ as $n\to\infty$. 
\medskip

\textbf{(II) Uniform lifespan on approximating solutions.}
For each $n$, \cite[Theorem~5.1]{Gancedo2008} gives a local $H^5$ solution $\gamma^n$ of the contour dynamics with constant-speed gauge \eqref{eq:endpoint-contour-equation}-\eqref{eq:endpoint-gauge-exact}, 
with initial datum $\gamma_0^n$. 
The first matter is, their lifespan $T_n$ may tend to zero as $n\to\infty$. 
To avoid this, we observe that the third-order derivative can be bounded double-exponentially, provided $H^2$ and arc--chord controls.
This follows from the second-order energy estimate \eqref{eq:SQG-second-order-energy}
and the logarithmic Gronwall inequality, as in \cite[Section~2.3]{Gancedo2022};
The proof we postpone to \autoref{app:SQG-smooth-continuation}.

\begin{lemma}\label{lem:SQG-logarithmic-continuation}
Let $\gamma\in C([0,T_*);H^5)$ be a regular simple constant-speed solution of \eqref{eq:endpoint-contour-equation}. Suppose that for some $M,F\geq1$, it holds
$$
\sup_{t<T_*}\|\gamma''(t)\|_{L^2}\leq M, \qquad \sup_{t<T_*}F(\gamma(t))\leq F.
$$
Then a constant $C=C(M,F)>0$ exists, such that
$$
\|\gamma'''(t)\|_{L^2}^2 \leq \|\gamma'''(0)\|_{L^2}^2 \e^{Ct}.
$$
In particular, the solution $\gamma$ can be continued beyond $T_*$ provided the $H^2$ and arc--chord bounds remain finite in $[0,T_*)$.
\end{lemma}

Thus, to obtain a uniform lifespan for the approximating solutions $\gamma^n$, it suffices to obtain uniform $H^2$ and arc--chord controls.
Let $T_n$ be the maximal $H^5$ lifespan of $\gamma^n$, and set $M=2\|\gamma_0''\|_{L^2}$. 
By \eqref{eq:SQG-initial-approximation}, we may assume for every $n$,
\begin{equation}\label{eq:SQG-approximate-initial-bounds}
\|(\gamma_0^n)''\|_{L^2}^2\leq\frac{M^2}{2}, \qquad F(\gamma_0^n)\leq2F(\gamma_0)
\end{equation}
by discarding finitely many terms. The definition of $F$ gives
$$
\begin{aligned}
\left|\frac1{F(\gamma^n(t))}-\frac1{F(\gamma^n(s))}\right|
&\le \sup_{x,h}\left|\frac{|\Delta_h\gamma^n(t,x)|-|\Delta_h\gamma^n(s,x)|}{|h|}\right| \\
&\le \|\left(\gamma^n(t)-\gamma^n(s)\right)'\|_{L^\infty}.
\end{aligned}
$$
For each $n$, the initial bounds \eqref{eq:SQG-approximate-initial-bounds} and $H^5$ continuity in time ensure the existence of a positive time $t_n\in (0,T_n]$ such that 
\begin{equation}\label{eq:SQG-approximate-spatial-bounds}
\|(\gamma^n)''(s)\|_{L^2}<M,\qquad F(\gamma^n(s))<4F(\gamma_0)\qquad(0\le s< t_n).
\end{equation}
Let $t_n$ be the maximal time with this property. Note that \autoref{lem:SQG-logarithmic-continuation} provides $T_n>t_n$ (unless they are both infinite).
For every $n$ and $0\leq t<t_n$, integrating \eqref{eq:SQG-second-order-energy} and \eqref{eq:SQG-velocity-bound} from $0$ to $t$ gives, with $C$ independent of $n$,
\begin{equation}\label{eq:SQG-existence-bootstrap}
\begin{aligned}
\|(\gamma^n)''(t)\|_{L^2}^2 &\leq\|(\gamma_0^n)''\|_{L^2}^2+CF(\gamma_0)^3M^4t,\\
\|\gamma^n(t)-\gamma_0^n\|_{L^\infty} &\leq\int_0^t\|\gamma_s^n(s)\|_{L^\infty}\,\d s\\
&\leq C\bigl(F(\gamma_0)M+F(\gamma_0)^2M^2\bigr)t.
\end{aligned}
\end{equation}
For the same time $t$, \eqref{eq:SQG-approximate-initial-bounds} and
\eqref{eq:SQG-approximate-spatial-bounds} give
$$
 \|(\gamma^n(t)-\gamma_0^n)''\|_{L^2}
 \leq\|(\gamma^n)''(t)\|_{L^2}+\|(\gamma_0^n)''\|_{L^2}
 \leq M+\frac M{\sqrt2}<2M.
$$
Applying \autoref{lem:periodic-interpolation} gives, for every $n$ and $0\leq t<t_n$,
\begin{equation}\label{eq:SQG-C1-time-increment}
\begin{aligned}
&\|(\gamma^n(t)-\gamma_0^n)'\|_{L^\infty}\\
\ls&\ \|\gamma^n(t)-\gamma_0^n\|_{L^\infty}^{1/3}\|(\gamma^n(t)-\gamma_0^n)''\|_{L^2}^{2/3}+\|\gamma^n(t)-\gamma_0^n\|_{L^\infty}\\
\leq&\ C(M,F(\gamma_0))(t^{1/3}+t).
\end{aligned}
\end{equation}
The constant $C(M,F(\gamma_0))$ is independent of $n$.
Choose $T\in(0,1]$, depending only on $M$ and $F(\gamma_0)$, small enough so that
\begin{equation}\label{eq:SQG-lifespan-time-choice}
\begin{aligned}
 &CF(\gamma_0)^3M^4T\leq\frac{M^2}{4},\\
 &C(M,F(\gamma_0))(T^{1/3}+T)\leq\frac1{8F(\gamma_0)}.
\end{aligned}
\end{equation}

We show that $t_n\ge T$, which will imply that $T_n\ge T$ for every $n$.
Indeed, if $t_n<T$, then the maximality of $t_n$ and continuity give
$$
\|(\gamma^n)''(t_n)\|_{L^2}^2=M^2 \qquad\hbox{or}\qquad F(\gamma^n(t_n))=4F(\gamma_0).
$$
However, for every $n$ and $0\leq t<\min\{t_n,T\}$, \eqref{eq:SQG-approximate-initial-bounds}, \eqref{eq:SQG-existence-bootstrap}, and \eqref{eq:SQG-lifespan-time-choice} give
$$
 \|(\gamma^n)''(t)\|_{L^2}^2
 \leq\frac{M^2}{2}+CF(\gamma_0)^3M^4T
 \leq\frac{3M^2}{4}<M^2,
$$
while \eqref{eq:SQG-C1-time-increment} and \eqref{eq:SQG-lifespan-time-choice} imply, for every $x$ and $h\ne0$,
$$
\begin{aligned}
 \frac{|\Delta_h\gamma^n(t,x)|}{|h|}
 &\geq\frac{|\Delta_h\gamma_0^n(x)|}{|h|}
       -\|(\gamma^n(t)-\gamma_0^n)'\|_{L^\infty}\\
 &\geq\frac1{2F(\gamma_0)}-\frac1{8F(\gamma_0)}
 =\frac3{8F(\gamma_0)}.
\end{aligned}
$$
Thus $F(\gamma^n(t))\leq8F(\gamma_0)/3<4F(\gamma_0)$.

If $t_n<T$, letting $t\uparrow t_n$ in the improved bounds above and using continuity, we see 
$$
\|(\gamma^n)''(t_n)\|_{L^2}^2\leq3M^2/4 \qquad\text{and}\qquad F(\gamma^n(t_n))\leq8F(\gamma_0)/3.
$$
This contradicts the estimates above, so $t_n\ge T$ for every $n$.

\textbf{(III) Boundedness and compactness.}
Now we are able to choose some $T>0$ depending only on $F(\gamma_0)$ and $\|\gamma_0''\|_{L^2}$, such that the approximating solutions $\gamma^n$ exist on $[0,T]$ and satisfy
\begin{equation}\label{eq:SQG-existence-uniform-bounds}
\sup_{t\in[0,T]}\|(\gamma^n)''(t)\|_{L^2}\leq M, \qquad \sup_{t\in[0,T]}F(\gamma^n(t))\leq4F(\gamma_0).
\end{equation}

By \autoref{lem:sqg-gauge-estimate} and \eqref{eq:SQG-velocity-bound}, for every $n$ and $t\in[0,T]$,
$$
\begin{aligned}
\|g_{\gamma^n(t)}\|_{L^\infty} +\|\lambda_{\gamma^n(t)}\|_{W^{1,\infty}}
&\leq CF(\gamma_0)^3M^2,\\
\|\gamma_t^n(t)\|_{L^\infty} &\leq C\bigl(F(\gamma_0)M+F(\gamma_0)^2M^2\bigr).
\end{aligned}
$$
Hence we have from \eqref{eq:SQG-initial-approximation} that
$$
\sup_n\|\gamma^n\|_{L^\infty([0,T];H^2)}<\infty.
$$

For every $n$ and $s,t\in[0,T]$, the velocity bound gives
$$
\|\gamma^n(t)-\gamma^n(s)\|_{L^\infty} \leq C(M,F(\gamma_0))|t-s|.
$$
Since $\|(\gamma^n(t)-\gamma^n(s))''\|_{L^2}\leq2M$, the $C^1$ continuity follows from \autoref{lem:periodic-interpolation}:
$$
\|\gamma^n(t)-\gamma^n(s)\|_{C^1}
\leq C(M,F(\gamma_0))\bigl(|t-s|^{1/3}+|t-s|\bigr).
$$
For each fixed $t\in[0,T]$, the uniform $H^2$ bound give relative compactness of $\{\gamma^n(t)\}_n$ in $C^1(\T)$, since embedding $H^2(\T)\Subset C^1(\T)$ is compact. 
Together with uniform equicontinuity, \autoref{lem:banach-valued-compactness}, applied with $X=C^1(\T;\R^2)$, gives a subsequence (still denoted by $\gamma^n$), such that
$$
\gamma^n\longrightarrow\gamma
\quad\hbox{in }C([0,T];C^1).
$$

For every $t\in[0,T]$, the uniform $H^2$ bound and $C^1$ convergence imply $\gamma^n(t)\rightharpoonup\gamma(t)$ in $H^2$. In particular, $\gamma(t)\in H^2$ with the same uniform bound.
If $t_j\to t$, again we have boundedness of $\|\gamma(t_j)\|_{H^2}$ and convergence to $\gamma(t)$ in $C^1$. This gives $\gamma(t_j)\rightharpoonup\gamma(t)$ weakly in $H^2$. Hence $\gamma\in C_w([0,T];H^2)$.
The convergence of initial datum claimed in \eqref{eq:SQG-initial-approximation} implies $\gamma(0)=\gamma_0$. 
Constant speed and the chord lower bound pass to the limit in $C^1$, preserving regularity, simplicity and orientation. Weak lower semicontinuity and \eqref{eq:SQG-existence-uniform-bounds} give \eqref{eq:SQG-local-geometric-control}.

By now we are still left to show that $\gamma\in C([0,T];H^2)$. The proof we take requires uniqueness of the solution and we postpone it to the end of this section. 

\textbf{(IV) Passage to the contour equation.}
We now show that $\gamma$ satisfies the contour equation \eqref{eq:endpoint-contour-equation}. We will show that the velocity terms converge in $C([0,T];C^0)$, which is sufficient to pass to the limit in  the equation that is integrated in time. 

For fixed $h\ne0$, $C_tC^1$ convergence and the chord lower bound give
convergence of the integrands
$$
\frac{\Delta_h(\gamma^n)'(t,x)}{|\Delta_h\gamma^n(t,x)|}
\longrightarrow\frac{\Delta_h\gamma'(t,x)}{|\Delta_h\gamma(t,x)|}.
$$
uniformly on $[0,T]\times\T$. Since they have a common majorant $C|h|^{-1/2}$, we have $\mathcal N(\gamma^n)\to\mathcal N(\gamma)$ in $C_tC^0$ by dominated convergence.

For the gauge term $\lambda_\gamma \gamma'$, we need the convergence of $g_\lambda$ first. Uniform  bounds on speed and chord imply
$$
\frac{|\Delta_h\gamma^n\cdot(\gamma^n)'|}{2|(\gamma^n)'|^2|\Delta_h\gamma^n|^3}
\leq\frac{CF(\gamma_0)^3}{h^2}.
$$
Thus,
$$
\begin{aligned}
&\left\|
\int_{|h|<\delta}
\frac{\Delta_h\gamma^n\cdot(\gamma^n)'} {2|(\gamma^n)'|^2|\Delta_h\gamma^n|^3}|\Delta_h(\gamma^n)'|^2\,dh
\right\|_{L_x^1}\\
\le&\ CF(\gamma_0)^3
\int_{|h|<\delta}
\frac{\|\Delta_h(\gamma^n)'(t)\|_{L_x^2}^2}{h^2}\,dh
\le 2CF(\gamma_0)^3 M^2\delta.
\end{aligned}
$$
Here, we used for $0<\delta<\pi$, the scaling property of the $L^2$ norm and the fundamental theorem of calculus to obtain
$$
\int_{|h|<\delta}\frac{\|\Delta_h(\gamma^n)'(t)\|_{L^2}^2}{h^2}\,\d h
\leq 2\int_{|h|<\delta}\|(\gamma^n)''\|_{L^2}^2\,\d h
\leq 2M^2\delta.
$$
This deal with the contribution from small $h$. The same estimate holds for $\gamma$. 
On $|h|\geq\delta$, the integrands converge uniformly by $C_tC^1$ convergence. Letting first $n\to\infty$ and then $\delta\to0$ gives $g_{\gamma^n}\to g_\gamma$ in $C_tL^1$.
By \eqref{eq:endpoint-lambda-definition}, we have $\|\lambda_{\gamma^n}-\lambda_\gamma\|_{L^\infty}\leq2\|g_{\gamma^n}-g_\gamma\|_{L^1}$.
Consequently $\lambda_{\gamma^n}\to\lambda_\gamma$ in $C_tC^0$, and $C_tC^1$ convergence of the curves gives $\lambda_{\gamma^n}(\gamma^n)' \to\lambda_\gamma\gamma'$ in $C_tC^0$.

The initial convergence \eqref{eq:SQG-initial-approximation} and the convergence of both velocity terms allow passage to the time-integrated equation:
$$
\gamma(t)=\gamma_0+
\int_0^t\bigl(\mathcal N(\gamma(s))+\lambda_{\gamma(s)}\gamma'(s)\bigr)\,\d s.
$$
The right-hand side has a continuous $C^0$-valued time derivative.
Thus $\gamma_t\in C([0,T];C^0)$ and
\eqref{eq:endpoint-contour-equation} holds for every $(t,x)$. The existence part of \autoref{prop:SQG-contour-wellposedness} is basically settled, except for the strong $H^2$ time continuity, which we will prove after uniqueness.

\subsection{Uniqueness}
\label{sec:SQG-uniqueness}

Now we compare constant-speed solutions $\gamma,\zeta$ of \eqref{eq:endpoint-contour-equation} we constructed above.
Since the solutions are in $C([0,T];C^1) \cap C_w([0,T];H^2)$, we note that the proof here shouldn't require strong $H^2$ time continuity.

Set
$$
M=8\left(1+\sup_{0\leq t\leq T}
\left\|\gamma(t)\right\|_{H^2}+\sup_{0\leq t\leq T}\left\|\zeta(t)\right\|_{H^2}\right).
$$
Let $F\geq1$ be their common arc--chord bound, as in \eqref{eq:gsqg_bound}. 
We make some conventions for the rest of this section.
All implicit constants in this section are absolute.  Dependence on $M$ and $F$ is displayed explicitly, and the estimates are uniform in the truncation. 
For $0<s\leq M/4$, define the modulus (related to $M$) as 
$$
\Phi(s)=s\log\frac{M}{s}.
$$
Set $\Phi(0)=0$.  Both $\Phi(s)$ and $s\Phi(s)$ are increasing on $[0,M/4]$,
and
$$
\int_{0^+}\frac{\d s}{\Phi(s)}=+\infty.
$$

\begin{lemma}\label{lem:critical-scale-splitting}
Let $f\in H^1(\T)$ has zero mean, then
\begin{align}
\int_{\T} \frac{\left\|\Delta_hf\right\|_{L^2}}{|h|}\,\d h
&\ls \left\|f \right\|_{L^2}\log\frac{\left\|f\right\|_{H^1}}{\left\|f\right\|_{L^2}}.
\label{eq:critical-linear-difference}
\end{align}

\end{lemma}

\begin{proof}
Without loss of generality, assume $f\not = 0$ and the length of $\T$ is 1 (by a scaling). It's been mentioned in \eqref{eq:basic-difference-L2} that for $h>0$,
$$
\left\|\Delta_hf\right\|_{L^2}
\leq\min\{2\left\|f\right\|_{L^2},\left\|f'\right\|_{L^2}h\}.
$$
Splitting the integral at some $\ell\in [0,1]$, we obtain
$$
\int_0^1\frac{\left\|\Delta_h f\right\|_{L^2}}{h}\,\d h \leq \left\|f'\right\|_{L^2} \ell - 2 \left\|f\right\|_{L^2}\log\ell.
$$
Take $\ell = \left\|f\right\|_{L^2}/\left\|f'\right\|_{L^2}$ proves \eqref{eq:critical-linear-difference}.
\end{proof}
With the modulus $\Phi$ defined above, we have from \autoref{lem:critical-scale-splitting} that
\begin{equation}\label{eq:sqg-difference-modulus}
\int_{\T} \frac{\left\|\Delta_h (\gamma'-\zeta')\right\|_{L^2}}{|h|}\,\d h \ls \Phi(\left\|\gamma'-\zeta'\right\|_{L^2}).
\end{equation}
This is the central tool for the energy estimate in the proof of uniqueness.
Now we can present the analogue of \autoref{prop:gauge-stability}. 
\begin{lemma}\label{lem:endpoint-gauge-difference}
Let $\gamma,\zeta\in H^2(\T;\R^2)$ be two constant-speed curves satisfying the uniform bounds above.  Then
\begin{align}
\left\|g_\gamma-g_\zeta\right\|_{L^1} &\ls M^4F^6\Phi(\left\|\gamma'-\zeta'\right\|_{L^2}),
\label{eq:sqg-g-difference}\\
\left\|\lambda_\gamma-\lambda_\zeta\right\|_{L^\infty} &\ls M^4F^6\Phi(\left\|\gamma'-\zeta'\right\|_{L^2}).
\label{eq:sqg-lambda-difference}
\end{align}
\end{lemma}

\begin{proof}
As in the proof of \autoref{prop:gauge-stability}, we split $g_\gamma - g_\zeta$ as \eqref{eq:gauge-difference-1-split} with $\beta=1$.
Since the power of $h$ in $K_1,K_2,K_3,K_5$ are strictly greater than $-1$ as shown in \eqref{eq:gsqg-gauge-K1-bound}-\eqref{eq:gsqg-gauge-K5-bound}, it suffices to control $K_4$.
Before we made estimate \eqref{eq:gsqg-gauge-K4-bound}, we actually got
$$
\left\||\Delta_h\gamma'|^2-|\Delta_h\zeta'|^2\right\|_{L^1} \leq\left\|\Delta_h\gamma'+\Delta_h\zeta'\right\|_{L^2} \left\|\Delta_h(\gamma'-\zeta')\right\|_{L^2}
$$
This, again conbined with $|\Delta_h\zeta \cdot \zeta'|\leq |\zeta'|^2 h$, lead us to
$$
\begin{aligned}
\left\|K_4\right\|_{L^1} & \le F^{3}h^{-2}|\zeta'|^2 \left\||\Delta_h\gamma'|^2-|\Delta_h\zeta'|^2\right\|_{L^1} \\ 
& \ls F^{3}h^{-1}|\zeta'|^2
\left(\left\|\gamma''\right\|_{L^2}+\left\|\zeta''\right\|_{L^2}\right)
\left\|\Delta_h(\gamma'-\zeta')\right\|_{L^2}.
\end{aligned}
$$

Using \eqref{eq:speed-arc-chord} and \eqref{eq:critical-linear-difference},
$$
\begin{aligned}
\left\|g_\gamma-g_\zeta\right\|_{L^1}
&\ls M^4F^6\left\|\gamma'-\zeta'\right\|_{L^2}
+M^4F^6\int_0^\pi \frac{\left\|\Delta_h(\gamma'-\zeta')\right\|_{L^2}}h\,\d h\\
&\ls M^4F^6\Phi(\left\|\gamma'-\zeta'\right\|_{L^2}).
\end{aligned}
$$
For $h<0$, translation gives the same bounds.
Finally, \eqref{eq:lambda-definition} gives $\left\|\lambda_\gamma-\lambda_\zeta\right\|_{L^\infty}
\leq2\left\|g_\gamma-g_\zeta\right\|_{L^1}$.
\end{proof}

Let $\gamma,\zeta$ be two solutions of \eqref{eq:endpoint-contour-equation}
on $[0,T]$, with positive constant speeds and the uniform $H^2$ and arc--chord bounds above. 
Assume
\begin{equation}\label{eq:SQG-weak-time-class}
\gamma,\zeta\in L^\infty([0,T];H^2(\T;\R^2))
\cap C^1([0,T];C^0(\T;\R^2))
\cap C^0([0,T];C^1(\T;\R^2)).
\end{equation}
This regularity is satisfied by the solutions constructed in \autoref{sec:SQG-local-existence}. We show for these solutions, the $H^1$ distance is controlled by Osgood-type inequality, which implies uniqueness and strong $H^2$ time continuity.
\begin{proposition}\label{prop:SQG-H1-stability}
For almost every $t\in[0,T]$,
\begin{equation}\label{eq:endpoint-Osgood-energy}
\frac{\d}{\d t}\left\|\gamma-\zeta\right\|_{H^1}^2
\ls M^5F^6\left\|\gamma-\zeta\right\|_{H^1} \Phi(\left\|\gamma-\zeta\right\|_{H^1}).
\end{equation}
Consequently, there is an absolute constant $C>0$ such that, for every
$t\in[0,T]$,
\begin{equation}\label{eq:SQG-H1-continuity}
\left\|\gamma(t)-\zeta(t)\right\|_{H^1}
\leq M\left(\frac{\left\|\gamma(0)-\zeta(0)\right\|_{H^1}}{M}\right)^{\exp(-CM^5F^6t)}.
\end{equation}
\end{proposition}

\begin{proof}
We again write \eqref{eq:zero-order-split} and \eqref{eq:first-order-split} with $\beta=1$.
The regularity \eqref{eq:SQG-weak-time-class} is enough to justify these differentiations on energy.

\textbf{(I) Zero-order estimate.} Split as \eqref{eq:normal-difference} with $\beta=1$.
We start with the term $J_{011}$. After the first step of \eqref{eq:J011-bound}, we apply \eqref{eq:critical-linear-difference} and get
\begin{equation}\label{eq:sqg-J011-bound}
\begin{aligned}
\left\|J_{011}\right\|_{L^2}&\le \int_{\T} \left\|\frac{\Delta_h(\gamma'-\zeta')}{|\Delta_h\gamma|}\right\|_{L^2}\,\d h\\
&\ls F \int_{0}^\pi \left\|\Delta_h(\gamma'-\zeta')\right\|_{L^2}\,\frac{\d h}{h} \ls F \Phi(\left\|\gamma'-\zeta'\right\|_{L^2}).
\end{aligned}
\end{equation}
$J_{012}$, as dealt with in \eqref{eq:J012-bound}, behaves well at $\beta=1$, hence it suffices to control $J_{02}$.
In the last step of \eqref{eq:J02-bound}, we use \eqref{eq:sqg-gauge-derivative} and \eqref{eq:sqg-lambda-difference} instead, which gives
$$
\left\|J_{02}\right\|_{L^2}
\ls M^5F^6\Phi(\left\|\gamma'-\zeta'\right\|_{L^2}).
$$
Combining these bounds with \eqref{eq:zero-order-split} and the monotonicity of $\Phi$ gives
\begin{equation}\label{eq:endpoint-L2-energy}
\frac12\frac{\d}{\d t}\left\|\gamma-\zeta\right\|_{L^2}^2
\ls M^5F^6\left\|\gamma-\zeta\right\|_{H^1}
\Phi(\left\|\gamma-\zeta\right\|_{H^1}).
\end{equation}

\textbf{(II) First-order estimate: Normal term.}
We start from \eqref{eq:first-order-split} now. 
For $J_{111}$, after the calculation in \eqref{eq:gsqg-J111-calculation}, we apply \eqref{eq:SQG-existence-kernel-x} and \eqref{eq:SQG-square-modulus} to obtain
$$
\begin{aligned}
J_{111}
&=-\frac14 \iint_{\T^2}|\Delta_h(\gamma'-\zeta')(x)|^2
\left(|\Delta_h\gamma(x)|^{-1}\right)'\,\d h\,\d x \\
&\ls F^3 \int_{\T}\frac{\omega_{\gamma'}(|h|)^2}{h^2}\left\|\Delta_h(\gamma'-\zeta')\right\|_{L^2}^2\,\d h\\
&\ls F^3 \left\|\gamma'-\zeta'\right\|_{L^2}^2 \int_{\T}\frac{\omega_{\gamma'}(|h|)^2}{h^2}\,\d h \ls M^2F^3\left\|\gamma'-\zeta'\right\|_{L^2}^2.
\end{aligned}
$$
For $J_{112}$, after the calculation in \eqref{eq:gsqg-J112-calculation}, we apply \eqref{eq:normalized-kernel-difference}, Cauchy-Schwarz, and
\eqref{eq:critical-linear-difference} to obtain
$$
\begin{aligned}
J_{112}&=\frac12\iint_{\T^2}\Delta_h(\gamma'-\zeta')\cdot\Delta_h\zeta''\left(|\Delta_h\gamma|^{-1}-|\Delta_h\zeta|^{-1}\right)\,\d h\,\d x \\
&\ls M F^3 \left\|\gamma'-\zeta'\right\|_{L^2} \iint_{\T^2} |\Delta_h(\gamma'-\zeta')||\Delta_h\zeta''||h|^{-1}\,\d h\,\d x \\
&\ls M^2 F^3\left\|\gamma'-\zeta'\right\|_{L^2}
\int_0^\pi\frac{\left\|\Delta_h(\gamma'-\zeta')\right\|_{L^2}}h\,\d h\\
&\ls M^2F^3\left\|\gamma'-\zeta'\right\|_{L^2}
\Phi(\left\|\gamma'-\zeta'\right\|_{L^2}).
\end{aligned}
$$

Now we turn to $J_{12}$.
Note that \eqref{eq:gsqg-J121-calculation} gives $J_{121}=2J_{111}$, so this term is already controlled. Also $J_{123}$ behaves well when $\beta=1$, as shown in \eqref{eq:J123-bound}. Hence it suffices to treat $J_{122}$.
For $J_{122}$, using \eqref{eq:gsqg-J122-calculation} allow us to make following estimate:
$$
\begin{aligned}
J_{122}&=\frac12\iint_{\T^2}\frac{\Delta_h(\gamma'-\zeta')\cdot\Delta_h\zeta'}{|\Delta_h\gamma|^{3}}\left(\Delta_h\gamma\cdot\Delta_h\gamma'-\Delta_h\zeta\cdot\Delta_h\zeta'\right)\,\d h\,\d x\\
&\ls F^3\iint_{\T^2}|h|^{-3}|\Delta_h(\gamma'-\zeta')\cdot\Delta_h\zeta'|
\left|\Delta_h\gamma\cdot\Delta_h\gamma'-\Delta_h\zeta\cdot\Delta_h\zeta'\right|\,\d h\,\d x\\
&\ls F^3\int_0^\pi
\frac{\left\|\Delta_h(\gamma'-\zeta')\right\|_{L^2}}{h^3}
(Mh^{1/2})(Mh^{3/2}\left\|\gamma'-\zeta'\right\|_{L^2})\,\d h\\
&= M^2F^3\left\|\gamma'-\zeta'\right\|_{L^2}
\int_0^\pi\frac{\left\|\Delta_h(\gamma'-\zeta')\right\|_{L^2}}h\,\d h\\
&\ls M^2F^3\left\|\gamma'-\zeta'\right\|_{L^2}
\Phi(\left\|\gamma'-\zeta'\right\|_{L^2}).
\end{aligned}
$$
In the second inequality, we used Cauchy--Schwarz, $\left\|\Delta_h\zeta'\right\|_{L^\infty}\leq M|h|^{1/2}$ and \eqref{eq:H-relative}; The last inequality is \eqref{eq:critical-linear-difference}.

Combining these bounds in \eqref{eq:normal-energy-decomposition}
and using $\Phi(s)\geq s$, we obtain
\begin{equation}\label{eq:endpoint-normal-energy}
|J_1|\ls M^4F^5\left\|\gamma-\zeta\right\|_{H^1}
\Phi(\left\|\gamma-\zeta\right\|_{H^1}).
\end{equation}

\textbf{(III) First-order estimate: Tangential term.}
Finally, \eqref{eq:tangential-decomposition}.
The first inequality in \eqref{eq:J21-bound} and \eqref{eq:sqg-gauge-derivative} give
$$
|J_{21}|\ls M^2F^3\left\|\gamma'-\zeta'\right\|_{L^2}^2.
$$
The first inequality in \eqref{eq:J22-bound} and
\eqref{eq:sqg-lambda-difference} give
$$
|J_{22}|\ls M^5F^6\left\|\gamma'-\zeta'\right\|_{L^2}\Phi(\left\|\gamma'-\zeta'\right\|_{L^2}).
$$

\textbf{(IV) Conclusion.}
Putting everything together, we have, for some absolute constant $C>0$,
$$
\frac{\d}{\d t}\left\|\gamma-\zeta\right\|_{H^1}^2
\leq CM^5F^6\left\|\gamma-\zeta\right\|_{H^1}^2
       \log\frac{M^2}{\left\|\gamma-\zeta\right\|_{H^1}^2},
$$
with the right-hand side defined as zero at $\gamma=\zeta$.
Integrating this inequality by Osgood's lemma gives \eqref{eq:SQG-H1-continuity}.
\end{proof}

\autoref{prop:SQG-H1-stability} gives uniqueness in the class
\eqref{eq:SQG-weak-time-class}. It also gives continuous dependence
in $C([0,T];H^1)$ under uniform $H^2$ and geometric bounds.

\subsection{Strong $H^2$ time continuity and patch solutions}
\label{sec:SQG-time-continuity}
Now we finally conclude the proof of \autoref{prop:SQG-contour-wellposedness} and \autoref{thm:SQG-wellposedness}.

\begin{proof}[Proof of \autoref{prop:SQG-contour-wellposedness}]
Let $\gamma$ be a solution constructed in \autoref{sec:SQG-local-existence}, and set $M=2\|\gamma_0''\|_{L^2}$.
The convergence of initial data \eqref{eq:SQG-initial-approximation} and weak lower semicontinuity in \eqref{eq:SQG-existence-bootstrap} give
\begin{equation}\label{eq:SQG-existence-limit-energy}
\|\gamma''(t)\|_{L^2}^2 \leq\|\gamma_0''\|_{L^2}^2+CF(\gamma_0)^3M^4t.
\end{equation}
Together with $C_wH^2$ and the strong continuity of the lower derivatives, this proves $\gamma(t)\to\gamma_0$ in $H^2$ at zero.

The key is, our existence and uniqueness results are valid for either velocity sign. That is, the solving procedure can be applied to the backward-in-time equation as well. 

Fix $t_0\in[0,T]$ and solve the forward and backward equations with initial data $\gamma(t_0)$. 
By uniqueness in \autoref{prop:SQG-H1-stability}, the resulting solutions agree with $\gamma(t_0+s)$ and $\gamma(t_0-s)$, respectively, for small $s\geq0$ whenever the corresponding time lies in $[0,T]$. 
Similarly, applying \eqref{eq:SQG-existence-limit-energy} in both directions and using \eqref{eq:SQG-local-geometric-control} at $t_0$ gives
$$
\|\gamma''(t)\|_{L^2}^2\leq\|\gamma''(t_0)\|_{L^2}^2
+CF(\gamma_0)^3\|\gamma_0''\|_{L^2}^4|t-t_0|
$$
for $t\in[0,T]$ sufficiently close to $t_0$. This norm continuity and weak $H^2$ continuity give strong $H^2$ continuity: $\gamma(t)\to\gamma(t_0)$ in $H^2$.

Together with \autoref{prop:SQG-H1-stability}, this proves
\autoref{prop:SQG-contour-wellposedness}. 
\end{proof}

\begin{proof}[Proof of \autoref{thm:SQG-wellposedness}]
Now we deal with the patch solutions. It's all about how to check the assumptions of \autoref{prop:patch-constant-speed-equivalence}.

For the existence part, we choose a constant-speed $H^2$ parametrization $\gamma_0$ of the initial boundary $\p\Omega_0$. The existence of \autoref{prop:SQG-contour-wellposedness} gives a solution $\gamma(t)$ of the contour equation with the uniform-in-time $H^2$ and arc–chord bounds, initiated by $\gamma_0$.
The converse direction of \autoref{prop:patch-constant-speed-equivalence} then addresses that the enclosed domains $\Omega(t)$ form a patch solution of the SQG equation.

Here, the strong $H^2$ time continuity also gives the boundary regularity required in the SQG part of \autoref{prop:patch-constant-speed-equivalence}.

Now we show uniqueness of the patch solution. 
Choose a constant-speed parametrization of the common initial boundary of two patch solutions $\Omega_1,\Omega_2$. 
The uniqueness part of \autoref{prop:SQG-contour-wellposedness} gives that the two parametrizations agree, hence $\Omega_1(t)=\Omega_2(t)$ for all $t\in[0,T]$.
\end{proof}

\subsection{Smooth approximation of constant-speed curves}
\label{app:constant-speed-approximation}

\begin{proof}[Proof of \autoref{lem:constant-speed-approximation}]
Assume $\gamma_0$ is positively oriented, simple, and has constant speed.
Let $\{ \rho_\varepsilon \}_{\varepsilon > 0} \subset C^\infty(\T)$ be an approximate identity. 
We first mollify $\gamma_0$ and then reparametrize it to be constant speed again. This leads to
$$
\begin{aligned}
s_\varepsilon(a)&:=-\pi+\frac{2\pi}{L(\rho_\varepsilon*\gamma_0)}\int_{-\pi}^{a}|(\rho_\varepsilon*\gamma_0)'(y)|\,\d y,\\
\gamma_\varepsilon(x)&:=(\rho_\varepsilon*\gamma_0)(r_\varepsilon(x)), \qquad r_\varepsilon:= s_\varepsilon^{-1}.
\end{aligned}
$$
where $L$ denotes length of the curve. For small $\varepsilon$, the estimate
$$
 |\Delta_h(\rho_\varepsilon*\gamma_0)|
 \geq |\Delta_h\gamma_0|
       -|h|\|(\rho_\varepsilon*\gamma_0-\gamma_0)'\|_{L^\infty}
 \geq\frac{|h|}{2F(\gamma_0)}
$$
keeps the mollified curve regular and simple. 
Thus $\gamma_\varepsilon$ is smooth, positively oriented, simple, and constant speed.
Since $\gamma_0$ has constant speed and $\rho_\varepsilon*\gamma_0\to\gamma_0$ in $H^2$,
$$
\begin{aligned}
(s_\varepsilon)'(a) &= \frac{2\pi}{L(\rho_\varepsilon*\gamma_0)} | (\rho_\varepsilon*\gamma_0)'(a) | \longrightarrow1\quad\hbox{in }C^0,\\
(s_\varepsilon)''(a) &= \frac{2\pi}{L(\rho_\varepsilon*\gamma_0)}
\frac{(\rho_\varepsilon*\gamma_0)' \cdot(\rho_\varepsilon*\gamma_0)''} {|(\rho_\varepsilon*\gamma_0)'|}
 \longrightarrow0\quad\hbox{in }L^2.
\end{aligned}
$$
Here, the numerator converges in $L^2$ to $\gamma_0'\cdot\gamma_0''=0$, while the denominators enjoy a common positive lower bound.
Thus, $r_\varepsilon$ converges to the identity in $C^1$, with its second derivative tending to zero in $L^2$. 
The chain rule now gives
$$
(\gamma_\varepsilon)''=\ ((\rho_\varepsilon*\gamma_0)''\circ r_\varepsilon)((r_\varepsilon)')^2+
((\rho_\varepsilon*\gamma_0)'\circ r_\varepsilon)r_\varepsilon''.
$$
The second term tends to zero in $L^2$. The first tends to $\gamma_0''$ by strong $L^2$ convergence of $(\rho_\varepsilon*\gamma_0)''$,
the uniform Jacobian bounds, and approximation of $\gamma_0''$ by continuous functions. 
This concludes the convergence in \eqref{eq:SQG-initial-approximation}.
Finally, the arc--chord bound $F(\gamma_\varepsilon) \le F(\gamma_0)$ for sufficiently small $\varepsilon$ follows from $|F(\gamma)^{-1}-F(\zeta)^{-1}| \leq\|\gamma'-\zeta'\|_{L^\infty}$.
\end{proof}

\subsection{Continuation of smooth SQG patches}
\label{app:SQG-smooth-continuation}

\begin{proof}[Proof of \autoref{lem:SQG-logarithmic-continuation}]
We adjust our proof of \autoref{prop:SQG-second-order}. Using $\dot H^1$ instead of $L^2$ inner product in \eqref{eq:SQG-second-order-decomposition}, we obtain
$$
\frac12\frac{\d}{\d t}\|\gamma'''\|_{L^2}^2=N_1+N_2+N_3+N_4,
$$
where
$$
\begin{aligned}
N_1&=-\iint_{\T^2}\frac{\gamma^{(4)}\cdot\Delta_h\gamma'''}
 {|\Delta_h\gamma|}\,\d x\,\d h,\\
N_2&=-2\iint_{\T^2}\gamma^{(4)}\cdot\Delta_h\gamma''
 (|\Delta_h\gamma|^{-1})'\,\d x\,\d h,\\
N_3&=-\iint_{\T^2}\gamma^{(4)}\cdot\Delta_h\gamma'
 (|\Delta_h\gamma|^{-1})''\,\d x\,\d h,\\
N_4&=-\int_\T\gamma^{(4)}\cdot(\lambda_\gamma\gamma')''\,\d x.
\end{aligned}
$$
The $H^5$ regularity justifies differentiation under the integral.

\textbf{(I) Estimate for $N_1$ and $N_2$.}
Symmetrization and estimate as in \eqref{eq:sqg-N1-calculation}-\eqref{eq:SQG-N1-bound} give
$$
\begin{aligned}
N_1&=-\frac12\iint_{\T^2}\frac{\Delta_h\gamma^{(4)}\cdot \Delta_h\gamma'''}{|\Delta_h\gamma|}\,\d x\,\d h\\
&=\frac14\iint_{\T^2}|\Delta_h\gamma'''|^2 (|\Delta_h\gamma|^{-1})'\,\d x\,\d h\\
&\ls \left\|\gamma''' \right\|_{L^2}^2 \int_\T \left\|(|\Delta_h\gamma|^{-1})' \right\|_{L^\infty}\,\d h \ls F^3 \left\|\gamma'' \right\|_{L^2}^2 \left\|\gamma''' \right\|_{L^2}^2,
\end{aligned}
$$
while
$$
\begin{aligned}
N_2&=-\iint_{\T^2}\Delta_h\gamma^{(4)}\cdot\Delta_h\gamma''
 (|\Delta_h\gamma|^{-1})'\,\d x\,\d h\\
&=\iint_{\T^2}|\Delta_h\gamma'''|^2
 (|\Delta_h\gamma|^{-1})'\,\d x\,\d h + \iint_{\T^2}(\Delta_h\gamma'''\cdot\Delta_h\gamma'')
 (|\Delta_h\gamma|^{-1})''\,\d x\,\d h .
\end{aligned}
$$
The first term is completely same as $N_1$.  
For the second integral, the kernel formula in \autoref{prop:SQG-second-order} gives
$$
(|\Delta_h\gamma|^{-1})''
=3\frac{(\Delta_h\gamma\cdot\Delta_h\gamma')^2}{|\Delta_h\gamma|^5}
 -\frac{|\Delta_h\gamma'|^2}{|\Delta_h\gamma|^3}
 -\frac{\Delta_h\gamma\cdot\Delta_h\gamma''}{|\Delta_h\gamma|^3},
$$
Follow the treatment of $N_3$ in \autoref{prop:SQG-second-order}, using \eqref{eq:arc-chord-lower}, and \eqref{eq:basic-difference-L2}, we get
$$
\begin{aligned}
\|(|\Delta_h\gamma|^{-1})''\|_{L^2}
&\ls \left\|\frac{|\Delta_h\gamma'|^2}{|\Delta_h\gamma|^3}\right\|_{L^2}
 +\left\|\frac{\Delta_h\gamma\cdot\Delta_h\gamma''}{|\Delta_h\gamma|^3}\right\|_{L^2}\\
&\ls F^3|h|^{-3} \left(\|\Delta_h\gamma'\|_{L^4}^2+\|\Delta_h\gamma\cdot\Delta_h\gamma''\|_{L^2}\right)\\
&\ls F^3\frac{\omega_{\gamma'}(|h|)}{h^2}\|\gamma''\|_{L^2}.
\end{aligned} 
$$
For the latter term we use the bound $\|\Delta_h\gamma\cdot\Delta_h\gamma''\|_{L^2}$ in \eqref{eq:SQG-Dh-Dh2-L2}. 
Thus the last integral for $N_2$ is bounded by
$$
CF^3\|\gamma'''\|_{L^2}\|\gamma''\|_{L^2}
 \int_\T\frac{\omega_{\gamma'}(|h|)\omega_{\gamma''}(|h|)}{h^2}\,\d h
\ls F^3\|\gamma''\|_{L^2}^2\|\gamma'''\|_{L^2}^2
$$
by Cauchy--Schwarz and \autoref{lem:SQG-square-modulus}.

\textbf{(II) Estimate for $N_3$.}
Using $(|\Delta_h\gamma|^{-1})'= -(\Delta_h\gamma\cdot\Delta_h\gamma')/|\Delta_h\gamma|^3$, we rewrite $N_3$ as
$$
\begin{aligned}
-N_3
&=-\iint_{\T^2}\gamma^{(4)}\cdot\Delta_h\gamma''
 (|\Delta_h\gamma|^{-1})'\,\d x\,\d h\\
&\quad+\iint_{\T^2}\gamma'''\cdot
 \left(\frac{\Delta_h\gamma'(\Delta_h\gamma\cdot\Delta_h\gamma')}
 {|\Delta_h\gamma|^3}\right)''\,\d x\,\d h.
\end{aligned}
$$
The first integral is just $N_2 / 2$ and has already been bounded.
Our focus is therefore on the second integral. By \autoref{lem:SQG-moser},
$$
N_5(h):=\left\|\left(\frac{
 \Delta_h\gamma'(\Delta_h\gamma\cdot\Delta_h\gamma')}
 {|\Delta_h\gamma|^3}\right)''\right\|_{L^2}
\ls N_{51}+N_{52}+N_{53},
$$
where
$$
\begin{aligned}
N_{51}&=\|\Delta_h\gamma'''\|_{L^2}
 \|\Delta_h\gamma\cdot\Delta_h\gamma'\|_{L^\infty}
 \||\Delta_h\gamma|^{-3}\|_{L^\infty},\\
N_{52}&=\|\Delta_h\gamma'\|_{L^\infty}
 \|(\Delta_h\gamma\cdot\Delta_h\gamma')''\|_{L^2}
 \||\Delta_h\gamma|^{-3}\|_{L^\infty},\\
N_{53}&=\|\Delta_h\gamma'\|_{L^\infty}
 \|\Delta_h\gamma\cdot\Delta_h\gamma'\|_{L^\infty}
 \|(|\Delta_h\gamma|^{-3})''\|_{L^2}.
\end{aligned}
$$
First we bound $N_{51}$.
By \eqref{eq:chord-tangent-identity} we have $\|\Delta_h\gamma\cdot\Delta_h\gamma'\|_{L^\infty} \leq |h|\omega_{\gamma'}(|h|)^2$, hence
$$
N_{51}\ls \|\gamma'''\|_{L^2}|h|\omega_{\gamma'}(|h|)^2\||\Delta_h\gamma|^{-3}\|_{L^\infty} \ls F^3\frac{\omega_{\gamma'}(|h|)^2}{h^2}\|\gamma'''\|_{L^2}.
$$
For $N_{52}$, we first bound $\|(\Delta_h\gamma\cdot\Delta_h\gamma')''\|_{L^2}$.
It suffices to take $h>0$.
We rewrite $\Delta_h\gamma\cdot\Delta_h\gamma'$ in \eqref{eq:chord-tangent-identity} as
$$
\Delta_h\gamma(x)\cdot\Delta_h\gamma'(x)
=\frac12\int_0^h\left(
 |\Delta_r\gamma'(x-h+r)|^2-|\Delta_r\gamma'(x)|^2\right)\,\d r,
$$
which allow us to differentiate in $x$ now. 
Use Minkowski's inequality, translation invariance of the $L^2$ norm, and \autoref{lem:SQG-moser} to obtain
$$
\begin{aligned}
\|(\Delta_h\gamma\cdot\Delta_h\gamma')''\|_{L^2}
&\leq\frac12\int_0^h\left( \|(|\Delta_r\gamma'(\cdot-h+r)|^2)''\|_{L^2} +\|(|\Delta_r\gamma'|^2)''\|_{L^2}\right)\,\d r\\
&=\int_0^h\|(|\Delta_r\gamma'|^2)''\|_{L^2}\,\d r\\
&\ls\int_0^h\|\Delta_r\gamma'\|_{L^\infty} \|\Delta_r\gamma'''\|_{L^2}\,\d r \ls h\omega_{\gamma'}(h)\|\gamma'''\|_{L^2}.
\end{aligned}
$$
Thus
$$
N_{52}\ls \omega_{\gamma'}(|h|) |h|\omega_{\gamma'}(|h|)\|\gamma'''\|_{L^2} \frac{F^3}{|h|^3}
=F^3\frac{\omega_{\gamma'}(|h|)^2}{h^2}\|\gamma'''\|_{L^2}.
$$

For $N_{53}$, we simply apply \autoref{lem:SQG-moser} and the bound of $\|(|\Delta_h\gamma|^{-1})''\|_{L^2}$ from $N_2$ above:
$$
\begin{aligned}
\|(|\Delta_h\gamma|^{-3})''\|_{L^2}
&\ls \||\Delta_h\gamma|^{-1}\|_{L^\infty}^2
 \|(|\Delta_h\gamma|^{-1})''\|_{L^2}\\
&\ls F^5|h|^{-4}\omega_{\gamma'}(|h|)\|\gamma''\|_{L^2}.
\end{aligned}
$$
Thus
$$
N_{53}\ls F^5\frac{\omega_{\gamma'}(|h|)^4}{|h|^3}\|\gamma''\|_{L^2}
\ls F^5M^2\frac{\omega_{\gamma'}(|h|)^2}{h^2}\|\gamma'''\|_{L^2},
$$
since $\omega_{\gamma'}(|h|)^2\leq |h|\|\gamma''\|_{L^2}^2$ and $\|\gamma''\|_{L^2}\leq\|\gamma'''\|_{L^2}$. By Cauchy--Schwarz and \autoref{lem:SQG-square-modulus},
$$
\begin{aligned}
|N_3+\tfrac12N_2|&\leq\|\gamma'''\|_{L^2}\int_\T N_5(h)\,\d h\\
&\ls F^5M^2\|\gamma'''\|_{L^2}^2 \int_\T\frac{\omega_{\gamma'}(|h|)^2}{h^2}\,\d h
\ls F^5M^4\|\gamma'''\|_{L^2}^2.
\end{aligned}
$$
Together with the bound for $N_2$, this gives $|N_3|\ls F^5M^4\|\gamma'''\|_{L^2}^2$.

\textbf{(III) Estimate for $N_4$.}
Integration by parts in $x$ and the product rule give
$$
\begin{aligned}
N_4
&=-\int_\T\gamma^{(4)}\cdot(\lambda_\gamma\gamma')''\,\d x =\int_\T\gamma'''\cdot(\lambda_\gamma\gamma')'''\,\d x\\
&=\int_\T\lambda_\gamma\gamma'''\cdot\gamma^{(4)}\,\d x
+3\int_\T\lambda_\gamma'|\gamma'''|^2\,\d x
+3\int_\T\lambda_\gamma''\gamma''\cdot\gamma'''\,\d x
+\int_\T\lambda_\gamma'''\gamma'\cdot\gamma'''\,\d x.
\end{aligned}
$$
The first term is
$$
\int_\T\lambda_\gamma\gamma'''\cdot\gamma^{(4)}\,\d x =\frac12\int_\T\lambda_\gamma(|\gamma'''|^2)'\,\d x =-\frac12\int_\T\lambda_\gamma'|\gamma'''|^2\,\d x.
$$
Also, the constant-speed condition gives $(|\gamma'|^2)''/2=|\gamma''|^2+\gamma'\cdot\gamma'''=0$, so the last term becomes
$$
\int_\T\lambda_\gamma'''\gamma'\cdot\gamma'''\,\d x
=-\int_\T\lambda_\gamma'''|\gamma''|^2\,\d x
=\int_\T\lambda_\gamma''(|\gamma''|^2)'\,\d x
=2\int_\T\lambda_\gamma''\gamma''\cdot\gamma'''\,\d x.
$$
Combining them we obtain
$$
N_4=\frac52\int_\T\lambda_\gamma'|\gamma'''|^2\,\d x+5\int_\T\lambda_\gamma''\gamma''\cdot\gamma'''\,\d x.
$$
By \autoref{lem:sqg-gauge-estimate}, $\|\lambda_\gamma'\|_{L^\infty} \ls F^3\|\gamma''\|_{L^2}^2$. Therefore the contribution of the first term is bounded by an absolute multiple of $F^3\|\gamma''\|_{L^2}^2\|\gamma'''\|_{L^2}^2$.
Note that $\lambda_\gamma''=-g_\gamma'$,  Cauchy--Schwarz gives
$$
\left|\int_\T\lambda_\gamma''\gamma''\cdot\gamma'''\,\d x\right|
=\left|\int_\T g_\gamma'\gamma''\cdot\gamma'''\,\d x\right|\leq\|\gamma''g_\gamma'\|_{L^2}\|\gamma'''\|_{L^2}.
$$
From \eqref{eq:g-definition}, note that $|\gamma'|$ is constant in our case, we have 
$$
|\gamma'|^2g_\gamma'=(\gamma'\cdot\mathcal N(\gamma)')' 
=\gamma''\cdot\mathcal N(\gamma)'+\gamma'\cdot\mathcal N(\gamma)''.
$$
Now we compute the derivatives of \eqref{eq:endpoint-normal-operator}:
$$
\begin{aligned}
\mathcal N(\gamma)'&=\left(\int_\T\frac{\Delta_h\gamma'}{|\Delta_h\gamma|}
 \,\d h\right)' =\int_\T\frac{\Delta_h\gamma''}{|\Delta_h\gamma|}\,\d h
 +\int_\T\Delta_h\gamma'(|\Delta_h\gamma|^{-1})'\,\d h,\\
\mathcal N(\gamma)''&=\left(\int_\T\frac{\Delta_h\gamma''}{|\Delta_h\gamma|}\,\d h +\int_\T\Delta_h\gamma'(|\Delta_h\gamma|^{-1})'\,\d h\right)'\\
&=\int_\T\frac{\Delta_h\gamma'''}{|\Delta_h\gamma|}\,\d h
+2\int_\T\Delta_h\gamma''(|\Delta_h\gamma|^{-1})'\,\d h
+\int_\T\Delta_h\gamma'(|\Delta_h\gamma|^{-1})''\,\d h.
\end{aligned}
$$
Combining these together gives $|\gamma'|^2g_\gamma'=N_{61}+N_{62}+N_{63}$, where
$$
\begin{aligned}
N_{61}&=\int_\T\frac{\gamma''\cdot\Delta_h\gamma'' +\gamma'\cdot\Delta_h\gamma'''}{|\Delta_h\gamma|}\,\d h,\\
N_{62}&=\int_\T (\gamma''\cdot\Delta_h\gamma'+2\gamma'\cdot\Delta_h\gamma'') (|\Delta_h\gamma|^{-1})'\,\d h,\\
N_{63}&=\int_\T(\gamma'\cdot\Delta_h\gamma') (|\Delta_h\gamma|^{-1})''\,\d h.
\end{aligned}
$$
For $N_{61}$, $\gamma'\cdot\gamma''=0$ gives
$
\gamma'(x)\cdot\Delta_h\gamma''(x)=-\gamma'(x)\cdot\gamma''(x-h)=-\Delta_h\gamma'(x)\cdot\gamma''(x-h).
$
Differentiating this yields
$$
\gamma''\cdot\Delta_h\gamma''+\gamma'\cdot\Delta_h\gamma'''
=-\Delta_h\gamma''\cdot\gamma''(x-h) -\Delta_h\gamma'\cdot\gamma'''(x-h).
$$
Using Cauchy--Schwarz in $h$, and \autoref{lem:SQG-square-modulus}, we reach the estimate
$$
\begin{aligned}
\|N_{61}\|_{L^\infty}
&\ls F\left( \left\|\frac{\omega_{\gamma''}(|h|)}{|h|}\right\|_{L^2_h}
 \|\gamma''\|_{L^2} +\left\|\frac{\omega_{\gamma'}(|h|)}{|h|}\right\|_{L^2_h} \|\gamma'''\|_{L^2}\right)
 \ls FM\|\gamma'''\|_{L^2}.
\end{aligned}
$$
For $N_{62}$, we again use
$\gamma'\cdot\Delta_h\gamma''=-\Delta_h\gamma'\cdot\gamma''(x-h)$ to get
$$
\gamma''\cdot\Delta_h\gamma'+2\gamma'\cdot\Delta_h\gamma''
=\Delta_h\gamma'\cdot(\gamma''(x)-2\gamma''(x-h)).
$$
By Minkowski's inequality and \eqref{eq:SQG-existence-kernel-x},
$$
\begin{aligned}
\|N_{62}\|_{L^2}
&=\left\|\int_\T\Delta_h\gamma'\cdot (\gamma''(x)-2\gamma''(x-h))(|\Delta_h\gamma|^{-1})' \,\d h\right\|_{L^2}\\
&\ls \int_\T\omega_{\gamma'}(|h|)\|\gamma''\|_{L^2} \|(|\Delta_h\gamma|^{-1})'\|_{L^\infty}\,\d h\\
&\ls F^3\|\gamma''\|_{L^2}\int_\T \frac{\omega_{\gamma'}(|h|)^3}{h^2}\,\d h\ls F^3M^4,
\end{aligned}
$$
where the last step uses $\omega_{\gamma'}(|h|)\ls M$ and \autoref{lem:SQG-square-modulus}.
For $N_{63}$, constant speed gives
$$
\gamma'\cdot\Delta_h\gamma'=\tfrac12|\Delta_h\gamma'|^2.
$$
It's been shown that $\|(|\Delta_h\gamma|^{-1})''\|_{L^2} \ls F^3h^{-2}\omega_{\gamma'}(|h|)\|\gamma''\|_{L^2}$ in the estimate for $N_2$. Hence, by Minkowski's inequality and \autoref{lem:SQG-square-modulus} again, 
$$
\begin{aligned}
\|N_{63}\|_{L^2}
&=\frac12\left\|\int_\T|\Delta_h\gamma'|^2 (|\Delta_h\gamma|^{-1})''\,\d h\right\|_{L^2}\\
&\ls \int_\T\omega_{\gamma'}(|h|)^2 \|(|\Delta_h\gamma|^{-1})''\|_{L^2}\,\d h\\
&\ls F^3\|\gamma''\|_{L^2}\int_\T \frac{\omega_{\gamma'}(|h|)^3}{h^2}\,\d h\ls F^3M^4.
\end{aligned}
$$

Combining the bounds and using \eqref{eq:speed-arc-chord}, we obtain
$$
\begin{aligned}
\|\gamma''g_\gamma'\|_{L^2}
&\leq |\gamma'|^{-2}\left(\|\gamma''\|_{L^2}\|N_{61}\|_{L^\infty}
 +\|\gamma''\|_{L^\infty}(\|N_{62}\|_{L^2}+\|N_{63}\|_{L^2})\right)\\
&\ls F^3M^2\|\gamma'''\|_{L^2}+F^5M^4\|\gamma''\|_{L^\infty}
 \ls F^5M^4\|\gamma'''\|_{L^2}.
\end{aligned}
$$

We finally conclude that, by the estimates for $N_1$--$N_4$ all above,
$$
\left|\frac{\d}{\d t}\|\gamma'''\|_{L^2}^2\right| \ls F^5M^4\|\gamma'''\|_{L^2}^2.
$$
Applying Gr\"onwall's inequality completes the proof.
\end{proof}

\appendix

\section{Justification of contour dynamics}\label{sec:patch-contour}

In this section we deal with the equivalence of the transport equation and the contour dynamics for patch solutions.
When we deal with the vortex patches for 2D Euler, the velocity field is bounded and log-Lipschitz, and the equivalence is classical. 
For gSQG with $0<\beta<1$, the velocity field is still H\"older continuous, 
while the SQG case $\beta=1$ only offers $u\in \text{BMO}$, in which case one may even cannot define the velocity on the boundary of the patch. 
Indeed H\"older continuity on velocity field is not enough to directly guarantee the uniqueness of particle trajectories
---as shown by the famous example $\dot{x}=\sqrt{|x|}$---
and hence a particle on the boundary initially may leave the boundary at any time. 
Before we apply the derivations of the contour dynamics in the beginning of \autoref{sec:contour-dynamics}, we need to overcome this difficulty.

We first prove the equivalence between weak transport and the boundary equation for a general velocity field, then derive the constant-speed gauge. 
In \autoref{sec:patch-contour-application} we verify the hypotheses for gSQG and SQG and prove \autoref{prop:patch-constant-speed-equivalence}.

\subsection{Weak transport and the boundary equation}\label{sec:patch-parametrization}
Let's review the general concept of transporting by a divergence-free velocity field $v$, especially when $v$ is not Lipschitz. 
As we discussed above, the lack of uniqueness of the flow map $X(t,x)$ may cause the contour dynamics invalid.
The theory of DiPerna--Lions \cite{DiPerna1989} and Ambrosio \cite{Ambrosio2004} shows that the flow map is unique in a measure-theoretic sense, with barely $W^{1,p}$ or BV velocity fields. 
This however doesn't help in our case too, since the boundary of patch is Lebesgue null. 
To make contour dynamics work again, we need to impose some additional structure on the velocity field. If everything is good we write the Green's formula
$$
\int_{\Omega(t)}v\cdot\nabla\varphi\,\d y
=\int_{\p\Omega(t)}\varphi(y)v(t,y)\cdot n(t,y)\,\d\sigma(y).
$$
But here on the boundary, only the normal component of $v$ matters. This leaves an opportunity to loosen the regularity requirement on $v$ if we can control $v\cdot n$ on the boundary. In the following, we denote by $V$ for some vector field which satisfies $V\cdot n=v\cdot n$ formally.

\begin{proposition}\label{prop:weak-transport-boundary-equivalence}
Let $\Omega(t)\subset\R^2$, $t\in[0,T]$, be a family of bounded simply connected domains with $\ch_{\Omega(t)}\in C_w([0,T];L^2(\R^2))$. 
Assume that for each $t\in[0,T]$, there exists a positively oriented parametrization of $\p\Omega(t)$, 
and the closure of this family in $C^1(\T;\R^2)$ is compact and consists of regular simple closed curves.

Let $v\in L^1_{\loc}([0,T]\times\R^2;\R^2)$ be divergence-free in the sense of distributions. 
Let $V(t,y)$ be a bounded measurable vector field on $y\in\p\Omega(t)$ with respect to $\d t\,\d\sigma(y)$, 
where $\d\sigma$ is arclength measure on $\p\Omega(t)$. 
Write $n(t,y)$ for the outward unit normal and assume that, for almost every $t\in(0,T)$
and every $\varphi\in C_c^\infty(\R^2)$,
\begin{equation}\label{eq:transport-normal-flux}
\int_{\Omega(t)}v\cdot\nabla\varphi\,\d y
=\int_{\p\Omega(t)}\varphi(y)V(t,y)\cdot n(t,y)\,\d\sigma(y).
\end{equation}

Then, the following statements are equivalent.
\begin{enumerate}[(i)]
\item The Patch solution solves the weak transport equation: 
For every $\varphi\in C_c^\infty([0,T)\times\R^2)$,
\begin{equation*}
\int_0^T\!\int_{\Omega(t)} (\p_t\varphi+v\cdot\nabla\varphi)\,\d y\,\d t+\int_{\Omega(0)}\varphi(0,y)\,\d y=0.
\end{equation*}

\item The boundary admits a parametrization $\gamma\in C([0,T];C^1(\T;\R^2))$ with derivative in time $\gamma_t\in L^\infty((0,T)\times\T)$, such that $\gamma(t,\cdot)$ is a $C^1$ embedding onto $\p\Omega(t)$ for every $t$. 
Every parametrization of this class satisfies, almost everywhere,
\begin{equation}\label{eq:transport-boundary-equation}
\gamma_t(t,x)\cdot\n_\gamma(t,x)=V(t,\gamma(t,x))\cdot\n_\gamma(t,x).
\end{equation}
\end{enumerate}
\end{proposition}

\begin{proof}
Assume (i). We derive the motion of the boundary in a fixed
transverse chart and then read off its normal velocity in any
parametrization of the class in (ii).

\textbf{(I) Continuity of the boundaries.}
The compactness assumption on $\p\Omega(t)$ asserts that $\Omega(t)$ is bounded uniformly in $t\in [0,T]$. 
Hence the weak $L^2$ continuity of $\ch_{\Omega(t)}$ gives that $|\Omega(t)|$ is continuous in time,\footnote{ Indeed the incompressibility of $v$ gives that $|\Omega(t)|=|\Omega(0)|$ for all $t\in[0,T]$.} 
which asserts that $\ch_{\Omega(t)}$ is actually strongly $L^2$ continuous in time. 

Now we need the continuity of the boundaries in Hausdorff distance, which allows us to choose a tubular neighborhood of $\p\Omega(t_0)$ such that it contains $\p\Omega(t)$ for all $t$ near $t_0$.
\begin{lemma}\label{lem:uniform-cone-compactness}\cite[Proposition~3.2 and Theorem~3.4]{Dalphin2018}
Let $\{ D_j \}$ be a sequence of open sets satisfying the following two conditions.
\begin{enumerate}[(a)]
\item Each $D_j$ is nonempty and contained in a bounded $B\subset \R^2$.
\item $\{ D_j \}$ satisfies Chenais's cone condition \cite{Chenais1975}: 
There is an $\eps\in(0,\pi/2)$ such that, for each $j$ and $p\in\p D_j$, there is a unit vector $\xi$ for which
\begin{equation*}
y+z\in D_j \quad\text{whenever}\quad y\in D_j\cap B_\eps(p),
\quad 0<|z|<\eps,\quad \text{ and } z\cdot\xi>|z|\cos\eps.
\end{equation*}
\end{enumerate}
Then there are a nonempty open set $D\subset B$ satisfying the same cone condition, and a subsequence $\{ D_{j_k} \}$ such that $|D_{j_k}\Delta D|\to 0$ and $\p D_{j_k}\to\p D$ in Hausdorff distance.
\end{lemma}

Fix any sequence $t_j\to t_0$ and set $D_j=\Omega(t_j)$. The two conditions of \autoref{lem:uniform-cone-compactness} can be easily verified from the $C^1$ compactness of the boundary parametrizations.
Now we can take a subsequence $\{ t_{j_k} \}$ such that $\p\Omega(t_{j_k})\to\p D$ in Hausdorff distance for some open set $D$. Note that $\ch_{\Omega(t_{j_k})}\to\ch_{\Omega(t_0)}$ in $L^2$, we have $\ch_D=\ch_{\Omega(t_0)}$ almost everywhere.
It follows from $C^1$ compactness again, that $D$ is simply connected and its boundary is a simple closed curve, thus $D=\Omega(t_0)$. By the arbitrariness of the sequence $\{ t_j \}$, we conclude that as $t\to t_0$, we have $d_{\text{H}}(\p\Omega(t),\p\Omega(t_0))\to 0$.

\textbf{(II) Coordinates near the boundary.}
Fix $t_0$ and let $\gamma_0$ be any parametrization of $\p\Omega(t_0)$.
Choose a smooth unit transverse vector field $w(x)$ with $w(x)\cdot\n_{\gamma_0}(x)>0$.
The tubular neighborhood theorem for $C^1$ curves gives some $\delta>0$ and a $C^1$ diffeomorphism $\Psi$ onto $U$, defined by
$$
\Psi(x,r)=\gamma_0(x)+r w(x),\qquad (x,r)\in\T\times(-\delta,\delta),\qquad U=\Psi\bigl(\T\times(-\delta,\delta)\bigr).
$$
For $t$ sufficiently close to $t_0$, continuity in Hausdorff distance gives
$$
\p\Omega(t)\subset\Psi\bigl(\T\times(-\delta/2,\delta/2)\bigr).
$$

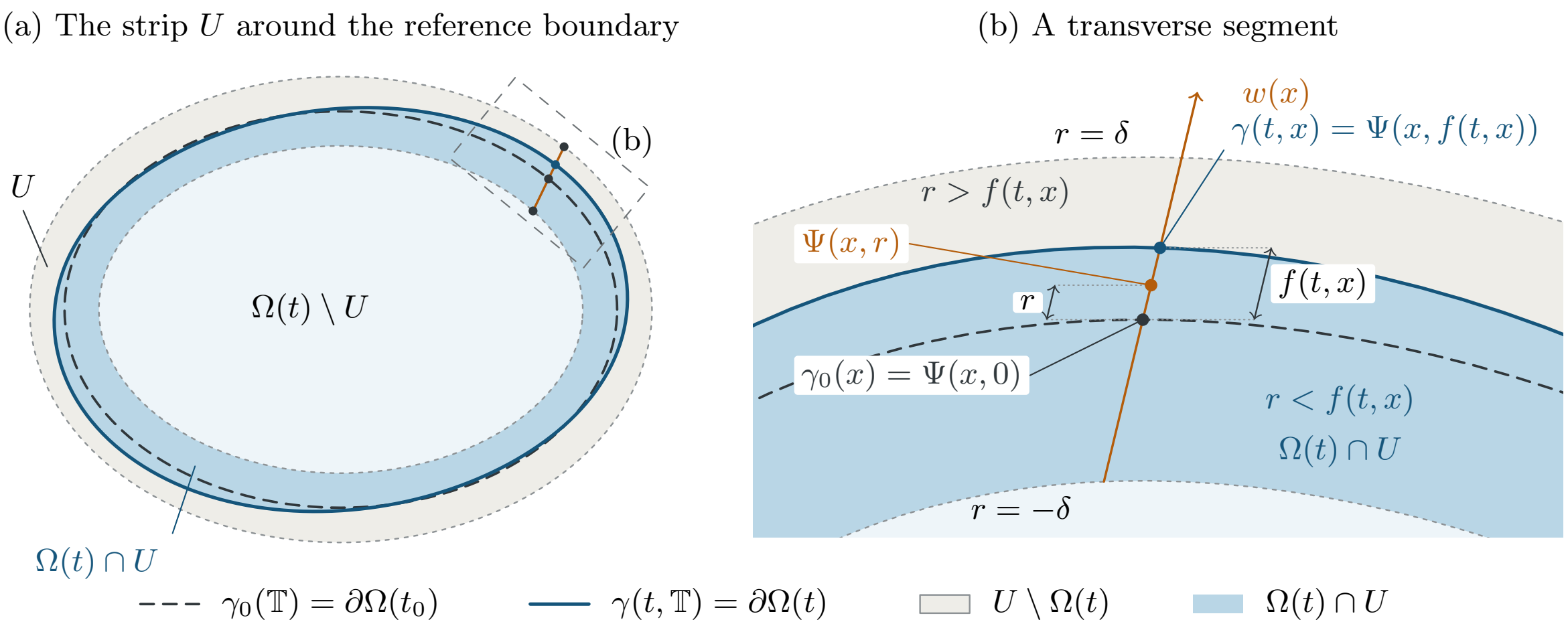
\begin{figure}[htbp]
\centering
\tikzsetnextfilename{boundary_coordinates}
\begin{tikzpicture}[font=\footnotesize, line cap=round, line join=round]
  \definecolor{referenceink}{HTML}{30373B}
  \definecolor{stripedge}{HTML}{8D9295}
  \definecolor{stripoutside}{HTML}{EEEDE8}
  \definecolor{patchblue}{HTML}{15557B}
  \definecolor{stripinside}{HTML}{B9D6E6}
  \definecolor{patchinterior}{HTML}{EEF5F9}
  \definecolor{transverseink}{HTML}{B45C0A}
  \tikzset{
    strip edge/.style={draw=stripedge, line width=.45pt, dash pattern=on 1pt off 1.7pt},
    reference boundary/.style={draw=referenceink, line width=.7pt, dash pattern=on 3.8pt off 2.6pt},
    moving boundary/.style={draw=patchblue, line width=.95pt},
    transverse segment/.style={draw=transverseink, line width=.65pt},
    zoom frame/.style={draw=referenceink!70, line width=.4pt, dashed},
    coordinate guide/.style={draw=stripedge, densely dotted, line width=.35pt},
    coordinate measure/.style={<->, draw=referenceink, line width=.45pt},
    point marker/.style={circle, minimum size=2.5pt, inner sep=0pt},
    point label/.style={fill=white, rounded corners=1pt, inner xsep=2pt, inner ysep=1pt},
    leader/.style={line width=.4pt, shorten >=1.5pt}
  }

  \def\referenceRadiusX{2.55}
  \def\referenceRadiusY{1.83}
  \def\stripHalfWidth{0.32}
  \def\transverseSlope{0.24}
  \tikzset{declare function={
    referenceSpeed(\a)=sqrt((\referenceRadiusX*sin(\a))^2+(\referenceRadiusY*cos(\a))^2);
    transverseX(\a)=(\referenceRadiusY*cos(\a)-\transverseSlope*\referenceRadiusX*sin(\a))/referenceSpeed(\a);
    transverseY(\a)=(\referenceRadiusX*sin(\a)+\transverseSlope*\referenceRadiusY*cos(\a))/referenceSpeed(\a);
    stripX(\a,\r)=\referenceRadiusX*cos(\a)+\r*transverseX(\a);
    stripY(\a,\r)=\referenceRadiusY*sin(\a)+\r*transverseY(\a);
    patchHeight(\a)=0.06+0.085*sin(2*\a+0.4*180/pi);
  }}
  \def\stripBoundaryPath#1{%
    plot[domain=0:360, samples=181, variable=\plotAngle]
      ({stripX(\plotAngle,#1)},{stripY(\plotAngle,#1)}) -- cycle%
  }
  \def\drawBoundaryStrip{%
    \filldraw[strip edge, fill=stripoutside] \stripBoundaryPath{\stripHalfWidth};
    \fill[stripinside] \stripBoundaryPath{patchHeight(\plotAngle)};
    \filldraw[strip edge, fill=patchinterior] \stripBoundaryPath{-\stripHalfWidth};
    \draw[reference boundary] \stripBoundaryPath{0};
    \draw[moving boundary] \stripBoundaryPath{patchHeight(\plotAngle)};
  }
  \pgfmathsetmacro{\zoomAngle}{0.72*180/pi}
  \pgfmathsetmacro{\zoomX}{stripX(\zoomAngle,0)}
  \pgfmathsetmacro{\zoomY}{stripY(\zoomAngle,0)}
  \pgfmathsetmacro{\tangentX}{-\referenceRadiusX*sin(\zoomAngle)/referenceSpeed(\zoomAngle)}
  \pgfmathsetmacro{\tangentY}{\referenceRadiusY*cos(\zoomAngle)/referenceSpeed(\zoomAngle)}
  \pgfmathsetmacro{\normalX}{\tangentY}
  \pgfmathsetmacro{\normalY}{-\tangentX}
  \pgfmathsetmacro{\boundaryHeight}{patchHeight(\zoomAngle)}
  \pgfmathsetmacro{\sampleHeight}{0.48*\boundaryHeight}
  \def\zoomLowerLeft{(-.77,-.43)}
  \def\zoomUpperRight{(.83,.53)}

  \node at (3.48,5.70) {(a) The strip $U$ around the reference boundary};
  \begin{scope}[shift={(3.48,2.91)}, scale=1.07]
    \drawBoundaryStrip
    \node at (-.28,0) {$\Omega(t)\setminus U$};
    \node (strip-label) at (-2.93,1.14) {$U$};
    \draw[leader, referenceink] (strip-label.south) -- (-2.69,.35);
    \node[patchblue] (inside-strip-label) at (-2.25,-2.34) {$\Omega(t)\cap U$};
    \draw[leader, patchblue] (inside-strip-label.north east) -- (-1.33,-1.40);

    \begin{scope}[cm={\tangentX,\tangentY,\normalX,\normalY,(\zoomX,\zoomY)}]
      \draw[zoom frame] \zoomLowerLeft rectangle \zoomUpperRight;
      \draw[transverse segment]
        ({-\transverseSlope*\stripHalfWidth},-\stripHalfWidth)
        -- ({\transverseSlope*\stripHalfWidth},\stripHalfWidth);
      \foreach \height in {-\stripHalfWidth,0,\stripHalfWidth}
        \node[point marker, fill=referenceink] at ({\transverseSlope*\height},{\height}) {};
      \node[point marker, fill=patchblue] at ({\transverseSlope*\boundaryHeight},\boundaryHeight) {};
    \end{scope}
    \node at ({\zoomX+.77},{\zoomY+.34}) {(b)};
  \end{scope}

  \node at (11.55,5.70) {(b) A transverse segment};
  \begin{scope}[shift={(11.40,2.81)}, scale=5]
    \begin{scope}
      \clip \zoomLowerLeft rectangle \zoomUpperRight;
      \begin{scope}[cm={\tangentX,\normalX,\tangentY,\normalY,(0,0)},
                    shift={(-\zoomX,-\zoomY)}]
        \drawBoundaryStrip
      \end{scope}
    \end{scope}
    \coordinate (reference-point) at (0,0);
    \coordinate (sample-point) at ({\transverseSlope*\sampleHeight},\sampleHeight);
    \coordinate (moving-point) at ({\transverseSlope*\boundaryHeight},\boundaryHeight);
    \draw[transverse segment, ->] ({-\transverseSlope*\stripHalfWidth},-\stripHalfWidth)
      -- ({.45*\transverseSlope},.45);
    \node[transverseink, anchor=west] at (.17,.445) {$w(x)$};

    \draw[coordinate guide] (-.19,0) -- (.22,0);
    \draw[coordinate guide] (sample-point) -- ({-.19+\transverseSlope*\sampleHeight},\sampleHeight);
    \draw[coordinate guide] (moving-point) -- ({.22+\transverseSlope*\boundaryHeight},\boundaryHeight);
    \draw[coordinate measure] (-.19,0) -- ({-.19+\transverseSlope*\sampleHeight},\sampleHeight)
      node[midway, left=2pt, point label] {$r$};
    \draw[coordinate measure] (.22,0) -- ({.22+\transverseSlope*\boundaryHeight},\boundaryHeight)
      node[midway, right=2pt, point label] {$f(t,x)$};

    \node[point label, anchor=west, text=referenceink] (reference-label) at (-.69,-.11)
      {$\gamma_0(x)=\Psi(x,0)$};
    \node[point label, anchor=west, text=transverseink] (sample-label) at (-.69,.15)
      {$\Psi(x,r)$};
    \node[point label, anchor=west, text=patchblue] (moving-label) at (.16,.375)
      {$\gamma(t,x)=\Psi(x,f(t,x))$};
    \draw[leader, referenceink] (reference-label.east) -- (reference-point);
    \draw[leader, transverseink] (sample-label.east) -- (sample-point);
    \draw[leader, patchblue] (moving-label.south west) -- (moving-point);
    \node[point marker, minimum size=3.5pt, fill=referenceink] at (reference-point) {};
    \node[point marker, minimum size=3.5pt, fill=transverseink] at (sample-point) {};
    \node[point marker, minimum size=3.5pt, fill=patchblue] at (moving-point) {};

    \node at (-.10,.37) {$r=\delta$};
    \node at (-.24,-.375) {$r=-\delta$};
    \node[referenceink] at (-.29,.25) {$r>f(t,x)$};
    \node[patchblue] at (.39,-.16) {$r<f(t,x)$};
    \node[patchblue] at (.39,-.26) {$\Omega(t)\cap U$};
  \end{scope}

  \draw[reference boundary] (1.50,0) -- (2.08,0);
  \node[anchor=west] at (2.17,0) {$\gamma_0(\T)=\p\Omega(t_0)$};
  \draw[moving boundary] (5.35,0) -- (5.93,0);
  \node[anchor=west] at (6.02,0) {$\gamma(t,\T)=\p\Omega(t)$};
  \filldraw[draw=stripedge, fill=stripoutside, line width=.45pt]
    (9.20,-.09) rectangle (9.69,.09);
  \node[anchor=west] at (9.78,0) {$U\setminus\Omega(t)$};
  \fill[stripinside] (11.90,-.09) rectangle (12.39,.09);
  \node[anchor=west] at (12.48,0) {$\Omega(t)\cap U$};
\end{tikzpicture}
\caption{The strip $U$ and a transverse segment. The dashed curve is the reference boundary $\p\Omega(t_0)$ and the solid curve is $\p\Omega(t)$. 
The blue part of the strip lies inside $\Omega(t)$, and the moving boundary meets this segment at $r=f(t,x)$.}
\label{fig:boundary-coordinates}
\end{figure}

Now we consider a $C^1$ simple closed curve $\gamma$ contained in the strip $U$, which parametrizes the boundary of a simply connected domain $\Omega$. 
If $\gamma$ is sufficiently close to $\gamma_0$ in $C^1$, then $\gamma$ can be written as a graph over $\gamma_0$, as a consequence of the implicit function theorem.
Indeed, we solve $(a,r)$ from the equation $$
F(\gamma,x;a,r):=\gamma(a)-\Psi(x,r)=0
$$
near $(\gamma, x,a,r)=(\gamma_0, x_*, x_*,0)$. Since $(\p_a F,\p_r F)(\gamma_0, x_*, x_*,0)=(\gamma_0'(x_*),-w(x_*))$ is invertible, the implicit function theorem gives a unique $C^1$ mapping $(\gamma,x)\mapsto(a,r)$ in a neighborhood of $(\gamma_0, x_*)$ such that $F(\gamma,x;a,r)=0$. The compactness of $\T$ allows us to obtain a unique $C^1$ mapping $(\gamma,x)\mapsto(a,r)$ for all $x\in\T$. We omit some functional setup details, but the argument here is standard.

For $t$ sufficiently close to $t_0$, the above procedure gives a unique $f(t,\cdot)$ such that
$$
\Omega(t)\cap U=\{\Psi(x,r):x\in\T,\ -\delta<r<f(t,x)\}.
$$
Clearly we have $f(t_0,x)=0$ and $|f(t,x)|<\delta/2$. The $C^1$ mapping above gives $f(t,\cdot)\in C^1(\T)$ and is continuous in $t$ in the $C^1$ topology.
We define the parametrization of $\p\Omega(t)$ by $\gamma(t,x):=\Psi(x,f(t,x))$.

\textbf{(III) Deriving the Contour Dynamics.}
Take $\rho\in C_c^\infty(0,T)$ supported sufficiently close to $t_0$, and $\varphi\in C_c^\infty(\R^2)$. 
Testing (i) with $\rho(t)\varphi(y)$ and using \eqref{eq:transport-normal-flux}, we obtain
\begin{equation*}
-\int_0^T\rho'(t)\int_{\Omega(t)}\varphi(y)\,\d y\,\d t
=\int_0^T\rho(t)\int_{\p\Omega(t)}\varphi(y)V(t,y)\cdot n(t,y)\,\d\sigma(y)\,\d t.
\end{equation*}
Clearly this identity can be extended to every $\varphi\in C_c(\R^2)$,
since the areas and boundary lengths of $\Omega(t)$ are uniformly bounded and $V$ is bounded. 

We choose $\varphi\in C_c(\R^2)$ so that the change of its integral on $\Omega(t)$ can be written in the graph height $f$. 
Denote $J(x,r)=|\det(\Psi_x(x,r),w(x))|$ the Jacobian of coordinates $(x,r)\mapsto\Psi(x,r)$.
Choose $\eta\in C_c^\infty(-\delta,\delta)$ equal to one on $[-\delta/2,\delta/2]$. 
For an arbitrary smooth periodic function $\psi$, we can now define
\begin{equation*}
\varphi(\Psi(x,r))=\frac{\psi(x)\eta(r)}{J(x,r)}
\end{equation*}
and extend $\varphi$ by zero outside $U$. 
Since $J$ is positive and continuous, this defines a function in $C_c(\R^2)$.
Changing variables $y=\Psi(x,r)$ and using $f(t_0,x)=0$, we obtain
\begin{equation*}
\begin{aligned}
\int_{\Omega(t)}\varphi(y)\,\d y-\int_{\Omega(t_0)}\varphi(y)\,\d y
&=\int_\T\psi(x)\int_0^{f(t,x)}\eta(r)\,\d r\,\d x\\
&=\int_\T\psi(x)f(t,x)\,\d x,
\end{aligned}
\end{equation*}
where the last equality follows from $|f(t,x)|<\delta/2$.

Now we compute the boundary integral in the right hand side. 
Since $\gamma(t,\cdot)\to\gamma_0$ in $C^1$ as $t\to t_0$, $w\cdot\n_\gamma$ has a uniform positive lower bound near $t_0$.
Taking determinants with $w(x)$ in $\gamma'(t,x)=\Psi_x(x,f(t,x))+f_x(t,x)w(x)$ and using $w(x)\cdot\n_\gamma(t,x)>0$ gives
\begin{equation*}
\begin{aligned}
J(x,f(t,x))&=|\det(\Psi_x(x,f(t,x)),w(x))|\\ 
&=|\det(\gamma'(t,x),w(x))|=|\gamma'(t,x)|\,w(x)\cdot\n_\gamma(t,x). \\
\end{aligned}
\end{equation*}
Changing variables $y=\gamma(t,x)$, using $\eta(f(t,x))=1$ and $\d\sigma(y)=|\gamma'(t,x)|\,\d x$, we obtain
\begin{equation*}
\int_{\p\Omega(t)}\varphi(y)V(t,y)\cdot n(t,y)\,\d\sigma(y)
=\int_\T\psi(x)\frac{V(t,\gamma(t,x))\cdot\n_\gamma(t,x)}
{w(x)\cdot\n_\gamma(t,x)}\,\d x.
\end{equation*}
This gives in the distributional sense in $t$,
\begin{equation*}
\frac{\d}{\d t}\int_\T\psi(x)f(t,x)\,\d x
=\int_\T\psi(x)\frac{V(t,\gamma(t,x))\cdot\n_\gamma(t,x)}{w(x)\cdot\n_\gamma(t,x)}\,\d x.
\end{equation*}
Since $\psi$ is arbitrary and $w\cdot\n_\gamma$ is bounded away from zero, the graph height has bounded time derivative
\begin{equation}\label{eq:weak-normal-law}
f_t(t,x)=\frac{V(t,\gamma(t,x))\cdot\n_\gamma(t,x)}{w(x)\cdot\n_\gamma(t,x)}.
\end{equation}
Since $\gamma(t,x)=\gamma_0(x)+f(t,x)w(x)$, we have
\begin{equation*}
\gamma_t(t,x)=f_t(t,x)w(x),\qquad
\gamma_t(t,x)\cdot\n_\gamma(t,x)=V(t,\gamma(t,x))\cdot\n_\gamma(t,x).
\end{equation*}

\textbf{(IV) Normal velocity in any parametrization.}
Now we complete the proof of (ii). Let $\tilde{\gamma}$ be a parametrization of the class as stated in (ii). 
Let $\gamma$ be the parametrization we choose above and write
$$
\tilde{\gamma}(t,x)=\gamma(t,a(t,x)),
$$
where $a(t,\cdot)=\gamma(t,\cdot)^{-1}\circ\tilde{\gamma}(t,\cdot)$.
The chain rule gives, almost everywhere\footnote{
A main technical point here is that $a(t,x)$ is only Lipschitz in $t$ and may not have a classical time derivative.
The chain rule requires $a_t(t,x)$ and $\gamma_t(t,a(t,x))$ to exist for almost every $(t,x)$. We verify this by showing that $(t,x)\mapsto(t,a(t,x))$ is locally bi-Lipschitz, so the null exceptional set for $\gamma_t$ remains null.
For the time estimate, we note that $a(t,x)$ is the first coordinate of $\Psi^{-1}(\tilde{\gamma}(t,x))$. Since $\Psi^{-1}$ is $C^1$ and $\tilde{\gamma}$ is Lipschitz in time in $C(\T)$, $a$ is locally Lipschitz in time in the uniform distance on $C(\T)$. 
The $C^1$ coordinate and $\tilde{\gamma}\in C_tC^1$ give joint continuity of $a'$. Moreover, $a'$ is positively lower bounded since both parametrizations are regular and positively oriented.
The change of variables $(t,x)\mapsto(t,a(t,x))$ is therefore locally bi-Lipschitz, hence the chain rule above holds almost everywhere by Rademacher's theorem.
},
$$
\tilde{\gamma}_t(t,x)=\gamma_t(t,a(t,x))
+\gamma'(t,a(t,x))a_t(t,x).
$$
Since the last term is tangent to the boundary, we obtain \eqref{eq:transport-boundary-equation} for $\tilde{\gamma}$.

\textbf{(V) From local to global.}
We have shown that for each $t_0$, there is a neighborhood $I(t_0)$ on which a parametrization $\gamma$ exists, satisfying \eqref{eq:transport-boundary-equation}. 
Now we can cover $[0,T]$ by finitely many of these neighborhoods and splice these local parametrizations together to obtain a global one. This completes the proof of (ii).

\textbf{(VI) Recovering the weak equation.}
For the converse, take one parametrization satisfying
\eqref{eq:transport-boundary-equation} with the regularity in (ii).
Take $\varphi\in C_c^\infty([0,T)\times\R^2)$. 
Since $\gamma\in C_tC^1$ and $\gamma_t\in L^\infty$, the boundary is locally a graph that is Lipschitz in time. 
Reynolds' transport formula on each graph follows from the Leibniz rule for an absolutely continuous integration limit. 
Thus $t\mapsto\int_{\Omega(t)}\varphi(t,y)\,\d y$ is absolutely continuous, and it follows from
\eqref{eq:transport-boundary-equation} and \eqref{eq:transport-normal-flux} that, we have for almost every $t$, 
\begin{equation*}
\begin{aligned}
\frac{\d}{\d t}\int_{\Omega(t)}\varphi(t,y)\,\d y
&=\int_{\Omega(t)}\p_t\varphi\,\d y + \int_{\p\Omega(t)}\varphi V\cdot n\,\d\sigma\\
&=\int_{\Omega(t)}(\p_t\varphi+v\cdot\nabla\varphi)\,\d y.
\end{aligned}
\end{equation*}
Integrating from $0$ to $T$ and using $\varphi|_{t=T}\equiv 0$ gives (i).
\end{proof}

\begin{remark}
The divergence-free hypothesis in \autoref{prop:weak-transport-boundary-equivalence} can be replaced by
$\div v\in L^1_{\loc}([0,T]\times\R^2)$. This may be useful in the study of compressible fluids.
\end{remark}

\subsection{The constant-speed parametrization}\label{sec:arclength-parametrization}

\begin{proposition}\label{prop:constant-speed-reparametrization}
Let $\Omega(t)$ and $V$ be as in the hypotheses and condition (i) of \autoref{prop:weak-transport-boundary-equivalence}. 
Assume additionally that the normalized arclength parametrization of $\Omega(t)$ is bounded in $H^2$ uniformly in $t$, 
and $V$ is space-time continuous on the moving boundary. 
Assume for almost every $t\in(0,T)$
\begin{equation}\label{eq:boundary-stretching-rate}
V(t,\cdot)\cdot\tau \in W^{1,1}(\p\Omega(t)),
\end{equation}
where $\tau$ is a unit tangent on $\p\Omega(t)$ and the derivative in $W^{1,1}$ is taken with respect to arclength.

Then, given any constant speed parametrization $\gamma_0$ of $\p\Omega(0)$, 
there exists $\gamma\in C_tC^1$ with $\gamma_t\in L^\infty((0,T)\times\T)$ such that $\gamma(0,\cdot)=\gamma_0$, and for each $t\in[0,T]$, $\gamma(t,\cdot)$ is a constant-speed parametrization of $\p\Omega(t)$.
Moreover, the map $\gamma$ satisfies
\begin{equation}\label{eq:transport-constant-speed-equation}
\gamma_t=V(t,\gamma)+\lambda_\gamma\gamma'
\end{equation}
almost everywhere on $(0,T)\times\T$, where
$$
\lambda_\gamma(t,x)=\frac{x+\pi}{2\pi}\int_\T g_\gamma(t,y)\,\d y -\int_{-\pi}^x g_\gamma(t,y)\,\d y, \qquad g_\gamma
=\frac{(\gamma'\cdot V(t,\gamma))'-\gamma''\cdot V(t,\gamma)}{|\gamma'|^2}.
$$
If the boundaries admit a regular parametrization in $C([0,T];H^2)$, then the constructed $\gamma$ belongs to $C([0,T];H^2)$.
\end{proposition}

\begin{proof}
\textbf{(I) Construction.}
Assume that we work with positive orientation, so that $\gamma'$ is tangent and $\n_\gamma$ is the outward normal.
Write $L(t) = |\p\Omega(t)|$ for the boundary length. Since we want a constant-speed parametrization $\gamma$, $\gamma(t)$ will be determined by a starting point, in our case $\gamma(t,-\pi)$. Since we want $\lambda_\gamma(t,-\pi)=0$, we will choose $\gamma(t,-\pi)$ to move with the velocity $V$ along the boundary. The rest of the parametrization is then determined by the constant-speed condition.

Let $\tilde{\gamma}(t,\cdot)$ be the parametrization of $\p\Omega(t)$ constructed in the proof of \autoref{prop:weak-transport-boundary-equivalence}. 
Since $\tilde{\gamma}(t,\T)=\p\Omega(t)$, we can write the starting point as $\gamma(t,-\pi)=\tilde{\gamma}(t,a(t))$ with $a(t)\in\T$. 
Now the task is to solve $a(t)$. Since this point moves with velocity $V$ along the boundary, we have
$$
V(t,\tilde{\gamma}(t,a(t))) =\frac{\d}{\d t}\tilde{\gamma}(t,a(t))=\tilde{\gamma}_t(t,a(t))+\tilde\gamma'(t,a(t))a'(t),
$$
hence $a(t)$ satisfies the scalar ODE
$$
a'(t)=\frac{\bigl[V(t,\tilde{\gamma}(t,a(t)))-\tilde{\gamma}_t(t,a(t))\bigr]
\cdot\tilde{\gamma}'(t,a(t))}{|\tilde{\gamma}'(t,a(t))|^2},
$$
with initial value satisfy $\tilde{\gamma}(0,a(0))=\gamma_0(-\pi)$.
Since for each $t$, $|\tilde{\gamma}'(t,a(t))|$ is bounded away from zero, the right-hand side is bounded for each fixed $t$. 
Applying Carath\'eodory's theorem, we obtain a Lipschitz solution $a(t)$ on $[0,T]$. Here $a(t)$ is understood as a continuous representative modulo $2\pi$.
After choosing $\gamma(t,-\pi)=\tilde{\gamma}(t,a(t))$, we can define $\gamma$ to a constant-speed parametrization of $\p\Omega(t)$.

\textbf{(II) The dynamics, and regularity of $\gamma$.}
Now we have constructed a constant-speed parametrization $\gamma$ of $\p\Omega(t)$.
Assume that $\gamma\in C_tC^1$, $\gamma_t\in L^\infty$. 
We have from \eqref{eq:transport-boundary-equation} that $\gamma_t=V(t,\gamma)+\lambda\gamma'$, where the scalar coefficient $\lambda$ is bounded, 
since $\gamma_t$ and $V$ are bounded and $|\gamma'|$ is bounded away from zero. 
Then the same calculation as in \autoref{sec:constant-speed-gauge} proves \eqref{eq:transport-constant-speed-equation}.

Now we only need to obtain the regularity of $\gamma$ and $\lambda_\gamma$.
The calculation in (II) requires $\gamma_t\in L^\infty$, $\gamma\in C_tC^1$, $g_\gamma(t,\cdot)\in L^1(\T)$ for almost every $t$. 
We verify them and finally prove the assertion about $C_tH^2$.

Denote, for $x\in \T$, $\gamma(t,x)=\tilde{\gamma}(t,a(t,x))$. Then $a(t,x)$ is determined by
$$
\int_{a(t,-\pi)}^{a(t,x)}|\tilde{\gamma}'(t,z)|\,\d z=\frac{x+\pi}{2\pi}L(t).
$$
For $\gamma_t\in L^\infty$ we need to show that $\gamma$ is uniformly Lipschitz in time. Directly calculation gives
$$
\begin{aligned}
&|\gamma(t,x)-\gamma(s,x)|=|\tilde{\gamma}(t,a(t,x))- \tilde{\gamma}(s,a(s,x))|\\ 
\le & |\tilde{\gamma}(t,a(t,x))-\tilde{\gamma}(t,a(s,x))|+|\tilde{\gamma}(t,a(s,x))-\tilde{\gamma}(s,a(s,x))|\\
\le & \|\tilde{\gamma}'\|_{L_{t,x}^\infty}|a(t,x)-a(s,x)|+\|\tilde{\gamma}_t\|_{L_{t,x}^\infty}|t-s|.
\end{aligned}
$$
Thus it suffices to show that $a(t,x)$ is Lipschitz in time.
Since $\tilde{\gamma}\in C_tC^1$ and $|\tilde{\gamma}'|>0$, we have $\inf_{t,z}|\tilde{\gamma}'(t,z)|>c$ for some $c>0$, and it suffices to show that 
$$
\int_{a(s,x)}^{a(t,x)} |\tilde{\gamma}'(t,z)|\,\mathrm{d}z \le C|t-s|, \qquad 0\leq t,s\leq T.
$$

We claim that the change in arc length between $s$ and $t$ satisfies
\begin{equation}\label{eq:appendix-arc-length-change}
\left|\int_{x_1}^{x_2}\bigl(|\tilde{\gamma}'(t,z)|-|\tilde{\gamma}'(s,z)|\bigr)\ \d z\right|\leq C|t-s|,
\qquad 0\leq x_2-x_1\leq 2\pi.
\end{equation}
This, however, doesn't immediately follow from $\tilde{\gamma}\in C_tC^1$. We need some extra regularity.
Write $\tau(t,z):=\tilde{\gamma}'(t,z)/|\tilde{\gamma}'(t,z)|$. 
Since
$$
|\tilde{\gamma}'(t,z)|-|\tilde{\gamma}'(s,z)| \leq \tau(t,z)\cdot \left(\tilde{\gamma}'(t,z)-\tilde{\gamma}'(s,z)\right),
$$
if we can bound $\|\tau'(t,\cdot)\|_{L^1(\T)}\le C$, then 
$$
\begin{aligned}
&\int_{x_1}^{x_2}\bigl(|\tilde{\gamma}'(t,z)|-|\tilde{\gamma}'(s,z)|\bigr)\,\d z\\
\leq &\int_{x_1}^{x_2}\tau(t,z)\cdot \bigl(\tilde{\gamma}'(t,z)-\tilde{\gamma}'(s,z)\bigr)\,\d z\\
=&\,\bigl[\tau(t,z)\cdot\bigl(\tilde{\gamma}(t,z)-\tilde{\gamma}(s,z)\bigr)\bigr]\Big|_{x_1}^{x_2} - \int_{x_1}^{x_2}\tau'(t,z)\cdot \bigl(\tilde{\gamma}(t,z)-\tilde{\gamma}(s,z)\bigr)\,\d z\\
\leq &\,\bigl(2+\|\tau'(t,\cdot)\|_{L^1(\T)}\bigr) \|\tilde{\gamma}(t,\cdot)-\tilde{\gamma}(s,\cdot)\|_{L_x^\infty}\\
\leq &\,C\|\tilde{\gamma}_t\|_{L_{t,x}^\infty}|t-s|.
\end{aligned}
$$
Interchanging $s$ and $t$ proves the claim. In particular, $x_2=x_1+2\pi$ gives $|L(t)-L(s)|\leq C|t-s|$.

It remains to show that $\|\tau'(t,\cdot)\|_{L^1(\T)}\leq C$ to bound the integral after integration by parts.
Differentiating the arc length formula gives $a'(t,x)|\tilde{\gamma}'(t,a(t,x))|=L(t)/2\pi$, hence $\inf_{x}a'(t,x)>0$, so $a(t,\cdot)$ admits a $C^1$ inverse on $\T$. 
By $\tau(t,a(t,x))=(2\pi/L(t))\gamma'(t,x)$, we obtain $\tau(t,\cdot)$ by composing $(2\pi/L(t))\gamma'(t,\cdot)\in H^1(\T)$ with the $C^1$ inverse of $a(t,\cdot)$. The chain rule now gives $\tau(t,\cdot)\in H^1(\T)$ and $\tau'(t,a(t,x))a'(t,x)=(2\pi/L(t))\gamma''(t,x)$. 
Changing variables by $z=a(t,x)$ implies
$$
\begin{aligned}
\int_\T|\tau'(t,z)|\,\d z
&=\int_\T|\tau'(t,a(t,x))|a'(t,x)\,\d x\\
&=\frac{2\pi}{L(t)}\int_\T|\gamma''(t,x)|\,\d x\leq C,
\end{aligned}
$$
where $L(t)\geq 2\pi c$ and the $H^2$ bound on $\gamma$.

With \eqref{eq:appendix-arc-length-change} in hand, we can now bound the integral between $a(s,x)$ and $a(t,x)$.
Indeed we have 
$$
\begin{aligned}
\int_{a(s,x)}^{a(t,x)} |\tilde{\gamma}'(t,z)|\,\d z
&=\frac{x+\pi}{2\pi}\bigl(L(t)-L(s)\bigr)\\
&\quad+\int_{a(s,-\pi)}^{a(s,x)}
   \bigl(|\tilde{\gamma}'(s,z)|-|\tilde{\gamma}'(t,z)|\bigr)\,\d z\\
&\quad+\int_{a(s,-\pi)}^{a(t,-\pi)}|\tilde{\gamma}'(t,z)|\,\d z\\
& \le C|t-s|.
\end{aligned}
$$
The lower bound for $|\tilde{\gamma}'|$ now gives $|a(t,x)-a(s,x)|\leq C|t-s|$, and hence $\gamma_t\in L^\infty$.

Now we prove $\gamma\in C_tC^1$. Indeed, as $t\to t_0$, we have $\gamma(t,\cdot)\to\gamma(t_0,\cdot)$ in $L^\infty$, and $\left\|\gamma(t) \right\|_{H^2}\leq C$. 
Interpolation gives $\gamma(t,\cdot)\to\gamma(t_0,\cdot)$ in $C^1$.

\textbf{(III) Regularity of $\lambda_\gamma$ and $g_\gamma$.}
Now we can apply \autoref{prop:weak-transport-boundary-equivalence} to write $\gamma_t=V(t,\gamma)+\lambda\gamma'$ for some bounded $\lambda$.
We justify the calculation in \autoref{sec:constant-speed-gauge} for $\lambda_\gamma$ and $g_\gamma$.

Again we calculate the following
\footnote{This can be verfied by mollification.} 
$$
\frac12\p_t|\gamma'|^2 = \gamma'\cdot\gamma_t' =(\gamma'\cdot V(t,\gamma))'-\gamma''\cdot V(t,\gamma)+|\gamma'|^2\lambda',
$$
Using $|\gamma'|=L/(2\pi)$, we have
$$
\frac{L_t}{L}=\lambda'+g_\gamma
$$
for almost every $t$. Here we denote $g_\gamma = \frac{(\gamma'\cdot V(t,\gamma))'-\gamma''\cdot V(t,\gamma)}{|\gamma'|^2}$.
Averaging this in $x$ gives
\footnote{
We need $g_\gamma\in L^1(\T)$ to rigorously take the averaging.
To show $(\gamma'\cdot V(t,\gamma))'\in L^1(\T)$, write
$$
\gamma'(t,x)\cdot V(t,\gamma(t,x)) =\frac{L(t)}{2\pi}(V(t,\cdot)\cdot\tau)(\gamma(t,x))
$$
and apply \eqref{eq:boundary-stretching-rate}.
The product $\gamma''\cdot V(t,\gamma)\in L^1(\T)$ since $\gamma''\in L^2(\T)$ and $V$ is bounded. 
This concludes $g_\gamma(t,\cdot)\in L^1(\T)$ since $|\gamma'|$ is bounded away from zero.
}
$$
\lambda'+g_\gamma=\frac{1}{2\pi}\int_\T g_\gamma(t,y)\,\d y,
$$
and integrating in $x$ gives the formula for $\lambda_\gamma$:
$$
\lambda_\gamma(t,x)=\frac{x+\pi}{2\pi}\int_\T g_\gamma(t,y)\,\d y - \int_{-\pi}^x g_\gamma(t,y)\,\d y + \lambda_\gamma(t,-\pi).
$$
This clearly gives $\lambda_\gamma$ is continuous.
To settle the constant $\lambda_\gamma(t,-\pi)$, we back to $\gamma_t=V(t,\gamma)+\lambda_\gamma\gamma'$. 
Though this holds only almost everywhere, now we have already shown that the right-hand side is continuous in $x$, hence $\gamma_t$ is continuous in $x$. 
Taking $x=-\pi$ gives $\gamma_t(t,-\pi)=V(t,\gamma(t,-\pi))+\lambda_\gamma(t,-\pi)\gamma'(t,-\pi)$. As in the construction, $\gamma(t,-\pi)$ moves with velocity $V$, hence $\lambda_\gamma(t,-\pi)=0$ since $\gamma'(t,-\pi)\neq 0$. This finally justify \eqref{eq:transport-constant-speed-equation}.

\textbf{(IV) Regularity of $\gamma$ in $C_tH^2$.}
It remains to prove continuity of $\gamma''$ in $L^2$, provided that $\eta$ is a regular parametrization of $\p\Omega(t)$ in $C_tH^2$.
Fix $t$ and let $\gamma(t,x)=\eta(t,y)$. The curvature at this point is $(\eta'/|\eta'|)'(t,y)/|\eta'(t,y)|$, while for $\gamma$ this is simply $\gamma''(t,x)/|\gamma'(t,x)|^2$. Hence we have
$$
\gamma''(t,x)=\left(\frac{L(t)}{2\pi}\right)^2
\frac{(\eta'/|\eta'|)'(t,y)}{|\eta'(t,y)|}.
$$
The arc length element gives $\d x=(2\pi/L(t))|\eta'(t,y)|\,\d y$, so
$$
\begin{aligned}
\|\gamma''(t,\cdot)\|_{L^2(\T)}^2
&=\left(\frac{L(t)}{2\pi}\right)^4\frac{2\pi}{L(t)} \int_\T\left|\frac{(\eta'/|\eta'|)'(t,y)}{|\eta'(t,y)|}\right|^2 |\eta'(t,y)|\,\d y\\
&=\left(\frac{L(t)}{2\pi}\right)^3 \int_\T\frac{|(\eta'/|\eta'|)'(t,y)|^2}{|\eta'(t,y)|}\,\d y.
\end{aligned}
$$
The norm is continuous in $t$: Since $L$ is continuous, $(\eta'/|\eta'|)'$ is in $L^2$, and $|\eta'|^{-1}$ is in $C^0$. The last two follow from the chain rule and embedding $H^1(\T)\hookrightarrow C^0(\T)$, using $\eta\in C_tH^2$ and $\inf_{t,y}|\eta'(t,y)|>0$. 
The $C^1$ continuity of $\gamma$ and the uniform $H^2$ bound give $\gamma''(t)\rightharpoonup\gamma''(t_0)$ weakly in $L^2$ as $t\to t_0$, so convergence of the norms implies strong convergence.
\end{proof}

\subsection{Application to gSQG and SQG}\label{sec:patch-contour-application}

\begin{proof}[Proof of \autoref{prop:patch-constant-speed-equivalence}]
We obtain the contour dynamic equation \eqref{eq:contour-equation}-\eqref{eq:lambda-definition}
by verifying the hypotheses of \autoref{prop:weak-transport-boundary-equivalence} and \autoref{prop:constant-speed-reparametrization}. 
The reverse part follows directly from \autoref{prop:weak-transport-boundary-equivalence}.

Let's start with a patch solution $\Omega(t)$ in the sense of \autoref{def:eulerian-patch}.
The $C^1$ compactness of the boundaries follows from the uniform $H^2$ bound and the arc--chord condition. Moreover, the limit curves are simple closed.
Now, we take care of the velocity $u(t)=\nabla^\perp(|\cdot|^{-\beta}*\ch_{\Omega(t)})$. 
Clearly it's divergence-free, and locally integrable in space-time since it's uniformly bounded in BMO. 
So it remains to verify the flux condition \eqref{eq:transport-normal-flux}. 

Denote $\t$ the positively oriented unit tangent,
integration by parts gives, in distributions,
$$
u(t,y)=-\int_{\p\Omega(t)}|y-z|^{-\beta}\t(t,z)\,\d\sigma(z).
$$

For each fixed $t$, we choose a positively oriented constant-speed parametrization $\gamma(t,\cdot)\in H^2(\T;\R^2)$ of $\p\Omega(t)$, and define
$$
V(t,\gamma(t,x))=\Nop(\gamma)(t,x)
=\int_\T\frac{\Delta_h\gamma'(t,x)}{|\Delta_h\gamma(t,x)|^\beta}\,\d h.
$$
Writing $n$ for the outward unit normal, we have for every
$\varphi\in C_c^\infty(\R^2)$
\begin{equation}\label{eq:patch-spatial-flux}
\int_{\Omega(t)}u\cdot\nabla\varphi\,\d y
=\int_{\p\Omega(t)}\varphi\,V\cdot n\,\d\sigma.
\end{equation}
For $0<\beta<1$, the Gauss--Green formula applies since $u$ is continuous. For $\beta=1$, the identity follows from the spatial integration by parts for $J$ in the proof of \cite[Proposition~3.2]{Gancedo2008}, with the normal trace given by (11)--(12) of that paper. 
Anyway in both cases, $\gamma'(t,x)\cdot\n_\gamma(t,x)=0$ allows $-\gamma'(t,x-h)$ to be replaced by $\Delta_h\gamma'(t,x)$ in the normal component.

Note that, combining \eqref{eq:arc-chord-lower} with \eqref{eq:basic-difference-local} for $f=\gamma'$ gives
\begin{equation}\label{eq:contour-operator-sup}
\left\|\Nop(\gamma)\right\|_{L^\infty(\T)} \leq F(\gamma)^\beta\left\|\gamma''\right\|_{L^2(\T)} \int_\T |h|^{1/2-\beta}\,\d h.
\end{equation}
we then have a uniform bound on $V$. 
Now we can apply \autoref{prop:weak-transport-boundary-equivalence} to obtain a $C_tC^1$ parametrization of the moving boundary, whose derivative in time is in $L^\infty((0,T)\times\T)$, and \eqref{eq:transport-boundary-equation} holds. 
The continuity of $V$ on moving boundaries then follows directly from \eqref{eq:contour-operator-sup} and dominated convergence.

Finally we give $g_\gamma$ needed in \autoref{prop:constant-speed-reparametrization}.
By \eqref{eq:g-single-F} for $0<\beta<1$ and \eqref{eq:sqg-gauge-derivative} for $\beta=1$, we have that $g_\gamma(t,\cdot)\in L^1(\T)$. Since $\gamma$ is constant-speed, it's easy to see that condition \eqref{eq:boundary-stretching-rate} is satisfied.
Finally we can apply \autoref{prop:constant-speed-reparametrization} to obtain a constant-speed parametrization of the moving boundary, 
and \eqref{eq:contour-equation}-\eqref{eq:lambda-definition} holds. 

The reverse part follows directly from \autoref{prop:weak-transport-boundary-equivalence}\footnote{\ 
With \eqref{eq:patch-spatial-flux}, it remains to verify $\gamma_t\in
L^\infty_{t,x}$, $\gamma\in C_tC^1$, $\ch_{\Omega(t)}\in C_wL^2$, and
\eqref{eq:transport-boundary-equation}. The contour equation gives
$\gamma_t=\Nop(\gamma)+\lambda_\gamma\gamma'\in L^\infty_{t,x}$.
Interpolating the resulting time Lipschitz bound with the uniform $H^2$
bound gives $\gamma\in C_tC^1$, hence $\ch_{\Omega(t)}\in C_tL^2$.
\eqref{eq:transport-boundary-equation} follows directly from taking normal component of the contour equation.
}. Also, in the case $\beta=1$, \autoref{prop:constant-speed-reparametrization} gives the desired regularity $\gamma \in C([0,T];H^2)$.
This completes the proof of \autoref{prop:patch-constant-speed-equivalence}.
\end{proof}

\noindent\textbf{Acknowledgements.}
This research is supported by the National Key Research and Development Program of China (No.2024YFA1013301) and NSFC (No. 12071211).

\noindent\textbf{AI Declaration.}
GPT-6.2 Sol is used to assist with language and organization, to solve questions about \LaTeX, and to draw TikZ figures like \autoref{fig:boundary-coordinates}.
The authors also used these tools to check the proof, improve the clarity of the exposition. The full article was written by the authors.
The main ideas are developed by the authors, who take full responsibility for the mathematical content.

\bibliographystyle{amsplain}
\bibliography{main}

\end{document}